\documentclass[10pt]{amsart}
\usepackage{amsmath,amsthm,amssymb,amsfonts}
\usepackage{hyperref}
\usepackage[nobysame,initials]{amsrefs}
\newcommand{\IGNORE}[1] {}
\newcommand{\eps}{\varepsilon}

\def\ov#1{\overline{#1}}
\newcommand{\refs}{\eqref}
\newcommand{\barmu}{\bar{\mu}}
\newcommand{\ra}{\rightarrow}

\newcommand{\res}{\mathrm{res}}
\newcommand{\Pic}{\mathrm{Pic}}

\newcommand{\bash}{\backslash}

\newcommand{\mbfx}{{{x}}}
\newcommand{\mbfy}{{{y}}}
\newcommand{\Cc}{\mathbb{C}}
\newcommand{\Zz}{\mathbb{Z}}

\newcommand{\Qq}{\mathbb{Q}}
\newcommand{\Rr}{\mathbb{R}}
\newcommand{\PSL}{\mathrm{PSL}}

\newcommand{\x}{{\mathbf{x}}} 

\newcommand{\imax}{i_{max}}

\newcommand{\imin}{i_{min}}

\newcommand{\question}{\beta}
\newcommand{\algebra}{A}
\newcommand{\tr}{\mathrm{tr}}

\newcommand{\disc}{\mathrm{disc}}
\newcommand{\reg}{\mathrm{reg}}

\newcommand{\epsilonvector}{\underline{\epsilon}}

\renewcommand{\v}{\mathbf{v}}
\newcommand{\B}{\mathcal{B}}
\newcommand{\bB}{\bar{\mathcal{B}}}

\newcommand{\Gal}{\mathrm{Gal}}

\renewcommand{\t}{\mathfrak{t}}
\newcommand{\T}{\mathbf{T}}

\newcommand{\Ad}{\mathrm{Ad}}

\newcommand{\h}{\mathfrak{h}}

\newcommand{\data}{\mathscr{D}}

\newcommand{\g}{\mathfrak{g}}
\newcommand{\C}{\mathbb{C}}

\newcommand{\vol}{\mathrm{vol}}

\newcommand{\G}{\mathbf{G}}
\newcommand{\Char}{\Omega(C_F)}

\newcommand{\CharK}{\Omega(C_K)}

\newcommand{\Cond}{\mathrm{Cond}}

\newcommand{\height}{\mathrm{ht}}

\newcommand{\GL}{\mathrm{GL}}

\newcommand{\PGL}{\mathrm{PGL}}
\newcommand{\pgl}{\mathfrak{pgl}}
\newcommand{\PO}{\mathrm{PO}}

\newcommand{\SL}{\mathrm{SL}}

\newcommand{\order}{\mathscr{O}}
\newcommand{\N}{\mathbb{N}}
\newcommand{\R}{\mathbb{R}}

\newcommand{\Z}{\mathbb{Z}}
\newcommand{\Q}{\mathbb{Q}}

\newcommand{\adele}{\mathbb{A}}
\newcommand{\adelef}{\mathbb{A}_{f}}
\newcommand{\Norm}{\mathrm{N}}

\newcommand{\dist}{\mathrm{dist}}

\newcommand{\Lie}{\operatorname{Lie}}
\DeclareFontFamily{OT1}{rsfs}{}
\DeclareFontShape{OT1}{rsfs}{n}{it}{<-> rsfs10}{}
\DeclareMathAlphabet{\mathscr}{OT1}{rsfs}{n}{it}

\begin{document}
\swapnumbers

\newtheorem{Theorem}[subsection]{Theorem}
\newtheorem*{Theorem*}{Theorem}
\newtheorem{Example}[subsection]{Example}
\newtheorem{Lemma}[subsection]{Lemma}
\newtheorem{Definition}[subsection]{Definition}
\newtheorem{Proposition}[subsection]{Proposition}
\newtheorem{Corollary}[subsection]{Corollary}
\newtheorem{Conjecture}[subsection]{Conjecture}
\newtheorem{Hypothesis}[subsection]{Hypothesis}
\theoremstyle{definition}
\newtheorem*{note}{Note}
\newtheorem*{Remark}{Remark}

\title[Duke's theorem for cubic fields.]{
Distribution of periodic torus orbits
and Duke's theorem for cubic fields.}
\author[M. Einsiedler]{Manfred Einsiedler}
\address[M. E.]{ETH Z\"urich, R\"amistrasse 101
CH-8092 Z\"urich
Switzerland}
\email{manfred.einsiedler@math.ethz.ch}
\author[E. Lindenstrauss ]{Elon Lindenstrauss}
\address[E. L.]{Department of Mathematics, Fine Hall, Washington Road, Princeton University 
Princeton, NJ 08544, USA \&
The Einstein Institute of Mathematics\\
Edmond J. Safra Campus, Givat Ram, The Hebrew University of Jerusalem
Jerusalem, 91904, Israel}
\email{elonl@math.princeton.edu \& elon@math.huji.ac.il}
\author[Ph. Michel]{Philippe Michel}
\address[Ph. M.]{EPF Lausanne, SB-IMB-TAN, Station 8, CH-1015 Lausanne, Switzerland}
\email{philippe.michel@epfl.ch}
\author[A. Venkatesh]{Akshay Venkatesh}
\address[A. V.]{Department of Mathematics, building 380, Stanford, CA 94305, USA}
\email{akshay@math.stanford.edu}
%\subjclass[2000]{Primary 37A17; Secondary 11F70}
\maketitle

\begin{abstract}
We study periodic torus orbits on spaces of lattices. Using the action
of the group of adelic points of the underlying tori,
 we define a natural equivalence relation on these orbits,
and show that the equivalence classes become uniformly distributed.
This is a cubic analogue of Duke's theorem about the distribution of closed geodesics on the modular surface: 
suitably interpreted, the ideal classes of a cubic totally real field are equidistributed in the modular $5$-fold
$\SL_3(\Z) \backslash \SL_3(\R) / \mathrm{SO}_3$.
In particular, this proves (a stronger form of) the folklore conjecture that the collection of maximal compact flats
in $\SL_3(\Z) \backslash \SL_3(\R) / \mathrm{SO}_3$ of volume $\leq V$ becomes equidistributed
as $V \rightarrow \infty$.

The proof combines subconvexity estimates, measure classification, and local harmonic analysis.
\end{abstract}

\tableofcontents

\section{Introduction}

\subsection{Historical perspective}
In the preface to his book
 ``Ergodic properties of algebraic fields,''  Linnik \cite{Linnik} writes
\begin{quote} \it
\dots In the present book other applications of the ergodic concepts are presented.
 Constructing ``flows'' of integral points on certain algebraic manifolds given by systems of integral polynomials,
  we are able to prove individual ergodic theorems and mixing theorems in certain cases. These theorems permit
   asymptotic calculations of the distribution of integral points on such manifolds and we arrive at results inaccessible
    up to now by the usual methods of analytic number theory. Typical in this respect is this theorem concerning the
     asymptotic distribution and ergodic behavior of the set of integral points on the sphere
\begin{equation*} \label{points on sphere}
x ^2 + y ^2 + z ^2 = m \tag{$*$}
\end{equation*}
for increasing $m$.
\end{quote}

\medskip

This presents what Linnik called ``the ergodic method'';
it enabled Linnik to show that solutions to ($*$) become equidistributed upon projection to the unit sphere -- 
at least, for the $m$ satisfying an explicit congruence condition. Subsequently,
using that method, Skubenko solved the related problem for the solutions of the equation 
 \begin{equation*} \label{hyperboloic}
y^2-xz = m \tag{$**$}
\end{equation*} also under similar congruence conditions on $m$. 

Both of these problems are related to the distribution of ideal classes of orders in quadratic fields: in the case of points on the sphere ($*$), one deals with imaginary quadratic fields, while ($**$) corresponds either to real or imaginary quadratic fields depending on the sign of $m$.
  The latter problem for $m>0$ is also equivalent  to the problem of the distribution of closed geodesics on the modular surface $\SL (2, \Z) \backslash \mathbb{H}$.
 
Since the time of Linnik's work, the tools of analytic number theory have developed tremendously. 
In particular W. Duke \cite{Duke}, using a breakthrough of H. Iwaniec, 
proved that the integer solutions of ($*$) as well as ($**$) become equidistributed as
 $m \to \infty$. 

In \cite [Chp. VI-VII]{Linnik}, Linnik considers in detail the 
corresponding questions for number fields of higher degree,  particularly cubic fields. However, he was able to prove
comparatively little compared to the quadratic case.
In modern terms, he established,
by a remarkable elementary calculation, a special case of the equidistribution of Hecke points.

In this paper, we revisit Linnik's problems for cubic (and higher degree) fields, settling them for totally real cubic fields (among other cases).

We give precise statements later in the introduction; for now, 
we note that
just as $(**)$ relates to the ideal class group of quadratic fields
and to closed geodesics on the modular surface, 
the higher rank analogues will pertain to ideal class groups of higher
degree fields and to periodic orbits of maximal tori in the space of lattices $X_{n}=\PGL_{n}(\Z)\bash\PGL_{n}(\Rr)$. We will, in fact,
explain our main result in this language in the Introduction;
there is also an interpretation analogous to $(**)$, which we postpone to
Corollary \ref{IPCor} in \S \ref{NT}.

Interestingly, we do not know how to prove our main equidistribution result using purely analytic techniques, nor using purely ergodic theoretic techniques, though each of these methods does give some partial information in this direction.
Our proof works by combining these two very different techniques;  to handle the case of non-maximal orders we also have to prove new estimates involving local Fourier analysis.\footnote{It is also conceivable that an ergodic
approach to these estimates may exist.}

This paper is part of a series of papers we have been writing on the distribution properties of compact torus orbits on homogeneous spaces. In \cite{ELMV1},
 we present a general setup for the study of the periodic orbits and prove results regarding the distribution of individual orbits as well as fairly arbitrary collections of periodic orbits.
  In \cite{ELMV2}, we give a modern reincarnation of Linnik's original argument, giving in particular a purely dynamical proof of equidistribution in the problem $(**)_{m>0}$;
   without an auxiliary  congruence condition. We still do not know how to give a purely dynamic proof of Duke's theorems regarding the equidistribution of the solutions to ($*$)
  or to $(**)_{m<0}$. Each of these papers is self-contained and can be read independently;
     related discussions can also be found in \cites{MV-ICM, EL-ICM}.

\subsection{Geometric perspective} \label{sec:GI}
For clarity, we shall continue to focus on the ``$\R$-split case'' of our main questions
(i.e.\  problem $(**)_{m>0}$, totally real fields, orbits of $\R$-split tori, etc.).  We introduce our main result in geometric terms. 
Later, in \S \ref{NT}, we discuss interpreting it in ``arithmetic'' terms (akin to the interpretations
$(*)$ and $(**)$). 

Let $M = \Gamma \backslash \mathbb{H} $ be a compact hyperbolic Riemann surface, and $$X := S^{1} M = \Gamma \backslash \mathrm{PSL}_2(\R)$$ the unit tangent bundle of $M$.
Bowen and Margulis, independently \cites{Bowen,Margulis-thesis}, proved that the set of geodesics of length $\leq L$, considered as closed
 orbits of the geodesic flow on $X$, are equidistributed w.r.t.\ Liouville measure
as $L\ra+\infty$. 

On the other hand, for a Riemann surface, if the latter is ``arithmetic" in a suitable sense (see \S \ref{sec:quaternion} below)
 the results of Bowen/Margulis are valid in a much stronger form.
 The most basic instance of an arithmetic surface
is the modular surface, $$M = M_2 = \PSL_2(\Z) \backslash \mathbb{H}.$$
The equidistribution theorem of Duke already mentioned implies\footnote{Strictly speaking, 
Duke's theorem establishes equidistribution on 
the modular surface; equidistribution at the level of the unit tangent bundle
was established in the unpublished PhD thesis of R. Chelluri, \cite{Chelluri}.}
that the collection of geodesics {\em of fixed length $\ell$} becomes equidistributed in $X = X_2$
as $\ell \rightarrow \infty$.
Note that the lengths of closed geodesics can have high multiplicity; indeed, the lengths  
are of the form $\log(d + \sqrt{d^2-1})$, for  $d \in \N_{>0}$, and the set of geodesics of this length
is parameterized by the class group of the quadratic order of discriminant $d$. 

Our main result establishes the analogues of both of these equidistribution theorems to
the setting of the rank two Riemannian manifold
 $$M_3 = \mathrm{PGL}_3(\Z) \backslash \mathrm{PGL}_3(\R)/\mathrm{PO}_3.$$
The role of ``closed geodesics'' is replaced by ``maximal compact flats'', and the role of ``quadratic order'' is replaced by ``cubic order''.

At least with our present 
understanding, the rank two case seems to be much more difficult than the rank one case.
This difficulty manifests itself from all the perspectives; crudely,
the smaller the {\em acting} group, relative to the {\em ambient} group,
the more difficult the question.
\subsection{Statement of results.}
Now let us give a precise statement of our main result,
at least in the ``$\R$-split, cubic field'' case. (Our most general
theorem is stated in Theorem \ref{RealTheorem}; what follows is a specialization of this.)

Let $H$ be the diagonal subgroup of $\PGL_3(\R)$.
In \cite{ELMV1}, we attached to each
closed orbit $x H$ on $$X_3 =  \PGL_{3}(\Z) \backslash \PGL_{3}(\R)$$ a discriminant $\disc(x H)$ as a way to measure the 
``arithmetic complexity" of that orbit; let us briefly recall how it is obtained.
 Writing $x H$ into the form $\Gamma\bash \Gamma gH$, we set $\mathbf{T}=\overline{\Gamma\cap gHg^{-1}}^{\mathrm{Zar}}$;
  $\mathbf{T}$ is a maximal (anisotropic) $\Qq$-torus. The discriminant is then closely related to the ``denominator" of the $\Qq$-point $\mathbf{T}$ inside the variety of maximal tori of
  $\PGL_{3}$. In Section \ref{Hsets}, we will review this construction in the adelic setting.
  We prove:

\begin{Theorem}\label{mainthm1}
The periodic orbits of $H$ on $X_3$ are grouped into equivalence classes,
equivalent orbits having the same volume and discriminant.
For each periodic orbit $x H$, let $Y_{x H}$ be the union all  compact orbits equivalent to $x H$. 

If $\{x_{i} H\}_{i}$ is a sequence of compact orbits, with $\disc(x_{i}H)\ra+\infty$ then:
\begin{enumerate}
\item {{$\displaystyle{\vol(Y_{x_{i}H})}={ \disc(x_{i}H)^{1/2+o(1)}},$}}
\item the $Y_{x_{i}H}$ become uniformly distributed in $X_3$.
\end{enumerate}
In particular, for $V>0$, let $Y(V)$ denotes the collection of all $H$-orbits with volume $V$; as $V_i \rightarrow \infty$ through any sequence for which $Y(V_i)
\neq \emptyset$, then the $Y(V_i)$ become uniformly distributed in $X_3$.
\end{Theorem}

Noting that the projection of $Y(V)$ to $M_{3}$ is the collection of all maximal compact flats on $M_{3}$ of volume $V$, Theorem \ref{mainthm1}
 implies (a stronger form of) the rank two analogue of the Bowen/Margulis theorem, indicated above.

\subsection{About the proof; adelization.}

Viewed from the classical point of view, the grouping of periodic torus orbits into packets is rather mysterious and can be quite tricky to define in the non-maximal case. It turns out that the adeles give a powerful and concise language to describe these equidistribution results, and we have written the bulk of this paper consistently in the adelic language.

In particular, as we shall see in \S\ref{data} the full equivalence class $Y_{x_{i}H}$ of a periodic $H$-orbit ($H$ being a maximal split torus) is essentially the projection to the infinite place of a single periodic orbit of an \emph{adelic} torus. The precise connection between packets and adelic tori
is contained in Theorem \ref{thmpacket};
Theorem \ref{mainthm1} is then immediate from the adelic results
Theorem \ref{Volumes} and Theorem \ref{RealTheorem}.

We will, indeed, go to some length to set equidistribution questions
in a genuinely adelic framework, which has the pleasing side effect that we are able to address simultaneously many different equidistribution questions (cf. \S \ref{Hsets}).

To aid the reader, we give an outline of the main ideas that enter into its proof
in purely classical terms \S \ref{outline}.
For the moment, we only observe an important
contrast between $X_3$ and $X_2$: while for $X_2$, the analogue of our main Theorem is ``purely'' a result about $L$-functions, for $X_3$ this is not so.
To fill this gap, we will need to combine results from measure rigidity,
$L$-functions, and harmonic analysis on Lie groups.

\subsection{Scope of the method}
We shall discuss certain natural generalizations of Theorem \ref{mainthm1}
and interesting questions associated to them.

\subsubsection{$S$-arithmetic variants.}
Theorem \ref{mainthm1} is derived from the underlying adelic result -- Theorem \ref{RealTheorem} -- and there is therefore no difficulty
in replacing $\PGL_3(\Z)$ by a congruence subgroup, or $\Q$ by a number field,
or passing to an $S$-arithmetic context.

However, although this is not apparent from the statement
of Theorem \ref{mainthm1}, the general statement Theorem \ref{RealTheorem}
is not as satisfactory as the corresponding statement for $\PGL_2$ (given
in Theorem \ref{CU}). Indeed, our general $\PGL_3$-theorem imposes local
conditions, akin to a Linnik-type condition, which happen to be automatically
satisfied in the setting
of Theorem \ref{mainthm1}.

To dispense with these local conditions seems to be a very interesting and fundamental question.

\subsubsection{Sparse equidistribution} \label{sparseequid}
Assuming a suitable  subconvex estimate on $L$-func\-tions\footnote {The specific case of subconvexity needed is subconvexity in the level aspect, on $\GL(3)$.} one can obtain the following ``sparse equidistribution'' result by our methods:

{\em Notations as in Theorem \ref{mainthm1}, there is $\alpha < 1/2$ so that, if $\vol(x_{i}H) >\disc(x_{i}H)^{\alpha}$ then the $x_i H$
become uniformly distributed on $X_3$.}

 Presumably assuming the full force of GRH would yield $\alpha = 1/4$, although
to prove this requires more careful
local analysis than we have done. 

Conjecturally, however, a much stronger statement should hold: we conjecture (\cite{ELMV1})
 that this equidistribution statement for single $H$-orbits remains true for {\em any} $\alpha > 0$. 
We refer to \cite{ELMV1} for discussion of this conjecture, some partial results towards it, and counterexamples to more optimistic conjectures. 

\subsubsection{Spaces of higher dimensional lattices}
Much of our analysis carries through
from $X_3$ to 
$$X_{n}=\PGL_{n}(\Zz)\bash\PGL_{n}(\Rr).$$
There are two obstacles, however, to obtaining a complete
generalization of Theorem \ref{mainthm1}:
\begin{enumerate}
\item  the lack of available subconvex bounds,
\item the lack of suitable technology to rule
out ``intermediate limit measures.'' \end{enumerate}

At the moment we have little to offer concerning the second point. In the case when $n$ is prime, the issue of intermediate measures does not occur; if
we suppose Hypothesis \ref{subconvex} (i.e., a subconvex bound for Dedekind $\zeta$-functions of degree $n$ number fields),  then {\em the analog of Theorem \ref{mainthm1} holds,} i.e. packets of periodic torus orbits become equidistributed on $X_n$.

\subsubsection{An almost-subconvexity bound for class group $L$-functions of cubic fields}
Let $K$ be a real cubic field, and let $\psi$ be a nontrivial character of the class group of $K$.  

We have an associated $L$-function $$L(K, \psi,s) = 
\sum_{\mathfrak{a}} \psi(\mathfrak{a}) \Norm_{K/\Q} (\mathfrak{a})^{-s},$$
the sum being extended over all integral ideals $\mathfrak{a}$, and $\Norm_{K/\Q}(\mathfrak{a})$ denoting 
the norm of $\mathfrak{a}$. This $L$-function
has conductor $D_K$, the discriminant of $K$. 

One corollary to our main result is that, for any fixed $\delta > 0$, 
\begin{equation} \label{c2}\frac{\left| \sum_{\Norm_{K/\Q}(\mathfrak{a})  \leq \delta \sqrt{D_K} } \psi(\mathfrak{a}) \right| }{
\left| \sum_{\Norm_{K/\Q}(\mathfrak{a}) \leq \delta \sqrt{D_K}} 1 \right| }= o(1).\end{equation}

To put this in perspective, a subconvex bound for the degree $3$ $L$-function $L(K, \psi, s)$
would guarantee that the same is true if we replace $\delta \sqrt{D_K}$ by $D_K^{0.499}$.
The result \eqref{c2} could therefore be considered as a ``non-quantitative'' form of subconvexity for this degree $3$ $L$-function. 

Since this paper was submitted, 
 K. Soundararajan \cite{Sound} proved a very general weak subconvex bound, valid for a wide class of $L$-functions. His result
does not imply \eqref{c2} with currently known bounds: what is needed is any improvement of Stark's result \cite{Stark}:
$$\res_{s=1}\zeta_{K}(s)\gg \log(\disc(K))^{-1}.$$
e.g.  any larger exponent would suffice.

\subsubsection{The cocompact case}
If one considers the quotient of $\PGL_{n}(\Rr)$ 
by a lattice associated to a $\Rr$-split division algebra,
one obtains a {\em compact} quotient. One certainly believes
the analogue of Theorem \ref{mainthm1} to be valid,  but in this case
the methods of the present paper which use Eisenstein series in an essential way do not apply.

This is an instance in which the cocompact case seems harder that the
non-compact case. We refer to \cite{ELMV1} where we obtain (weaker) results in the cocompact case by different methods.

\subsection{Connection to existing work.}
In the rank one case of $\PGL_2$, the analogs
of the questions we consider have been intensively studied
from many perspectives, both from the perspective of the work of Linnik
and that of Iwaniec and Duke.

Concerning $\PGL_n$ for $n \geq 3$,
the direct ancestor of our work is that of Linnik,
who devotes several chapters of his book \cite{Linnik}
to the question of distribution of the packets $Y_{x_{i} H}$.
The paper of Oh and Benoist \cite{Benoist-Oh} considers problems similar to those we consider.
Both \cite{Linnik} and \cite{Benoist-Oh} give results about the problem
in the special case when the $\Q$-torus attached to $x_iH$
remains constant.

\subsection{Organization of the paper.}

In \S \ref{outline}, we present an outline of the proof of Theorem \ref{mainthm1} in entirely classical language.

In \S \ref{NT}, we discuss some of the arithmetic manifestations of our result.

In \S \ref{Hsets}, we present a systematic framework for thinking
about adelic equidistribution problems.
We then explain, in this context, our main results: Theorem \ref{Volumes}
and Theorem \ref{RealTheorem}. These imply the first and second assertions of Theorem \ref{mainthm1}.

In \S \ref{data}, we explain the grouping of periodic orbits into packets.
This uses the setup of \S \ref{Hsets}.

In \S \ref{HsetsGLn}, we specialize to the case of the group $\G = \PGL_n$
and explain the parameterization of packets of periodic orbits.

In \S \ref{local}, we give a brief recollection of properties of the building
of $\PGL_n$, over a local field.

In \S \ref{toruslocal}, we explain the local harmonic analysis
that will be needed.

In \S \ref{Notation}, we set up notational conventions about number fields, ideles
and adeles (especially: normalizations of measures).

In \S \ref{EisSec}, we set up general notation about Eisenstein series
(these are the generalization of the functions ``$E_f$'' discussed in \S \ref{outline}).

In \S \ref{core}, we prove, in adelic language,
the estimates for the integral of an Eisenstein series over a torus orbit
 (this is the adelic version of \eqref{siegel} from \S \ref{outline}). 

In \S \ref{explication},  we translate the results of \S \ref{core}
from the adelic to the $S$-arithmetic context, obtaining
Proposition \ref{siegelgen} (this is the $S$-arithmetic form of \eqref{siegel}
from \S \ref{outline}). 

In \S \ref{Proofs}, we complete the proofs of Proposition \ref{Volumes}
and Theorem \ref{RealTheorem}, and therefore also of Theorem \ref{cubicequid}.

In \S \ref{subconvex-a}, we briefly recall basic facts about the subconvexity problem for $L$-functions.
\subsection{Acknowledgements}

We would like to express deep gratitude to Peter Sarnak
for his interest in the results of this paper, as well as his encouragement
in preparing it. We would also like to thank L. Clozel and E. Lapid and D. Rama\-krishnan for helpful discussions as well as J.-P. Serre for his careful reading
 of the manuscript.

The present paper is part of a project that began on the
occasion of the AIM workshop ``Emerging applications of measure rigidity" on June
2004 in Palo Alto. It is a pleasure to thank the American Institute of Mathematics,
as well as the organizers of the workshop.

We would like to thank
Clay Mathematics Institute, the Sloan Foundation,
and the National Science Foundation for their financial support
during the preparation of this paper. A.V. would like to thank the Institute for Advanced Study
for its hospitality (and cookies) during the year 2005--2006. Ph. M. would like to thank Princeton University, the Institut des Hautes Etudes Scientifiques, and Caltech  for their hospitality on the
occasion of visits during various stages of elaboration of this paper.
 
M.E.\ was supported by the NSF grant DMS-0622397.
Ph.M. was partially supported by the ERC Advanced research Grant 228304.
A.V.\ was supported by a Sloan research fellowship.
M.E, E.L. and A.V. were supported by the collaborative grant DMS-0554365.

\section{An overview of the proof for $\PGL_3(\Z) \backslash \PGL_3(\R)$}
\label{outline}
The majority of this paper is presented in an ``adelized'' framework,
and the results presented are substantially more general than Theorem \ref{mainthm1}. Nonetheless, in this section, we would like to explain the ideas that go into Theorem
\ref{mainthm1} in as classical a setting as possible.

\subsection{Parameterizing compact orbits of maximal tori.}
Let us briefly recall how, to a number field $K$ with $[K:\Q]=n$,
we may associate a compact orbit of a maximal torus inside $\PGL_n(\R)$,
on the space $\PGL_n(\Z) \backslash \PGL_n(\R)$. 
If the field $K$
is totally real, the torus will be $\R$-split, 
we will be in the situation described in Theorem \ref{mainthm1},
and the construction will specialize to that discussed in our prior
paper \cite[Corollary 4.4]{ELMV1}. 

Fix a subalgebra $A_{r,s} \subset M_n(\R)$ isomorphic to $\R^r \oplus \C^s$.
Then $H_{r,s} := A_{r,s}^{\times}/\R^{\times}$ is a maximal torus
in $\PGL_n(\R)$; as $(r,s)$ vary through pairs satisfying $r+2s=n$,
they exhaust maximal tori, up to conjugacy.
In the case $(r,s) = (n,0)$, $H_{n,0}$ is conjugate to the diagonal subgroup $H$. 

Let $[K:\Q]=n$. Let $r$ and $s$ be the number of real and complex
embeddings of $K$, respectively. We shall say $K$ has {\em signature} $(r,s)$. 

Suppose given data $(K, L, \theta)$, where $K$ has signature $(r,s)$, 
 $\theta: K \otimes \R \rightarrow A_{r,s}$ is an algebra isomorphism,
 and $L$ is a ``$K$-equivalence class of lattices'' in $K$, i.e.
a free $\Z$-submodule of rank $[K:\Q]$, up to multiplication by $K^{\times}$. 

We associate to this data the $H_{r,s}$-orbit
 $\ov{\iota(L)}H_{r,s}$; here $\iota$ is any map
$K \otimes \R \rightarrow \R^n$ 
satisfying $(ab)^{\iota} = a^{\iota} . \theta(b)$. 
The resulting orbit is independent of choice of $\iota$. 
The stabilizer in $H_{r,s}$ of any point in the orbit is
 $\theta(\order_L^\times)$ where $\order_L = \{ \lambda \in K: \lambda L \subset L \}$, and the volume of the orbit
is $\reg(\order_L)$, the regulator of $\order_L$.

As explained in \cite[Corollary 4.4]{ELMV1}, in the totally real case $(r,s) = (n,0)$, {\em all} compact $H$-orbits in $X_{n}$ are obtained in that way.
This does not hold for the other signature (think of an imaginary quadratic field); to recover a form of this property, one has to consider more general $S$-arithmetic quotients
of $\PGL_{n}$.

{\em For clarity, we specialize in the remainder of this section to the case $(r,s) = (n,0)$
and  $H=H_{n,0}$ a maximal $\R$-split torus. Interpreted in an appropriate
$S$-arithmetic context (cf. \S \ref{explication}), most of the discussion
carries over to general signature so long as the field $K$
admits a fixed split place, and this is how our main result is proven in general.
}

\subsection{Packets for $\PGL_n$}\label{packetsPGLn}
Consider an order $\order$ inside a totally real field $K$; assume
we have fixed an identification $\theta$ as above.

 Let us denote by $\tilde{Y}_{\order}$ the set of $H$-orbits associated
to data $(K, L, \theta)$ such that $\order_L = \order$.
Varying $K$, $\order$ and $\theta$, the collections $\tilde{Y}_{\order}$
define a partition of the set of compact $H$-orbits. 
As it turns out, these form a slightly {\em coarser}
partition of the compact $H$-orbits than the one alluded to in Theorem \ref{mainthm1}; i.e., there is a {\em further} natural equivalence relation on the 
set of lattices with $\order_L = \order,$
 which is not trivial in general\footnote{However, when $\order$ is a {\em Gorenstein} ring - e.g.\ when  $n=2$, or $\order$ is the maximal order,
or $\order$ is monogenic -- this further equivalence is trivial, 
i.e.\ $\tilde{Y}_{\order}$ is a single packet in the sense of Theorem~\ref{mainthm1}.}.
  This equivalence is discussed and explicated in more detail in \S \ref{glnpackets},
 and we use the term {\em packets} to refer to its equivalence classes,
or to the associated collections of $H$-orbits.

 Assuming this for now, let $Y_{\order}$ be any packet of compact $H$-orbits
contained in $\tilde{Y}_{\order}$, and let $\mu_{\order}$ the corresponding measure.

It is not difficult to verify that $\reg(\order)$ goes to infinity with $\disc(\order)$. In the totally real case, the equidistribution statement of
Theorem \ref{mainthm1} is thereby equivalent to the statement,
in the $n=3$ case:
\begin{equation} \label{cubicequid} \mbox{As $\disc(\order)\ra+\infty$, $\mu_{\order}$ approaches Haar measure on $X_3$.}
\end{equation}

\subsection{Overview of proof for $X_2$ via analytic number theory}

To put things in perspective, it will be useful to recall the principle of Duke's proof for $X_{2}$.

Duke verifies {\em Weyl's equidistribution criterion}, i.e.
for a suitable basis $\{ \varphi \}$ for the functions in $L^2(\PGL_{2}(\Zz)\bash \PGL_{2}(\Rr))$with integral zero,
he shows:\begin{equation}\label{dukedecay}\mu_{\order}(\varphi):=\int_{X_{2}}\varphi(g)d \mu_{\order}(g)\ra 0,\ \disc(\order)\ra+\infty.
\end{equation}
The basis is chosen to consist of automorphic forms -- either cusp forms or Eisenstein series.
Duke proved \refs{dukedecay} by interpreting
the period integral $\mu_{\order}(\varphi)$ in terms of the $\disc(\order)$-th Fourier coefficient of an half-integral weight form
$\tilde \varphi$, proving non-trivial bounds for such coefficients
by generalizing a method of Iwaniec \cite{Iwaniec}. In its most general form, the formula relating the period integral to Fourier coefficients is due to Waldspurger \cite{Waldspurger1}.

Soon thereafter, another proof emerged that turned out to work in greater generality. Namely, by a result of Wal\-dspurger \cite{Waldspurger2}, one has the relation
$$|\mu_{\order}(\varphi)|^2=I_{\order}(\varphi)\frac{L(\pi,1/2)L(\pi \otimes \chi_{K},1/2)}{\disc(\order)^{1/2}}$$
where
\begin{enumerate}
\item $\pi$ is the automorphic representation to which $\varphi$ belongs;
\item  $L(\pi,s)$ and $L(\pi \otimes \chi_{K},s)$ are, respectively, the Hecke $L$-function of $\pi$ and the Hecke $L$-function of the twist of $\pi$ by the quadratic
character associated with $K$;
\item  $I_{\order}(\varphi)$ is a product of local integrals supported at the place at $\infty$ and at the primes dividing $\disc(\order)$.
\end{enumerate}
Then $\refs{dukedecay}$ is a consequence of the estimates
\begin{equation}\label{subconvex2}
L(\pi \otimes \chi_{K},1/2)\ll \disc(\order_{K})^{1/2-\eta},\
I_{\order}(\varphi)\ll \bigl(\frac{\disc(\order)}{\disc(\order_{K})}\bigr)^{1/2-\eta}
\end{equation}
for some absolute $\eta>0$; here $\order_{K}$ denote the maximal order of $K$.

The first bound in \refs{subconvex2}
is called a {\em subconvex} bound and is due to
Duke, Friedlander and Iwaniec \cite{DFI1}; it is a special instance of the so-called {\em subconvexity problem} for automorphic $L$-functions (see \cite{IS} for a discussion of that problem).

The second bound in greater generality is due to Clozel and Ullmo \cite{CU} and we will call it a {\em local} subconvex bound.
It is somewhat easier that the first, but it addresses the issue
of $\order$ non-maximal.

It is tempting to try to generalize this approach to the space $X_{n}$ of lattices of higher rank. However this does not seem within reach of the current technology. In particular, it is not expected that the corresponding ``Weyl sums'' are related to ($\GL_{n}$) $L$-functions. Even were this the case,
we do not know how to solve the corresponding subconvexity problems.\footnote{
Another possibility, more in line with Duke's original proof
would be to use results of Gan, Gross and Savin \cite{GGS} which relate the Weyl's sums to Fourier coefficient to automorphic forms on $G_{2}$; unfortunately our state of knowledge concerning bounds for  these Fourier coefficients is rather limited. }
There is however an exception to this which will turn out to be crucial for our coming argument.

\subsection{Overview of the proof for $X_3$}

In summary, our strategy is to check Weyl's equidistribution criterion against test functions taken from a {\em tiny} portion of $L^2(X_{3})$,
by using/proving
global and local subconvex bounds,
and to {\em bootstrap} this information to a full equidistribution statement using
results on measures invariant under higher rank torus actions.

The input we use from ergodic theory is the following measure classification 
 result regarding invariant measures in rank $\geq 2$ --- as well as a $p$-adic variant of it -- by the first two authors and A. Katok \cite{EKL}:

\begin{Theorem}\label{EKL} Let $n \geq 3$ and let $\mu$ be an ergodic
$H$-invariant probability measure on $X_{n}$ where $H$ denotes a maximal $\R$-split torus in $\PGL_{n}(\R)$. Assume that for at least one $a \in H$
the ergodic theoretic entropy $h_{\mu}(a)$ is positive.
Then $\mu$ is homogeneous: there exists a reductive group $H \subsetneq L \subset \PGL_{n}(\Rr)$ such that $\mu$ is the $L$
invariant probability measure on a single periodic $L$-orbit. 
In particular, if $n$ is prime, $\mu$ is Haar measure on $X_{n}$.
\end{Theorem}

The main content of this paper will be to show that assumptions of this theorem are satisfied when $n=3$.

It is conjectured that the following substantially stronger statement holds:

\begin{Conjecture} [{Furstenberg, Katok-Spatzier \cite{Katok-Spatzier}, Margulis \cite{Margulis-conjectures}}] \label{FKSM}
Let $n \geq 3$ and let $\mu$ be an ergodic
$H$-invariant probability measure on $X_{n}$, \ $H$ as above. Then $\mu$ is  homogeneous.
\end{Conjecture}

Note that in Conjecture~\ref{FKSM}, the measure $\mu$ can certainly be the natural measure on a periodic $H$-orbit, a possibility that is ruled out in Theorem~\ref{EKL}.
We refer the interested reader to \cite{EL-ICM} or to the original paper \cite{EKL} for a historical background and for an exposition of some of the ideas that enter into the proof.

Let $\mu_{\infty}$ denote a $\mathrm{weak}^*$ limit\footnote{recall that a sequence of probability measures
 $\{\mu_{i}\}_{i}$ $\mathrm{weak}^*$ converge to some measure $\mu_{\infty}$ if, for any compactly supported function $f$,
$\mu_{i}(f)\rightarrow \mu_{\infty}(f)$ as $i\rightarrow+\infty$.} of the $\{\mu_{\order}\}_{\order}$. There are two main issues to verify:
\begin{enumerate}
\item[(A)]\label{compact} The  measure $\mu_{\infty}$ is a probability measure (i.e.\ the sequence of measures $\{ \mu_{\order}\}_{\order}$ is {\em tight}).
\item [(B)] \label{posentropy} Almost every ergodic component of $\mu_{\infty}$ has positive entropy w.r.t some $a \in H$.
\end{enumerate}

Even assuming the stronger conjectured measure classification given by Conjecture~\ref{FKSM}  one needs to overcome pretty much the same obstructions;
 in that case the following weaker form of (B) would suffice:
\begin{enumerate}
\item[(B$'$)] $\mu _ \infty (xH) = 0$ for any periodic $H$-orbit $xH$.
\end{enumerate}

\noindent

In the context of this paper (B$'$) does not seem to be much easier to verify than the weaker statement in (B).

\subsection{Weyl's equidistribution criterion}

Our method for verifying both (A) and (B) is by checking {\em Weyl's equidistribution criterion} for a special class of functions from which follows {\em a priori bounds} for the
$\mu_{\infty}$-volumes of certain sets. We shall be able
to obtain such bounds on the mass of neighborhoods of the cusp in $X_{3}$ (which addresses (A)) or of $\eps$-balls around {\em any} $x \in X_{3}$
(which addresses (B)).

\subsubsection{The Siegel-Eisenstein series}
Let us recall that we can identify $X_3$ with the space of lattices in $\R^3$ of co-volume $1$. We shall make use of this identificaition throughout what follows.
Let $f$ be any continuous, compactly supported function on $\R^3$ and let $E_f$
be the Siegel-Eisenstein series $$E_f(L) = \sum_{\lambda \in L-\{ 0 \}} f(\lambda); $$ we shall prove
\begin{equation} \label{siegel} \mu_{\infty}(E_f) :=\lim_{\disc(\order)\ra+\infty}\mu_{\order}(E_f)= \int_{\R^3} f(x) dx.\end{equation}
Observe that, by {\em Siegel's formula}, (see eg. \cite{Weil}), $\int_{\R^3} f(x) dx=\mu_{\operatorname{Haar}}(E_{f})$; in particular \refs{siegel} is consistent with \refs{cubicequid}.

By taking suitable choices of $f$, \eqref{siegel}
yields the necessary {\em a priori} bounds.
Indeed, take $v \in \Rr^3$. Take $f$ to be a smooth non-negative function supported in the $2 \eps$-ball $B(v,2\eps)$,  which takes value $1$ on $B(v,\eps)$.
When $v=0$,  $E(f)$ dominates a neighborhood of the
cusp (in fact approaches infinity near the cusp); 
when $v \neq 0$, $E(f)$ dominates the characteristic function of an $\eps$-neighborhood
of any lattice class $x$ that contains $v$. In the latter case,
we deduce that for, $\eps$ small enough
$$\mu_{\infty}(B(x,\eps))\ll \eps^3.$$
This improves over the trivial bound $\mu_{\infty}(B(x,\eps))\ll \eps^2$,
arising from the fact that $\mu_{\infty}$ is invariant by
the two parameter group $H$.

This improvement\footnote{In passing,  we note that the test functions ``$E_f$'' considered are {\em not} invariant under the maximal compact $K=\PO_{3}(\Rr)$,
in general (unlike many problems using the classical theory of modular forms). This feature is essential to improve the trivial bound $\mu_{\infty}(B(x,\eps))=O(\eps^2)$  to $O(\eps^3)$. }
 from $2$ to $3$ already shows that $\mu_{\infty}$ cannot be supported along a compact $H$-orbit. More importantly, this implies that for a generic $a \in H$, almost every ergodic component of $\mu_{\infty}$ has positive
entropy with respect to the action of $a$ from which we deduce the full equidistribution by Theorem \ref{EKL}. To finish this section, we remark that the principle of testing Weyl's criterion
against Einsenstein series appears in other contexts, for example,
 \cites{EMM, Veech}.

\subsubsection{Connection to $L$-functions}\label{ConnectionToLFunctions}
The key point, for establishing \eqref{siegel}, is that the $\mu_{\order}(E_{f})$ is indeed related to $L$-functions.
One has the following formula, which goes back to Hecke \cite{Heckewerke}\footnote{Hecke proved that way the analytic continuation and the functional equation
of Gr\"ossencharacter $L$-functions.}
\begin{equation}\label{HeckeFormula}\mu_{\order}(E_{f})=\int_{\Re s \gg 1}\widehat f(s)\frac{I_{\order}(f,s)\zeta_{K}(s)}{\disc(\order)^{1-s}}ds,\end{equation}
where $\widehat f(s)$ is a certain Mellin-like transform of $f$, $\zeta_{K}(s)$ is the Dedekind zeta function of $K$ and $I_{\order}(f,s)$ is a product of local integrals supported at the place $\infty$ and at the primes dividing $\disc(\order)$.

Shifting the contour to $\Re s=1/2$, we pick up a residue at $s=1$ which equals $\mu_{Haar}(E_{f})$;
the fact that the remaining integral goes to $0$ as $\disc(\order)\ra+\infty$ follows from the global and local subconvex bounds
\begin{equation}\label{globalsubconvex}\zeta_{K}(s)\ll_{s}\disc(\order_{K})^{1/2-\eta},
\end{equation}
\begin{equation}\label{localsubconvex}\ I_{\order}(f,s)\ll_{s}\bigl(\frac{\disc(\order)}{\disc(\order_{K})}\bigr)^{1/2-\eta}
\end{equation}
for some absolute $\eta>0$ and $\Re s=1/2$.

For $n=3$ the  bound \refs{globalsubconvex} follows from the work of Burgess \cite{Bu} if $K$ is abelian and (essentially)
from the deep work of Duke, Friedlander and Iwaniec if $K$ is cubic not abelian (\cites{DFI8,BHM,MV}).

The local bound  \refs{localsubconvex} is new and occupies a good part of the present paper; moreover, it is valid for any $n$. Let us describe how it is proved.

\subsubsection{The local estimates: harmonic analysis on $p$-adic homogeneous spaces} \label{sketch}
Our approach to bounding the local integrals $I_{\order}$ is based on
\footnote{In fact, the estimates needed can be proved in a direct and elementary way;
this was the original approach of the paper. However, the approach carried out, although requiring more input, has the advantage of being very general.}
, first of all, relating $I_{\order}$ to integrals of matrix coefficients (inspired by ideas of Waldspurger and Ichino-Ikeda) and then bounding
the integrals of matrix coefficients using the local building (inspired by ideas of Clozel and Ullmo).

Let us make a remark in representation theory to explain this.
Let $V$ be an irreducible, unitary representation of a group $G$, and let $H \subset G$ be a subgroup. Suppose that there exists a unique scaling class of invariant functionals $L: V \rightarrow \C$ invariant by $H$.

It is sometimes possible to understand something about functionals $L$
simply by studying matrix coefficients. Indeed, if convergent,  $\int_{h \in H} \langle h v_1, v_2 \rangle$ defines a functional (on $v_1$) invariant
by $H$ and a conjugate-linear functional (on $v_2$) invariant by $H$. We conclude by the uniqueness assumption that:
\begin{equation} \label{Basic} \int_{h \in H} \langle h v_1, v_2 \rangle dh = c_V L(v_1) \overline{L(v_2)} \end{equation}
for some constant $c_V$. Thereby, $L$ can be studied through
matrix coefficients.

It turns out that computing $I_{\order}$ amounts to computing
with such functionals $L$,
when $G = \GL(n,k)$ and $H$ is a maximal torus inside $G$;
the representations $V$ we are concerned with are those
that occur inside $L^2(k^n)$. Thereby, \eqref{Basic} allows
us to reduce understanding $I_{\order}$ to computations with matrix coefficients.

\section{Number theoretic interpretation} \label{NT}
We discuss how our main result can be interpreted in terms
of integral points on varieties, generalizing the equations
$(*)$ and $(**)$ from the introduction.

\subsection{Integral points on homogeneous varieties}

In this section, we interpret our results in terms of distribution of algebraically defined sets of integral matrices,
which was one of Linnik's original motivations. This is part of a more general problem of studying the structure of the set of integral
 points $\mathbf{V}(\Z)$ on an algebraic variety $\mathbf{V}$.

A particularly structured situation occurs when $V$ is homogeneous, i.e.\ $\mathbf{V}(\C)$ possesses a transitive action of a linear algebraic group $\G$.
In that case, it is expected that there are many points which are rather well distributed.
 Similar results are found in  \cites{Linnik,Sarnak,EMS,EO,EV}.

\subsection{Solving polynomial equations in matrices}
To motivate what follows, note that for $(x,y_1,y_2, z) \in \Zz^4$ and $m$ not a perfect square, 
$$\left( \begin{array} {cc} y_1 & z \\ x & y_2 \end{array}\right)^2 = m \, \mathrm{Id}
 \Longleftrightarrow  y_1 = -y_2, y_1^2 + xz = m.$$
Thereby, $(**)$ is a statement concerning $2 \times 2$ integral
matrices satisfying a prescribed quadratic equation. 

Given $P$ a monic integral irreducible polynomial of degree $n$ with integral coefficients, we let
$$Z_{P}=\{M\in M_{n},\ P(\lambda)=\det(M-\lambda I)\}.$$
Thus integral points $Z_P(\Zz)$ can be identified simply with integral solutions to $P=0$ in $n \times n$ matrices. 

The {\em signature} $(r,s)$ of such a polynomial
will be the number of real roots, resp. conjugate pairs of complex roots;
thus $r+2s=n$.
Let $Z_{r,s}$ be the space of
all splittings of $\R^n$
 into $r$ real lines and $s$ complex planes.

If $P$ has signature $(r,s)$, the space $Z_P(\Rr)$ is identified with
$Z_{r,s}$ by, first, fixing an ordering of the real and complex roots of $P$;
and then associating to a matrix $M \in Z_P(\Rr)$ its eigenspaces.\footnote{For a complex eigenvalue, we take the eigenspaces corresponding to that eigenvalue and its complex conjugate, and intersect their sum with $\Rr^{n}$.}

The spaces $Z_{r,s}$ carry a $\PGL_n(\R)$-invariant measure, unique up to scaling (indeed, $Z_{r,s} \cong \PGL_n(\R) / H_{r,s}$), which we denote by $\vol(\cdot)$.

\begin{Corollary}\label{IPCor}
Let $\{P_{i}\}_{i}$ be a sequence of cubic, monic, integral, irreducible
 polynomials of signature $(r,s)$ and of
discriminant satisfying $\disc(P_{i})\ra+\infty$. 

\begin{enumerate}
\item If $(r,s) = (3,0)$, then
 $Z_{P_i}(\Zz)$ becomes uniformly distributed on $Z_{3,0}$.
\item If $(r,s) = (1,1)$, and there exists a fixed prime number $p$
so that $P_i$ has three $p$-adic roots, then
$Z_{P_i}(\Zz)$ becomes uniformly distributed on $Z_{1,1}$. 
\end{enumerate}
\end{Corollary}
 
Here, we say, a sequence of discrete sets $\mathcal{Z}_i \subset Z_{r,s} \cong G/H_{r,s}$
 becomes {\em equidistributed on $Z_{r,s}$} if, for any compact sets $\Omega_{1},\Omega_{2}\subset Z_{r,s}$ with  boundary measure zero and $\vol(\Omega_{2})>0$, one has
 $$\frac{|\mathcal{Z}_i \cap \Omega_{1}|}{|\mathcal{Z}_i \cap  \Omega_{2}|}\ra \frac{\vol(\Omega_{1})}{\vol(\Omega_{2})}, \mbox{ as }i \ra \infty.$$
 Implicit in this statement is the fact that $|\mathcal{Z}_i \cap  \Omega_{2}|$ is non-zero if $i$ is large enough:
  for instance, in Corollary \ref{IPCor}, one can show that for any $\eps>0$ and $i$ sufficiently large (depending on $\Omega_{2}$)
$$|{Z}_{P_{i}}(\Zz) \cap  \Omega_{2}|\gg_{\eps}\disc(P_{i})^{1/2-\eps}\vol(\Omega_{2}).$$

\subsection{Cube roots of integers in $3 \times 3$ matrices.}
Let us specialize further, to make this even more concrete. For $d>0$ not a perfect cube,
consider the polynomials $P_{d}(X)=X^3-d$  and set $Z_{d}(\Zz)=Z_{P_{d}}(\Zz)$. 

Here, it is convenient to explicitly interpret $Z_{1,1} \cong G/H_{1,1}$ as the space of "matrix cube roots of unity": $$Z_{1,1}=\{M\in M_{3}(\Rr),\ M^3=\mathrm{Id},M\neq \mathrm{Id}\}.$$ 
This being so,
our previous Corollary can be stated in terms of the ``radial projections''
to the latter space:

\begin{Corollary}\label{corconj}Let $p>3$ be a fixed prime. As $d\ra+\infty$ amongst the integers which are not perfect cubes and such that
$$\hbox{$p$ is {\em totally split} in the field $\Qq(\sqrt[3]{d})$}$$
 then the sequence of sets $\frac{1}{d^{1/3}}.\{M \in M_3(\Zz), M^3 = d\}$ becomes equidistributed in the space
$\{M \in M_3(\Rr), \ M^3 = \mathrm{Id}\}.$
\end{Corollary}

\subsection{Translations}
Let us explain how the above
corollaries follow, indeed, from our
main theorems. Set $G=\PGL_{n}(\Rr)$ and $\Gamma=\PGL_{n}(\Zz)$.

Let $P$ have signature $(r,s)$; let $\order_P$
be the ring $\Z[t]/P$ and $K_P = \order_P \otimes \Q$.   By a {\em coarse ideal class for $\order_P$}, we understand a lattice $L \subset K_P$
so that
 $\order_P . L \subset L$, considered up to multiplication by $K_P^{\times}$. With this convention, there are maps:
\begin{equation} \label{bij}
\mbox{$\Gamma$-orbits on $Z_P(\Zz)$ $\leftrightarrow$ coarse ideal classes for $\order_P$
 $\rightarrow$ compact $H_{r,s}$-orbits on $X_n$.}
\end{equation}

The first map is a bijection. 
The composite of the two arrows amounts to the identification
between $\Gamma$-orbits on $G/H_{r,s}$, and $H_{r,s}$-orbits on $\Gamma \backslash G$. 

In arithmetic terms, we can understand the maps as follows:
\begin{enumerate}
\item If we fix $M\in Z_{P}(\Zz)$, the map $t \mapsto M$
makes $\Zz^n$ into a $\order_P$-module; there is a unique coarse
ideal class $L$ so that $L$ and $\Zz^n$ are isomorphic as $\order_P$-modules. 
(This is very classical; see, e.g.\ \cite{LMD}). 
\item
The second injection associates to the class of
$L \subset K_P$ the compact orbit associated to $(K_P, L, \theta)$,
in the notation of \S \ref{packetsPGLn}.
Here $\theta: K_P \rightarrow A_{r,s}$ is the identification
arising from the chosen ordering of the roots of $P$. 
\end{enumerate}
The composite map associates to the set
of $\Gamma$-orbits on $Z_P(\Zz)$ a set of $H_{r,s}$-orbits on $X_n$; in the notation of \S \ref{packetsPGLn},
this set is:
\begin{equation} \label{ypdef} Y_{P}=\bigcup_{\order_P \subset \order\subset\order_{K_{P}}}\tilde Y_{\order},\end{equation} corresponding to all packets whose associated order in $K_P$
contains $\order_P$.  The possibility of intermediate orders
corresponds to the fact that the ideals that arise need not be {\em proper} $\order_P$-ideals. 

Under these bijections, the equidistribution assertion about $Z_P$
translates to an equidistribution assertion about $Y_P$, as we now recall.

\subsection{Integral-points interpretations} 
As is well known (cf. \cite[\S 8]{Benoist-Oh}) the (tautological) equivalences
$$\Gamma\bash\Gamma gH_{r,s}\longleftrightarrow \Gamma g H_{r,s}\longleftrightarrow\Gamma gH_{r,s}/H_{r,s}$$
can be used to transfer equidistribution results about periodic $H_{r,s}$-orbits on $\Gamma \backslash G$,
to equidistribution of the corresponding $\Gamma$-orbits on $G/H_{r,s}$. 
Taking into account \eqref{bij}, 
the first assertion of Corollary \ref{IPCor}
reduces to the following:

\begin{Corollary}\label{totrealcor} Let $\{P_{i}\}_{i}$ be a sequence of cubic, monic, integral, irreducible, $\Rr$-split polynomials of discriminant $\disc(P_{i})\ra+\infty$. Then the set of compact $H$-orbits defined by \eqref{ypdef}
 becomes equidistributed on $X_{3}$ (here $H=H_{3,0}$).
\end{Corollary}
This is indeed a corollary to Theorem \ref{mainthm1}, 
taking into account the fact that the total number of compact $H$-orbits on $X_3$ with bounded
volume is finite. 

Similarly, the other assertion of Corollary \ref{IPCor}
follows from the more general
adelic Theorem \ref{RealTheorem}.

\section{Homogeneous subsets in the adelic context} \label{Hsets}
This paper has been written consistently in the {\em adelic} framework.
It is therefore appropriate for us to discuss adelic equidistribution problems. 
We confine ourselves to equidistribution problems
associated to tori, although much of the discussion applies in greater generality. 

Let us emphasize that the adeles are simply a linguistic tool: 
all statements and results could be readily stated in the $S$-arithmetic context. The advantage of the adeles, rather, is that they provide a unified
approach to broad classes of questions.  For instance,
consider the following equidistribution questions on the modular surface:
\begin{enumerate}
\item  Equidistribution of CM points;
\item Equidistribution of large hyperbolic circles, centered
at the point $i \in \mathbb{H}$;
\item Equidistribution of closed geodesics (see \S \ref{sec:GI}). 
\item On $\PGL_2(\Z) \backslash \PGL_2(\R)$,
equidistribution of the translate
  of a {\em fixed}
closed $H$-orbit by a ``large'' group element in $\PGL_2(\R)$.\footnote{See \cite{EMS}
for very general theorems concerning this setting, and also \cite{Benoist-Oh}
for related results.}
\end{enumerate}

These situations are all closely related, although
they are often treated separately, and our aim
is to discuss them as specializiations of a single adelic context. 

Similarly, in \cite{CU}, two classes of equidistribution problems are considered: ``equidistribution on the group'', and ``equidistribution on the symmetric space;'' these two problems again become unified in our presentation.

Explicitly, the goals of this section are as follows.
We define ``homogeneous toral sets'';
roughly, as will be clarified in Theorem \ref{thmpacket},
these generalize the groupings $Y_{x_{i}H}$ of compact
orbits discussed in Theorem \ref{mainthm1}. 
  We then define 
two important invariants (``volume'' and ``discriminant'')
for homogeneous toral sets, 
formulate the main question about their distribution (\S \ref{desid})
and state our main theorems -- Theorem \ref{Volumes} and Theorem \ref{RealTheorem} -- in these terms.  These theorems imply immediately Theorem \ref{mainthm1},
and their proofs comprise most of this paper.

\subsection{Homogeneous sets: definitions}
Let $F$ be a number field with adele ring $\adele$. 
Let $\G$ be a $F$-group with Lie algebra $\g$. 
Set $X = \G(F) \backslash \G(\adele)$. 

A {\em homogeneous toral subset}
of $X_{\adele}$ will be, by definition, one of the form
 $$Y = \T(F) \backslash \T(\adele). g_{\adele},$$
when $g_{\adele} \in \G(\adele)$ and $\T \subset \G$ is a maximal torus. 
We shall consider only the case where the torus $\T$ is anisotropic over $F$. 

 Then $Y$ supports a natural probability measure $\mu_Y$: the pushforward of the Haar probability measure on $\T(F) \backslash \T(\adele)$ by $h \mapsto hg_{\adele}$.

We shall associate to $Y$ 
two additional invariants: a {\em discriminant} $\disc(Y)$,
measuring its arithmetic complexity, and a {\em volume} $\vol(Y)$, measuring
how ``large'' it is.

\subsection{Discriminant.} \label{subsec:disc}
Let $r = \dim \T$. Let $V$ be the affine space $(\wedge^{r} \g)^{\otimes 2}$. 
We fix a compatible system of norms $\|\cdot\|_v$ on $V \otimes_F F_v$, for each place $v$ (for a discussion of norms, see \S \ref{local}; ``compatible'' means
that, for almost all $v$, the unit balls of the norms coincide
with the closure inside $V \otimes F_v$ of a fixed $\order_F$-lattice within $V$).

 To the Lie algebra $\t$ of $\T$, we associate
a point in the affine space $(\wedge^{r} \g)^{\otimes 2}$:
 \begin{equation} \label{embedding}\iota(\mathfrak{t}) = (e_1 \wedge \dots \wedge e_r)^{\otimes 2}
(\det B(e_i, e_j))^{-1}.\end{equation}
where $e_1, \dots, e_r$ is a basis for $\mathfrak{t}$, and $B$ the Killing form. 
We set $$\disc_v(Y) = \|\Ad(g_v^{-1}) \iota(\mathfrak{t})\|_v.$$ The discriminant $\disc(Y)$ is defined to be the product
$$\prod_{v} \disc_v(Y).$$

\subsection{Volume.} \label{Voldef}
The definition of ``volume,''  for a homogeneous toral subset, will depend on the choice of a compact
neighborhood  $\Omega_0 \subset \G(\adele)$ of the identity. 
 This notion depends on $\Omega_0$, but the notions arising from two different choices of $\Omega_0$
are comparable to each other, in the sense that their ratio is bounded above and below. We define
\begin{equation} \label{volume}
\vol(Y) := \vol \bigl(
\{t \in \T(\adele):g_{\adele}^{-1} t g_{\adele} \in \Omega_0\}
\bigr)^{-1},
\end{equation}
where we endow $\T(\adele)$ with the measure that assigns the quotient
$\T(F) \backslash \T(\adele)$ total mass $1$.

\subsection{Desideratum.}\label{desid}
We shall say that a measure on $X$ is {\em homogeneous} if it is supported on a single orbit of its stabilizer. 

The kind of problem we are interested in is the following:

{\em When $\disc(Y_i) \rightarrow \infty$ (equivalently $\vol(Y_i) \rightarrow \infty$),
show that $\mu_{Y_i}$ converges to an homogeneous measure.}

There are certain cases of this problem which are easier.
For instance (the ``depth'' aspect) we might consider
a sequence of homogeneous toral sets for which
there exists a fixed place $v$ with $\disc_v(Y_i) \rightarrow \infty$. In this case, a limit of the $\mu_{Y_i}$s will be invariant under a unipotent subgroup.
This special case is interesting from the point of view
of many applications, as for instance in the work of Vatsal \cite{Va}.
 Eskin, Mozes and Shah,
\cite{EMS}, and also Benoist and Oh \cite{Benoist-Oh}, also study this aspect.

\subsection{The case of a quaternion algebra.} \label{sec:quaternion}
The current state of knowledge concerning quaternion algebras implies the following theorem:

\begin{Theorem} \label{CU}
Let $\G$ be the projectivized group of units in a quaternion algebra over
a number field $F$. 
Let $\{Y_{i}\}_{i}$ be a sequence of homogeneous toral sets whose discriminant 
approaches $\infty$ with $i\ra+\infty$.  Then
$$\vol(Y_i)=\disc(Y_i)^{1/2+o(1)},\ i\ra+\infty.$$
Moreover, any $\mathrm{weak}^*$ limit of the measures $\mu_{Y_i}$ is a homogeneous probability
measure on $\G(F) \backslash \G(\adele)$, invariant under the image of $\widetilde{\G}(\adele) \mapsto \G(\adele)$. 
\end{Theorem}

Here, $\widetilde{\G}$ denotes the simply-connected covering-group of $\G$. 
Indeed, one even knows this in a quantitative form: if $f \in C^{\infty}(\G(F) \backslash \G(\adele))$
generates an irreducible, infinite-dimensional $\G(\adele)$-representation, then $|\mu_{Y_i}(f) - \mu(f)|$ is bounded by $O_f (\disc(Y_i)^{-\delta})$
for a positive $\delta$.  

The reason that we cannot simply assert that $\mu_{Y_i}$ converge to the Haar measure has to do with ``connected component issues.'' 

Theorem \ref{CU} is a consequence of works by several authors:
\begin{enumerate}
\item Siegel's lower bounds (for the statement concerning volumes);
\item Iwaniec \cite{Iwaniec}, Duke \cite{Duke}, Duke--Friedlander--Iwaniec \cite{DFI1}, Cogdell--\-Piatetski-Shapiro--Sarnak
\cite{CPSS}, and the fourth-named author \cite{Ve} (for the pertinent subconvexity bounds). 
\item Waldspurger \cite{Waldspurger2} (see also \cite{Katok-Sarnak}), Clozel--Ullmo \cite{CU}, Popa \cite{Popa}, S.-W. Zhang \cites{Zhang2} and P. Cohen \cite{Cohen}. 
\end{enumerate}
  It conceals a unified statement
of a large number of ``instances'' of that theorem, corresponding
to varying the parameter $\disc$ in different ways:
e.g.\ \cites{Chelluri, MV, CU}. For instance, if the quaternion algebra is defined over a totally real field and the quaternion algebra is ramified at one place,
 one obtains in that way, equidistribution results for
closed geodesics on an arithmetic Riemannian surface.  Another example: it 
implies the solution (outlined by Cogdell--Piatetsky-Shapiro--Sarnak in \cite{CPSS}) 
to the representability question for ternary quadratic forms over number fields. 

\subsection{Results for $\G=\PGL_n$.}

Let $\{Y_i\}_{i}$ be any sequence of homogeneous (maximal) toral sets on $X=\PGL_{n}(F)\bash\PGL_{n}(\adele)$ whose discriminant approaches $\infty$ with $i\ra+\infty$. 
Let $\T_i$ be the associated tori; then $$\T_i = \mathrm{Res}_{K_i/F}\mathbb{G}_{m,K_{i}}/\mathbb{G}_{m,F},$$
 for a field extension $K_i/F$ of degree $n$, unique up to isomorphism.   
We show:
\begin{Theorem}\label{Volumes} For $\{Y_{i}\}_{i}$ as above, one has
$$\vol(Y_{i})=\disc(Y_{i})^{1/2+o(1)},\ as\ i\ra+\infty$$  
\end{Theorem}
This result is easy given well-known (but difficult) bounds
on class numbers.  It shows that the definitions of adelic volume and discriminant proposed are compatible.

Now let us describe our result on the distribution of homogeneous toral sets. First
suppose there exists a fixed place $v$ with one of the following properties:
\begin{enumerate}
\item The local discriminant $\disc_v(Y_i) \rightarrow \infty$. 
\item Every $\T_i$ is split at $v$ and a sub-convexity result {\em in the discriminant aspect} is known
for values, along the critical line, of the Dedekind-$\zeta$-functions associated
to the fields $K_i$. (See \eqref{Hypmod} for the precise requirement.)
\end{enumerate}
Then our results establish that any $\mathrm{weak}^*$ limit of the measures $\mu_{Y_i}$ is a convex combination
of homogeneous probability measures. However, the precise shape of such a homogeneous probability measure appears to be somewhat 
complicated in the adelic setting.

Rather than attempt a precise statement of the above,
let us simply give the result in the simple case $n=3$. It is simple for two reasons: first of all, the necessary subconvexity is known; secondly,
the fact that $n$ is prime forces there to be very few intermediate measures. 
We prove:\footnote{It should be observed that this relies on an extension of \cite{DFI8}
that has been announced by the latter two authors \cite{MV-ICM},
and a theorem announced in \cite{EL-ICM}, but neither of the proofs
have yet appeared.  With $F=\Q$, the proofs exist in print and are
 contained in \cites{EKL, DFI8, BHM}. }

\begin{Theorem} \label{RealTheorem}Suppose $n=3$ and let $\{Y_{i}\}_{i}$ be a sequence of homogeneous toral sets
such that $\disc(Y_{i})\ra+\infty$ with $i$; suppose there exists a place $v$ so that $\disc_v(Y_i) \rightarrow \infty$ or so that each $\T_i$ is split at $v$. 
Then any $\mathrm{weak}^*$ limit of the $\mu_{Y_i}$, as $i \rightarrow +\infty$,
is a homogeneous probability measure on $X$,  invariant by the image of $\SL_3(\adele)$. 
\end{Theorem}

The equidistribution assertion of
Theorem \ref{mainthm1} is a consequence of Theorem \ref{RealTheorem}, applied
with $F=\Q, v= \infty$; translation from Theorem \ref{RealTheorem}
is provided in the next section. 

 We conclude this section by observing that it remains a very interesting problem
to remove the usage of the place $v$ in Theorem \ref{RealTheorem},
i.e., obtaining for $\PGL_3$ a result as strong as Theorem \ref{CU} (even without a rate).

%%%%%%%%%%%%%%%%%%%%%%%%%%%%%%%%%%%%%%%%%%%%%%%%%%%%%%%%%%%%%%%%%%%%%%%%%%%%%%%%%%%%%%%%%%%%%
\section{Packets} \label{data}
In this section we will clarify the relationship between
the adelic perspective of \S \ref{Hsets}, and the classical perspective
of \cite{ELMV1}. 

We will therefore exhibit a natural equivalence relation
on the set of compact $H$-orbits on $\Gamma \backslash G$
for which the equivalence classes are (almost) finite abelian groups. 

The equivalence classes will be called packets, and the union of compact
torus orbits in a packet corresponds, roughly speaking, to an adelic torus orbit.

\subsection{Notation.}
We recall the data prescribed in \cite{ELMV1}. 
Let $\G$ be a semisimple group over $\Q$ that is $\R$-split;  $G = \G(\R), \Gamma \subset G$ a congruence
lattice, and $H$ a Cartan subgroup of $G$.  To simplify notations we will write $\Gamma g H$ for the right $H$-orbit
 $\Gamma\bash\Gamma g H\subset \Gamma\bash G$.

We fix a lattice $\g_{\Z} \subset \g$ that
is stable by the (adjoint) action of $\Gamma$, as well as a $G$-invariant bilinear form
$B(\cdot, \cdot)$ on $\mathfrak{g}$ with $B(\g_\Z,\g_\Z)\subset\Z$.  Finally, we fix
a Euclidean norm on $\mathfrak{g}_{\R}$. 

Let $r$ be the rank of $G$.
Take $V = \wedge^r \mathfrak{g}$. 
For all finite $p$, we endow $V\otimes\Q_p$ with the norm which has as unit ball the closure
of $\bigl(\wedge^r\g_\Z\bigr)^{\otimes 2} \otimes\Z_p$. 
For $p=\infty$, we give $V \otimes \R$ the Euclidean norm derived
from that on $\mathfrak{g}$. 
These choices allow one to define the notion of discriminant of a homogeneous toral set, as in \S \ref{Hsets}. 

 In addition to this, we shall take {\em as given} one further piece of data. Let $\adelef$ be the ring of finite adeles.
Choose a compact open subgroup
$K_f \subset \G(\adelef)$ so that $K_f \cap \G(\Q) = \Gamma$ and so that $\g_{\Z}$ is stable
under the (adjoint) action of $K_f$. 

Let $X_{\adele}$ denote the double quotient $X_{\adele} := \G(\Q) \backslash \G(\adele) / K_f$. Clearly $G$ acts on $X_{\adele}$
 and we shall refer to its $G$-orbits as the {\em components} of $X_{\adele}$   and
to  the orbit of the identity double coset  as the {\em identity component} (these need not be topologically connected);
 the identity component is identified with $\Gamma \backslash G$.
  The set of components  is finite and is parametrized by the double quotient $\G(\Q) \backslash \G(\adelef) / K_f$.

\begin{Theorem}\label{thmpacket}

\begin{enumerate}
\item 
Each compact $H$-orbit $\Gamma g H \subset \Gamma \backslash G \subset X_{\adele}$
is contained in the projection to $X_{\adele}$ of a 
homogeneous toral set $Y \subset \G(\Q) \backslash \G(\adele)$.
The set $Y$ is unique up to translation by $K_f$; in particular, all such 
$Y$s have the same projection to $X_{\adele}$. 
Moreover, the discriminant of $\Gamma g H$, in the sense
of \cite{ELMV1}, and the discriminant of $Y$, in the sense of \S \ref{Hsets},
coincide up to a positive multiplicative factor; the latter factor depending only on $H$, on the choice of $B(\cdot,\cdot)$ and on the 
Euclidean norm on $\g_\R$.

\item Declare two compact $H$-orbits to be {\em equivalent} 
if they both are contained in the projection to $X_{\adele}$
of a homogeneous toral set $Y$. 

An equivalence class of compact $H$-orbits we refer to as a {\em packet}.
Packets are finite; indeed, the
  packet of $\Gamma g H$ is parameterized by the fiber of the map
$$\T(\Q) \backslash \T(\adelef)/ (K_f \cap \T(\adelef)) \rightarrow \G(\Q) \backslash \G(\adelef) / K_f$$
above the identity double coset; here $\T$ is the unique $\Q$-torus
so that $\T(\R) = g H g^{-1}$.

In particular, if $\G(\Q) \backslash \G(\adele)/K_f$ has a single component, every packet naturally has the structure of a principal
 homogeneous space for a finite abelian group. 

\item Compact orbits in the same packet have the 
same stabilizer and the same discriminant.
\end{enumerate}
\end{Theorem}

The proof of this theorem is straightforward. However, its content
is quite beautiful: the collection of compact Cartan orbits on $\Gamma \backslash G$ group themselves into equivalence classes, each (almost -- see assertion (2) of the Theorem) parameterized by finite abelian groups, and the latter are themselves closely related to ideal class groups in number fields.  

As such, this is a natural generalization of the situation
described in the introduction to our paper: the set of geodesics of fixed length on $\SL_2(\Z) \backslash \mathbb{H}$
is parameterized by the class group of a real quadratic field. 

\subsection{Proofs.}
In order to keep notation minimal, for $g \in \G(\adele)$ we shall write $[g] = \G(\Q). g . K_f \in X_{\adele}$
for the associated double coset.  If $\mathfrak{S} \subset \G(\adele)$ is a subset, we often write simply $[\mathfrak{S}]$ for $[g]: g \in \mathfrak{S}$.

Let us recall from \cite{ELMV1}:
\begin{Proposition}[Basic correspondence]  \label{Basic bijection}
There is a canonical bijection between
\begin{enumerate}
\item [(1)] periodic $H$-orbits $\Gamma g H$ on $\Gamma \backslash G$
\item [(2)] $\Gamma$-orbits on pairs $(\T, g)$ where $\T$ is an anisotropic torus defined over $\Q$
and $g \in G/H$ is so that $g H g^{-1} = \T(\R)$. 
\end{enumerate}
\end{Proposition}

The bijection associates to $\Gamma g H$ the pair $(\T, g)$, where $\T$ is the unique $\Qq$-torus
whose real points are $g H g^{-1}$.

{\em Proof of the first statement of the theorem.} Given $(\T, g)$, clearly $\Gamma g H$
is contained in the projection to $X_{\adele}$
of the homogeneous toral set
 $Y_0 := (\T(\Q) \backslash \T(\adele)). (g,1)$. (Here $(g,1) \in \G(\adele)$ is the element that equals $g$ at the real place, and is the identity elsewhere;
in what follows we abbreviate it simply to $g$.)

Now let us show that $Y_0$ is the only such homogeneous toral set, up to $K_f$. 
Take any homogeneous toral set $Y = (\T'(\Q) \backslash \T'(\adele)) g'_{\adele}$
whose projection to $X_{\adele}$ contains $\Gamma g H = \Gamma \T(\R) g$. Therefore,
$$\T(\R) \subset \bigcup_{\delta \in \G(\Q)} \delta \T'(\adele) g'_{\adele}
K_f 
 g^{-1}.$$
It follows that there exists $\delta \in \G(\Q)$ so that
$$\T(\R)^{(0)} \subset \delta \T'(\adele) g'_{\adele} K_{f} g^{-1},$$
Therefore  there exists $t' \in \T'(\adele)$ and
$k \in K_{f}$ so that $\delta t' g'_{\adele} k g^{-1}=1$ and $\T= \delta \T'\delta^{-1}$ 
and, moreover,. 
We conclude that
$$(\T'(\Q) \backslash \T'(\adele)) g'_{\adele} =
\delta^{-1} (\T(\Q) \backslash \T(\adele)) \delta t' g'_{\adele} 
= (\T(\Q) \backslash \T(\adele)) g k^{-1}$$
where we treated these as subsets of $\G(\Q)\backslash \G(\adele)$. 

It follows that any homogeneous toral set $Y$, whose projection to $X_{\adele}$ contains $\Gamma g H$, is necessarily of the form $(\T(\Q) \backslash \T(\adele)). (g,1)$, up to modification by $K_f$. 

The equality of discriminants asserted in the Theorem is a direct consequence
of the definitions of the discriminant (see \cite[(2.2)]{ELMV1} 
for the definition for compact Cartan orbits).
For this note that the $p$-adic discriminant measures the power of $p$ in the denominator of $\iota(\mathfrak{t})$, while the discriminant at $\infty$ equals the norm of $\iota(\mathfrak{h})$
and so is constant.

{\em Proof of the second assertion of the theorem.}
Now let us observe that the packet
of the compact orbit $\Gamma gH$ consists of all compact orbits
$\Gamma \delta g H$ with $\delta \in \G(\Q) \cap K_f \T(\adelef)$ (where the intersection is taken in $\G(\adelef)$). 
Here $\T$ is the torus corresponding to $\Gamma g H$.
We shall verify that the packet of $\Gamma g H$ is parameterized by the fiber of the map
$$\T(\Q) \backslash \T(\adelef)/ K_f \cap \T(\adelef) \rightarrow \G(\Q) \backslash \G(\adelef) / K_f$$
above the identity double coset.  The finiteness assertion is then an immediate consequence of the finiteness of class numbers for
algebraic groups over number fields. 

Notations being as above, two $\delta, \delta' \in \G(\Q) \cap K_f \T(\adelef)$
define the same compact orbit $\Gamma \delta g H$ if and only if 
$\Gamma \delta \T(\Q) = \Gamma \delta' \T(\Q)$; therefore compact orbits are parameterized
by the double quotient
$$ \G(\Q) \cap K_f \backslash (\G(\Q )\cap K_f \T(\adelef))/ \T(\Q) =
K_f \backslash (K_f \G(\Q) \cap K_f \T(\adelef))/ \T(\Q) $$
but this is (after inverting) precisely what is described by the theorem. 

{\em Proof of the final assertion of the theorem.}
Let $\Gamma g H, \Gamma \delta g H$ be in the same packet. We have already verified the equality of discriminants.  To verify the equality of stabilizers, we check
that $\T(\R) \cap \Gamma = \T(\R) \cap \delta^{-1} \Gamma \delta$. 
By assumption, $\Gamma = \G(\Q) \cap K_f$. 
 Thus $\T(\R) \cap \Gamma = \G(\Q) \cap \T(\R) K_f$ and
$\T(\R) \cap \delta^{-1} \Gamma \delta = \G(\Q)  \cap \T(\R) (\delta_f^{-1} K_f \delta_f)$, 
where $\delta_f$ is the image of $\delta$ under $\G(\Q) \hookrightarrow \G(\adelef)$.
There exists $t_f \in \T(\adelef), k_f \in K_f$ so that $\delta_f = k_f t_f$; therefore, we need to prove
$$\G(\Q) \cap \T(\R) K_f= \G(\Q) \cap \T(\R) (t_f^{-1} K_f t_f).$$
An element on the left-hand side belongs to $\T(\Q)$, so automatically commutes with $t_f$, and so belongs to the right-hand side also.  Reversing this reasoning shows the equality.  \qed

\subsection{Example: packets for $\GL_n$}\label{glnpackets}
Let us explain, by way of illustration, the equivalence relation explicitly in the case of $\PGL_n(\Z) \backslash \PGL_n(\R)$. 

More precisely, we take $\G = \PGL_n$,
$\Gamma = \PGL_n(\Z)$, $H$ the diagonal subgroup, and $K_f$
the closure of $\Gamma$ in $\PGL_n(\adele_f)$. We may identify
the Lie algebra $\mathfrak{pgl}_n$ with the quotient of $n \times n$ matrices by diagonal matrices;  for $\mathfrak{g}_{\Z}$, we take then
the projection of $M_n(\Z)$ to $\mathfrak{pgl}_n$.
For $B$
we take simply the Killing form, and we take the Euclidean norm on $\pgl_{n,\R}$ to be that derived from the Hilbert-Schmidt norm on $M_n(\R)$. 
 (This differs by a factor of $2$ from that chosen in \cite[\S 4]{ELMV1}.)

As discussed in \cite{ELMV1}, and recalled previously,
compact $H$-orbits are parameterized by data $(K, L, \theta)$. 

The equivalence corresponding to {\em packets} can be described in elementary terms as follows: Declare $(K, L, \theta) \sim (K, L', \theta)$
whenever $L, L'$ are {\em locally homothetic}; here, 
we say that two lattices $L, L' \subset K$ are {\em locally homothetic}
for every prime number $p$, there exists $\lambda_p 
\in (K \otimes \Q_p)^{\times}$ so that $L_{p} = \lambda_p L'_{p}$;
here $L_{p},L'_{p}$ denote their respective closures in $K \otimes \Q_p$.

Observe that $L\sim_{loc} L'$ implies that $\order_{L}=\order_{L'}$. However, the converse is a priori not true (unless $\order_{L}$ is Gorenstein, see e.g. \cite[(6.2),(7.3)]{Bass}: for instance if $\order_{L}$ is the maximal
order). Thus, the grouping into packet refines several plausible cruder groupings,
e.g., grouping compact orbits with the same volume, or grouping
compact orbits for which the order $\order_L$ is fixed. 

Moreover, one can see at a heuristic level why the ``packet" grouping has better formal properties than the ``fixed order'' grouping. Namely
 the set of proper $\order_L$ ideals  up to homothety (ie. up to multiplication by $K^{\times}$) is not a priori a group, nor a principal homogeneous space under a group.
 By contrast, the  set of lattices up  to homothety, within the local-homothety class of $L$, {\em does} form a principal homogeneous space for a certain abelian group, the {\em Picard group $\Pic(\order_{L})$}, ie. the group of homothety classes of ideals locally homothetic to $\order_{L}$.  Of course, if $\order_L$ happens to be Gorenstein,
the packet is parameterized simply by $\Pic(\order_{L})$.

\section{Homogeneous toral sets for $\GL_n$}  \label{HsetsGLn}
We shall now consider explicitly the case of
$\G = \PGL_n$ over a global field $F$,
introducing data $\data$ which parameterizes homogeneous toral sets,
which we term global torus data. 

By localization, such data will give rise to certain {\em local} data over each place of $v$,
which we term local torus data. 
We shall explain how to compute the discriminant of the homogeneous toral
set -- in the sense of \S \ref{Hsets} --  very explicitly in terms of the local torus data.

\subsection{The local data.} \label{LTD}
{\em Local torus data} $\mathscr{D}$, over a local field $k$, consists simply of
an {\'e}tale $k$-algebra $A \subset M_n(k)$ of dimension $n$: i.e.\ a direct sum $A=\oplus_{i}K_{i}$ where the $K_{i}$ are field extensions of $k$.

Set $\g = M_n(k)/k$ identified as the Lie algebra of $\PGL_{n}$.
Let $V$ be the affine space $$V=(\wedge^{n-1} \g)^{\otimes 2}$$
To the data $\data/k$,
we associate a point in $x \in V$ via:
 \begin{equation} \label{embeddinglocal} x_{\data} =  [(f_1 \wedge \dots \wedge f_{n-1})^{\otimes 2}
(\det B(f_i, f_j))^{-1}]\end{equation}
where $f_1, \dots, f_{n-1}$ is a basis for $A/k \subset M_n(k)/k$,
and $B$ is the Killing form on $\g$.

We set $\disc(\data/k) := \|x_{\data}\|_{V}$.
Here $\|\cdot\|_{V}$ is (see \S \ref{local}) \begin{itemize}
\item[-]  the norm  on $V$ whose unit ball is the closure of $\wedge^{r} M_n(\order_k)/\order_k$, for $k$ non-archimedean;
\item[-]  the norm  on $V$ which descends from the Hilbert-Schmidt norm on $M_n(k)$ for $k$ archimedean.
\end{itemize}

Explicitely: let $f_0, f_1, \dots, f_{n-1}$ be a $k$-basis of for $A$ with $f_{0}\in k$ and
which span $\Lambda=A \cap M_n(\order_k)$ as an $\order_k$-module (when $k$ is nonarchimedean)
or which is orthonormal w.r.t. the Hilbert-Schmidt
norm on $M_n(k)$ (when $k$ is archimedean). If this is so,
one may compute $\disc(\data/k)$ by the rule:
\begin{equation} \label{explicit}\disc(\data/k) := |(2n)^{1-n} \det( \tr(f_i f_j) )^{-1}|.\end{equation}
In particular, for $k$ nonarchimedean, the discriminant of $\data/k$ 
differs, by a constant, from the discriminant of the ring $\Lambda$. 

\subsection{Proof of equivalence between \eqref{explicit} and \eqref{embeddinglocal}.} \label{Proof-equiv}

This is, in essence, proved in \cite{ELMV1}. Let us reprise it here. 
Let $f_i$ be a basis for $A$, as chosen above. 
Let $\bar{f}_i$, for $1 \leq i \leq n-1$, be the projection of $f_i$
to $\g$. 
Norms as above, $$\|(\bar{f}_1 \wedge \dots \wedge \bar{f}_{n-1})^{\otimes 2}\|_{V} = 1.$$
In fact, in the nonarchimedean case, we may extend
$\bar{f}_i \ (1 \leq i \leq n-1)$ to an $\order_k$-basis
for $M_n(\order_k)/\order_k$, which makes the result obvious.
In the archimedean case, the claim is an immediate consequence
of the fact that $f_0, \dots, f_n$ are orthonormal
w.r.t. the Hilbert-Schmidt norm. 

The Killing form on $\g$ evaluates to:
$B(\bar{f}_i, \bar{f}_j) = 2 (n \tr(f_i f_j) - \tr(f_i) \tr(f_j))$.
Therefore, the determinant
$\det B(\bar{f}_i, \bar{f}_j)_{1 \leq i,j \leq n-1}$
equals $(2n)^{n-1} \det(\tr(f_i f_j)_{0 \leq i,j \leq n-1}$.
The discussion of \S \ref{subsec:disc} 
associates to $\mathfrak{t} = A/k \subset \g$
the point
$$\iota(\mathfrak{t}) = (2n)^{1-n} \det(\tr(f_i f_j)_{0 \leq i \leq n-1})^{-1}
 (\bar{f}_1 \wedge \dots \bar{f}_{n-1})^{\otimes 2}$$
The required compatibility follows. 
\subsection{Global data} \label{GTD}
Let $F$ be a global field.
We define {\em global torus data} $\mathscr{D}$ to consist of a subfield $K \subset  M_n(F)$ and an element $g_{\adele} = (g_{\infty}, g_f)
\in \adele_K^{\times} \backslash \GL_n(\adele)$. 

To the global data we may associate  (cf. \S \ref{Hsets}):
\begin{enumerate}
\item A homogeneous toral set $Y_{\mathscr{D}} =
(\T_K(F) \backslash \T_K(\adele))g $, and the probability
 measure $\mu_{\mathscr{D}}$ on $Y_{\mathscr{D}}$. 
Here $\T_K$ is the unique subtorus of $\PGL_n$
with Lie algebra $K/F$. 
\item A (global) discriminant $\disc(\mathscr{D}):=\disc(Y_{\data})$, depending on $K, g_{\adele}$. 
\end{enumerate}

The global discriminant $\disc(\data)$ can be computed in terms of our discussion above. Indeed, the global torus data $\data = (K,g)$ also gives
rise to a collection of local torus data $(\data_{v})_{v}$: for every place $v$, 
$\data_v$ consists of the subalgebra $$A_{v}= g_v^{-1} (K \otimes F_v) g_v \subset M_n(F_v)$$ and it follows from the definitions that
$$\disc(\data)= \prod_{v} \disc(\data_v).$$

\section{The local building} \label{local}
In this section, we are going to recall the basic theory of the
building attached to the general linear group, over a local field $k$. 
We will follow the beautiful old ideas
of Goldman and Iwahori, \cite{GoldmanIwahori}, interpreting
this by norms on a $k$-vector space.

\subsection{Notation concerning local fields} \label{NormCon}
We denote by $k$ a local field.  Let us normalize
once and for all an absolute value on it. If $k=\R$ or $\Cc$, $|\cdot|$ denote the usual absolute value and
 if $k$ is non-archimedean, we normalize $|\cdot|$ to be the module of $k$: $|\cdot|=q^{-v_{\pi}(\cdot)}$, where 
$q$ is the cardinality of the residue field and $\pi\in k$ any uniformizer, and $v_{\pi}(\cdot)$,
 the corresponding valuation. 

\subsection{Definition of the building by norms.}
Let $k$ be a local field and $V$ vector space over $k$ of finite dimension.
A norm on $V$ is a function $N:V\ra\Rr^+$ into the non-negative reals
that satisfies 
\begin{gather*}N(v)=0\Leftrightarrow v=0_{V},\ \ 
N(\lambda x) = |\lambda| N(x),\ \lambda\in k\\
N(x+y) \leq \begin{cases}N(x)+N(y)& \mbox{(if $k$ is archimedean)}\\
 \max(N(x), N(y))& \mbox{(if $k$ is nonarchimedean). }
\end{cases}
\end{gather*}
For $N$ a norm, we denote its homothety class by $[N]$:
$$[N]=\{ \mu N,\ \mu \in \R_{>0}\}$$

If $k$ is real (respectively complex),
we call a norm on $V$ {\em good} 
if it is quadratic (respectively Hermitian). If $k$ is nonarchimedean,
we shall refer to any norm as {\em good}. 

 We let $\B(V)$ and $\bB(V)$ be the building of $\GL(V)$
and $\PGL(V)$ respectively; specifically
$\B(V)$ is the set of good norms on $V$, and $\bB(V)$ the set of such norms up to homothety. 

\subsection{Action of the group}
The group $\GL(V)$ acts transitively on $\B(V)$ 
via the rule $$g.N(x) = N(xg).$$
This induces a transtive action of $\PGL(V)$ on $\bB(V)$.

If $k$ is archimedean, 
the stabilizer of a good norm is a maximal
compact subgroup, and any (good) norm is determined by its unit ball.

If $k$ is nonarchimedean, neither of these are true:
consider, for instance, the norm 
on $\Q_p^2$ given by $(x,y) \mapsto \max(|x|, p^{-1/2} |y|)$. 
We say that a norm is {\em standard} if is satisfies
$$N(x) = \inf\{|\lambda|: \lambda \in k, N(x) \leq |\lambda|\}.$$
A standard norm is determined by its unit ball, and its stabilizer
is a maximal compact subgroup of $\GL(V)$. Standard norms, are, equivalently ($q$ the cardinality of the residue field of $k$):
\begin{enumerate}
\item Those which take values in $q^{\Z}$;
\item  Those which look like $x=(x_1, \dots, x_n) \mapsto
\max_i |x_i|$ in suitable coordinates;
\item Those that correspond to (special) vertices of the building. 
\end{enumerate} 
The action of $\GL(V)$ preserves standard norms.
\subsection{Direct sums. Apartments.} \label{apartment}
Given norms $N_V$ on $V$ and $N_W$ on $W$, they determine a norm
$N_V \oplus N_W$ on $V \oplus W$, defined by $\sqrt{N_V^2+N_W^2}$ if $k$ is archimedean,
and $\max(N_V(v), N_W(w))$ if $k$ is nonarchimedean.

Any splitting of $V$ into one-dimensional subspaces determines an {\em apartment}.
This consists of all norms that are direct sums of norms on the one-dimensional subspaces. 
Any two norms belong to an apartment.  Apartments are in bijection with
split tori within $\GL(V)$: a splitting of $V$ into one-dimensional spaces
determines a split torus, namely, those automorphisms of $V$ preserving each one-dimensional space. 

If $H$ is a split torus 
with co-character lattice $X_{\star}(H)$, then the vector space $\h := X_{\star}(H) \otimes \R$
acts simply transitively on the apartment.
 It suffices to explicate this when $V$ is one-dimensional;
in that case, the action of the tautological character $\G_m \rightarrow \GL(V)=H$ 
on the set of norms is multiplication by cardinality of the residue field if $k$ is nonarchimedean, and multiplication by (e.g.) $e=2.718$, when $k$ is archimedean.

In explicit terms, we can phrase this as follows:
the apartment in the building of $k^n$, corresponding to the diagonal torus in $\GL(n,k)$, consists of all norms of the following form,
:
$$N(x_1, \dots, x_n) = \begin{cases} \max_i q^{t_i}|x_i|, 
\mathrm{nonarchimedean},
\ q:= \mbox{size of residue field.} \\ \left( \sum_i  e^{2 t_i} |x_i|^2 \right)^{1/2}, \mathrm{archimedean}.,\end{cases}$$
for  some $(t_1, \dots, t_n) \in \R^n$.
Therefore, this apartment is parameterized by the affine space $\R^n$.

\subsection{The canonical norm on an algebra.} \label{CNA}
Let $A$ be a (finite dimensional) {\'e}tale algebra over $k$, i.e.\ $A= \oplus_{i} K_i$ is a direct sum of field extensions of $k$, 
we equip it with a norm which we shall call the {\em canonical} norm.

\begin{itemize}
\item[-]If $k$ is nonarchimedean: let $\order_{A}=\oplus_{i}\order_{K_{i}}$
denote the maximal compact subring of $A$:
\begin{equation}\label{canonicalnonarchimedeandef} N_{A}(t) = \inf \{ |\lambda|,\ \lambda\in {k}: t \in \lambda \order_A.\}\end{equation} 
That norm is standard and has unit ball $\order_A$. Moreover, for $t$ decomposing as $t=(t_{1},\dots,t_{i},\dots)$, $t_{i}\in K_{i}$,
one has
\begin{equation}\label{maxNA}
N_{A}(t)=\max_{i}N_{K_{i}}(t_{i}).
\end{equation}
Here $N_{K_{i}}$ denote again the standard norm of the $k$-algebra $K_{i}$.

\item[-]If $k$ is non-archimedean, $A= \oplus K_i$, for certain subfield $K_i \subset \C$.  We define
$$N_{A}(\sum_{i} x_i)=
\left(\sum_i |x_i|^2 \right)^{1/2}.$$
To keep our notations consistent with the non-archimedean case, we {\em define} $\order_{A}$ to be the unit ball of $N_{A}$:
\begin{equation} \label{orderdef} \order_{A}=\{x\in A,\ N_{A}(x)\leq 1\}.\end{equation}
\end{itemize}

When $A=k^n$, we denote the canonical norm by $N_{0}$.

\subsection{The metric and operator norms.}
We may equip $\B(V)$ with a $\GL(V)$-invariant metric 
by using the notion of operator norm:

 If $N_1, N_2$ are any two norms,  
 we let $\exp(\dist(N_1, N_2))$ be the smallest constant $\alpha\geq 1$ such that $N_2$ satisfies $\alpha^{-1} N_1 \leq N_2 \leq \alpha N_1$.

Given any two norms $N_1, N_2$, there exists an apartment that contains them both; thus, to understand $\dist$, it suffices to understand it on each apartment. 
In \S \ref{apartment}, we explicitly parameterized each apartment
by an affine space $(t_1, \dots, t_n) \in \R^n$. In terms of that parameterization, 
$$\dist((t_1, \dots, t_n), (t_1', \dots, t_n')) := \log(q) \max_{i}
|t_i - t_i'|$$
Here we understand $q = e$ for $k$ archimedean. Therefore, this amounts
to an $L^{\infty}$-metric on each apartment.

   We equip $\bB(V)$ with the quotient metric\footnote{Recall that if $X$ is a metric space and $G$ acts by isometries on $X$, we may define the metric on $X/G$ via
 $d(x_1, x_2) = \inf_{g \in G} d(x_1 g, x_2)$.}.
   In particular, if $N_1, N_2$ are two norms and $[N_1], [N_2]$ the corresponding 
   elements in $\bB(V)$, then, for any $v_1, v_2 \in V$:
   \begin{equation} \label{note} \dist([N_1], [N_2]) \geq  \frac12\log \frac{N_1(v_1) N_2(v_2)}{N_2 (v_1) N_1(v_2)}.\end{equation}
Indeed, one may define $\dist$ as the supremum of the quantities appearing on the right-hand side.

\subsection{Harish-Chandra spherical function} 
We shall make use of the Harish-Chandra spherical function on 
$\GL(V)$.  It is defined with reference to a maximal compact subgroup $K \subset \GL(V)$, which we take to be the stabilizer of a standard norm $N$. 

Fix an apartment containing $N$, and let $H \subset \GL(V)$ be the split torus
corresponding to this apartment.  Let $B \supset H$ be Borel containing $H$, with unipotent radical $U$ corresponding
to all positive roots. We have a decomposition $\GL(V) = U H K$.
Let $\mathrm{H}: \GL(V) \rightarrow H$ be the projection according to this decomposition
and let $\rho: H \rightarrow \R_{+}$ be defined by
$$\rho: a \mapsto \prod_{\alpha\in\Phi^+}|\alpha(a)|^{1/2}
$$ be the ``half-sum of positive roots'' character w.r.t $B$. 

The Harish-Chandra spherical function is defined as:
$$\Xi(g) := \int_{k \in K} \mathrm{H}(k g)^{\rho} dk,$$
where the measure on $K$ is the Haar measure with total volume $1$. 
We will be needing the following bound: for any $\alpha < 1$,
\begin{equation} \label{HCS} \Xi(g) \ll_\alpha \exp(-\alpha. \dist([g N], [N]))\end{equation}

Indeed, it suffices to prove \eqref{HCS} when $gN$ belongs to a fixed apartment containing $N$. Identifying this with an affine space, with $N$ as origin,
the point $gN$ has coordinates $(t_1, \dots, t_n)$; without
loss of generality, we may assume that $t_1 \leq t_2 \leq \dots t_n$, and (cf. \cite[Proposition 7.15]{Knapp} for real semisimple Lie groups) 
$$\Xi(g)  \ll_\alpha \bigl(q^{- \frac{1}{2}\sum_{i<j} t_i - t_j}\bigr)^{-\alpha}$$
for any $\alpha<1$. On the other hand, by \refs{note} $$\dist([gN], [N]) = \frac{1}{2}\log q. (t_{n}-t_{1})
,$$ whence our conclusion. 

\subsection{Action of invertible linear maps} \label{normconventions}
If $\iota: V \rightarrow W$ is an invertible map between 
vector spaces, which we understand as acting on the right (i.e.\ $v.i$ or $v^\iota \in W$), and $N$
is a norm on $W$, we denote by $\iota N$ the norm $v \mapsto N(v.\iota)$, a norm on $V$.

\section{Notation} \label{Notation}

In this section, we set up some fairly standard notation concerning number fields. 
We set up local notation first, and then global notation. In \S \ref{toruslocal} we shall
use only the local notation; in the rest of the paper, we make use of the global notation.

\subsection{Local notation and normalizations}$\phantom{.}$ \label{phantom}
Let $k$ be a local field and $A \subset M_n(k)$ local torus data, with $[A:k]=n$. 
The absolute value on $k$ is normalized as in \S \ref{NormCon}. 

We denote by $|\cdot|_A$ the ``module'' of $A$, i.e., the factor by which the map
$y \mapsto yx$ multiplies Haar measure on $A$ if $x\in A^\times$ and $|x|_{A}=0$ if $x\in A-A^\times$. Consequently, writing $A=\oplus_{i}K_{i}$, 
$$|x|_{A}=\prod_{i}|x_{i}|_{K_{i}},\ x=(\dots,x_{i},\dots),\ x_{i}\in K_{i};$$
in particular, for\footnote{At this point we must observe a small ugliness of notation:
for $x \in k$, 
the ``module'' $|x|_k$ coincides with our normalization of $|x|$ from
\S \ref{NormCon}
if $k \neq \C$.  If $k=\C$, however, $|x|_{\C} = |x|^2$. This unfortunate
notational clash seems somewhat unavoidable, for the module of $\C$ does not coincide with what is usually termed the absolute value. } $x\in k \subset A$, $|x|_{A}=|x|_k^n$.
 Observe also that the module $|x|_{A}$ coincide with $ |\det(x)|$ when we view $x$ as an element of $M_{n}(k)$.

We fix an additive character $e: k \rightarrow \C$;
it induces the additive character on $A$
$$a\ra e_{A}(a) := e(\tr_{A/k}(a))=e(\tr(a)),\ a\in A.$$ 
Observe that for $a\in A$, the $A/k$-trace coincide with the restriction to $A\subset M_{n}(k)$ to the matrix trace, thus there is no ambiguity in refereing to the trace.

  We fix Haar measure $dx, da$ 
on $k, A$; for definiteness, we normalize them to be self-dual w.r.t the characters $e(\cdot)$ and $e(\mathrm{tr}(\cdot))$ respectively. Sometimes we will write $d_{k}x$ and $d_{A}x$ to emphasize 
the measures on $k$ and $A$ respectively. 
We will often write $\vol_k$ or $\vol_A$ for volume of a set with respect to these measures.

Even though we have normalized $\vol_k$ and $\vol_A$ to be self-dual, it is occasionally more conceptually clear and helpful -- for instance, when working with Fourier transforms -- to introduce a separate notation for the dual measures.  Thus,we shall denote by $\widehat{\vol_A}$ the Haar measure dual to $\vol_A$, w.r.t. the character $e_{A}$; $\widehat{\vol_k}$ is defined similarly .   Our normalizations are so that $\widehat{\vol} = \vol$, but
we try to keep the notions conceptually separate.

We normalize multiplicative Haar measure $d_{k}^{\times}x, d_{A}^{\times}x$
on $k^{\times}$ and $A^{\times}$, respectively, by the rules
$$d_{k}^{\times}x = \zeta_k(1) |x|_k^{-1}d_{k}x, d_{A}^{\times}a =  \zeta_A(1) |x|^{-1}_{A}d_{A}a.$$ 
(See \S \ref{discdef} below, for a recollection of the definition of $\zeta$). 

These normalizations have the following effect:
for $k$ nonarchimedean,
\begin{equation} \label{mult}
\vol_{A^\times}(\order_{A}^{\times}) = \vol_{A}(\order_A), \vol_{k^\times}(\order_k^{\times})
= \vol(\order_k). \end{equation}
For $k$ archimedean, we shall not define $\order_A^{\times}$ or $\order_k^{\times}$; however,  it will be convenient (to uniformize notations) to {\em define} their volumes $\vol_A(\order_A^{\times}):=\vol_A(\order_{A})$
and $\vol_k(\order_k^{\times}):=\vol_k(\order_{k})$ so that the equality \eqref{mult} remains valid.
Recall that $\order_A, \order_k$ in the archimedean case are defined by the convention \eqref{orderdef}.

\subsection{Local $\zeta$-functions} \label{discdef}
The local $\zeta$-function of the field $k$ and of $A$ are denoted
$\zeta_k(s), \zeta_{\algebra}(s)$ ($s\in \C$): for $A=\oplus_{i}K_{i}$,
$$\zeta_{A}(s)=\prod_{i}\zeta_{K_{i}}(s).$$
We recall that the local $\zeta$-function of a local field $k$ is defined by $\zeta_k(s) = \pi^{-s/2} \Gamma(s/2)$ if $k = \R$,
$2 (2\pi)^{-s} \Gamma(s)$ if $k= \C$, and finally
$\zeta_k(s) = (1-q^{-s})^{-1}$, where $q$ is the size of the residue field,
if $k$ is nonarchimedean.

More generally for $\psi$ a character of $A^\times$ ($\psi=(\dots,\psi_{i},\dots)$, $\psi_{i}$ a character of $K_{i}^\times$), we denote by $L(A,\psi)$ the local $L$-function of $\psi$:
$$L(A,\psi)=\prod_{i}L(K_{i},\psi_{i}).$$
See \cite[Chapter 3]{Bump} for definition and discussion. 

For $s\in\C$, we will also write
$$L(A,\psi,s):=L(A,\psi|\cdot|_{A}^s).$$
In particular $\zeta_{A}(s)=L(A,|\cdot|^s_{A})$.

If $k$ is nonarchimedean, we attach to $\psi$ a discriminant $\disc(\psi)$. 
This may be defined directly as follows:
we write $A = \oplus K_i$ and $\psi=(\dots,\psi_{i},\dots)$; for each $i$, let $t_i$ be the largest integer so that $\psi_{i}$ is trivial
on $1 + \mathfrak{q}_{K_i}^{t_i}$;
here $\mathfrak{q}_v$ is the prime ideal of the ring of integers in $K_i$. 
Then $\disc(\psi) := \prod_{i} q_i^{t_i}$, where $q_i$ is the residue field size of $K_i$. 

For $k$ archimedean, we set by definition $\disc(\psi) \equiv 1$.

\subsection{Number fields} \label{ssnf}
We now pass to a global setting. 

Let $F$ be a number field. 
 We denote by $\adele$ and $\adele_{F,f}$ the ring of adeles and of finite adeles, respectively.

We will work with global data $\data$, as in \S \ref{GTD},
which will consist of:
 $K \subset M_n(F)$ and
 $g_{\data} \in \adele_K^{\times} \backslash \GL_n(\adele)$. 
We shall fix an identification $\iota: K \rightarrow F^n$ of right $K$-modules; this means
that $(a b)^{\iota} = a^{\iota} . b$, where, on the right hand side, we understand $ b \in M_n(F)$. 

We obtain from this data an embedding of the $F$-torus $$\T_K := \mathrm{Res}_{K/F} \mathbb{G}_m / \mathbb{G}_m \hookrightarrow
\PGL_{n,F}$$
Here $\PGL_{n,F}$ denotes the algebraic group $\PGL_n$ over the field $F$. 

 We will use the letter $v$ for a place of $F$ and $w$ for a place of $K$. If $v$ is a place of $F$,
we denote by $F_v$ the completion of $F$ at $v$, and $K_v := K \otimes_F F_v$.

By localization, the global data $\data$ gives rise to local data
$\data_v$ (in the sense of \S \ref{LTD}) for each place $v$ of $F$, 
i.e.\ we take $A_{v} = g_{\data,v}^{-1} K_v g_{\data,v} \subset M_n(F_v)$. We will 
write $A$ instead of $A_{v}$ when the dependence on $v$ is clear.

Let us note that the map 
\begin{equation}\label{iotavdef}\iota_v :
a \in A_v \mapsto  (g_{\data,v} a g_{\data,v}^{-1})^{\iota} . g_{\data,v} \end{equation}
from $A_v$ to $F_v^n$ is then an identification for the $A_{v}$-module structures.

 \subsection{Adeles, ideles and their characters.}
There is a natural norm map, the module, $\adele^{\times} \stackrel{|\cdot|}
{\rightarrow} \R_{>0}$. We write it $|x|_{\adele}$
or sometimes simply $|x|$; this will cause no confusion so long
as it is clear that the variable $x$ belongs to $\adele^{\times}$. 
We denote by $\adele^{(1)}$ the kernel of the norm map,
 and similarly for $\adele_K$.

Let $\Char$ resp. $\CharK$ denote the group of homomorphisms from $\adele^{\times}/F^{\times}$ resp.
$\adele_K^{\times}/K^{\times}$ to $\C^{\times}$. 
  For $\psi \in \CharK$ we shall denote by $\psi|_F$ the restriction of $\psi$ to $\adele^{\times}/F^{\times}$.

The group $\C$ is identified with the connected component of $\Char$, via
identifying $s \in \C$ with the character $x \mapsto |x|_{\adele}^s$. 
Given $\chi \in \Char$, we set $\chi_s(x) := \chi(x) |x|_{\adele}^s$. 
This $\C$-action on $\Char$ gives $\Char$ the structure of a complex manifold. 

For any $\chi \in \Char$, there exists unique $s \in \R$
so that $|\chi(x)| = |x|_{\adele}^s$. We shall denote this $s$ by $\Re \chi$, the ``real part'' of $\chi$. Thus $\Re \chi =0$ if and only if $\chi$ is unitary. 

Finally, there is a natural map $\R_{>0} \rightarrow C_{\Q}$ (inclusion
at the infinite place). Thus there is also a map $\R_{>0} \rightarrow  C_K, C_F$.
We say that a character of $C_K$ or $C_F$ is {\em normalized} if its
pullback to $\R_{>0}$ is trivial.

If $\omega: F_v^\times \rightarrow \C^{\times}$ is a multiplicative character,
we denote by $L(F_v, \omega, s)$ the corresponding $L$-factor,
and by $L(F_v, \omega)$ its value when $s=0$. 
In particular, when $\omega = |x|_v^s$, we get the local $\zeta$-factor:
$L(F_v, \omega) = \zeta_{F_v}(s)$. 
Corresponding definitions also hold for $K$.
If $v$ is a place of $F$, we will write
$\zeta_{K,v} := \prod_{w|v} \zeta_{K,w}.$

\subsection{Discriminants}
Suppose $\psi \in \CharK$. 
  For $v$ any place of $F$, we have discriminants:
 $$\disc_v(F), \disc_v(K/F), \disc_v(K), \disc_v(\data_v), \disc_v(\psi|_{K_v})$$

Namely, $\disc_v(F)$ is the discriminant of $F_v/\Q_p$ (where $p$ is the prime of $\Q$ below $v$) and $\disc_v(K)$ the discriminant of $K_v/\Q_p$.  We set
  $$\disc_v(K/F)  = \disc_v(K) \disc_v(F)^{-[K:F]}.$$  By convention, we shall understand $\disc_v(F) = \disc_w(K) = 1$ if $v$ or $w$ are archimedean. 
 $\disc_v(\data_v)$ is as in \S \ref{HsetsGLn} and $\disc_v(\psi|_{K_v})$ is defined in \S \ref{discdef}. 

For any of the objects above, we set $\disc(\dots) = \prod_{v} \disc_v(\dots)$. 
We note that $\disc(\data) \gg_{F, n} \disc(K)$;
this will follow from Lemma \ref{dvexpl} and the
fact that $\disc_v(\data)$ is bounded from below at each archimedean place. 

\subsection{Measure normalizations}
Let $e_{\Q}: \adele_{\Q} /\Q \rightarrow \C$
be the unique character whose restriction to $\R$ is $x \mapsto e^{2 \pi i x}$.  Set $e_F = e_{\Q} \circ \mathrm{tr}_{F/\Q}$. 
We then normalize local measures according to the prescription of \S \ref{phantom}
with $k=F_v, A = K_v$.

 Let us explicate this to be precise. 

We choose for each $v$ the measure on $F_v$ that is self-dual w.r.t. $e_F$.
The product of these measures, then, assigns volume $1$ to $\adele/F$. 
We define
 a measure on $F_v^{\times}$ by $d^{\times}x :=\zeta_{F,v}(1) \frac{dx}{|x|_v}$. 
We make the corresponding definitions for $K$, replacing the character $e_F$
by the character $e_K := e_F \circ \mathrm{tr}_{K/F}$. 
Taking the product of these measures yields measures on $\adele, \adele_K, \adele^{\times}, \adele_K^{\times}$. 
This fixes, in particular, a quotient measure on $\T_K(\adele) = \adele_K^{\times}/\adele^{\times}$. 

We obtain a measure on $A_v$ through the identification
$x \mapsto g_{\data,v}^{-1} x g_{\data,v}$ from $K_v$ to $A_v$.

We will often denote by $d_{K}x$
the measure on $K_v$ and by $d_{F}x$
the measure on $F_v$ or $F_v^n$. 

With these definitions, it is not difficult to verify that
for finite places $v,w$:
\begin{eqnarray}\label{measure}\int_{ \order_{K,w}^{\times}} d_K^{\times}x = \disc_w(K)^{-1/2}, 
\int_{\order_{F,v}^{\times}} d_{F}^{\times} x = \disc_v(F)^{-1/2} \\
\label{measure2} \int_{ \order_{K,w}} d_{K}x = \disc_w(K)^{-1/2}, 
\int_{\order_{F,v}} d_F x = \disc_v(F)^{-1/2} \end{eqnarray}

Moreover, \eqref{measure2} remains valid for $v$ archimedean, if we replace
equality by $\asymp$, and we interpret $\order_{K,w}$ resp. $\order_{F,v}$
as the unit balls for the canonical norms associated to $K_w$ resp. $F_v$.

\section{Local  theory of torus orbits} \label{toruslocal}
In this section, we are going to explicate certain estimates over a local field,
which are what are needed for ``local subconvexity'' discussed in \S \ref{outline}.  Indeed, the main result of this section,  Proposition \ref{localbounds},
is designed precisely to bound the (general version of the) quantities
``$I_{\order}$'' that occur in \eqref{HeckeFormula} and \eqref{localsubconvex}.

  Our methodology is quite general (i.e., would apply to estimates for more general period integrals of automorphic forms) and is inspired, in part, by the paper of Clozel and Ullmo \cite{CU}.  

\subsection{Explanation in classical terms}
A simple case of our estimates is the following result:

Let $Q$ be a positive definite quadratic form on $\R^3$. Consider:
\begin{equation}\label{fdef} I(Q) = \frac{ \int_{Q(x,y,z) \leq 1} |xyz|^{-1/2} dx \, dy \, dz.}{\vol(\x\in\R^3:\ Q(\x) \leq 1)^{1/2}}
\end{equation}

It will transpire  -- this will follow from Lemma \ref{lem:int} that we prove later -- that the quantity $I_{\order}$, defined in \eqref{HeckeFormula}, will be be bounded by products of integrals like \eqref{fdef}
and nonarchimedean analogues.
(More precisely, in the notation of \eqref{HeckeFormula}: if the function $f$ is invariant by the rotation group of $Q$, then we will be able to bound $I_{\order}(f,s)$ in terms of \eqref{fdef}; the general case will reduce to this). 
Our goal will be to get good bounds for $I(Q)$.
Evidently, the function $I$ is invariant under coordinate dilations,
and thus our bounds should also be so.

\subsection{The generalization to local fields.}
The general setting we consider will be: we replace $\R, \R^3$
by a local field $k$ and an etale $k$-algebra $A$ of dimension $n$ (say). 

For an arbitrary norm $N$ on $A$, we consider the integral
\begin{equation} I(N) :=\frac { \int_{x \in A: N(x) \leq 1}
  |x|_{A}^{-1/2} d_{A}x}{\vol(x \in A: N(x) \leq 1)^{1/2} }
\end{equation}
where $d_{A}x$ denote a Haar measure on $A$ and $|x|_{A}$ denote the ``module'' of $A$ (the factor by which $d_{A}x$
 is transformed under $y\mapsto xy$).
We consider the variation of $I(N)$ with $N$. Again, $I(N)$ is invariant by scaling and so 
defines a function on the building of $\PGL(A)$.\footnote{Throughout this paper,
we use $\GL(A)$ to denote $k$-linear automorphisms of $A$, thought of as a $k$-vector space; similarly, $\PGL(A)$.}
We shall show that:
\begin{enumerate}
\item $I(N)$ decays exponentially fast with the distance of $N$
to a certain subspace in the building, viz. the $A^{\times}$-orbit of the norm $N_A$;
\item The distance to this subspace measures the discriminant
of local torus data:

 Choose an identification $\iota: A\rightarrow k^n$
which carries the unit ball for $N(x)$, to the unit ball of the
standard norm on $k^n$. The identification $\iota$, together
with the action of $A$ on itself by multiplication, can be used
to embed $A \subset M_n(k)$. Thus, we have specified local torus data. 
Roughly speaking (Lemma \ref{Deltabound}), the discriminant
of this local torus orbit, is a measure of the distance
of $N$ to our distinguished subspace of the building. 
\end{enumerate}

Our final result is presented in Proposition \ref{localbounds};
the reader may find it helpful to interpret it in terms of the language above,
in order to better absorb its content.

\subsection{Local torus data} \label{subsec:Disc2}

In the rest of this section, we fix local torus data $\data$
consisting of $A \subset M_n(k)$; we shall follow the notations of \S \ref{Notation}. 
Throughout
this section, we allow the notation
 $\ll$ and $O(\cdot)$ to indicate an implicit constant that
remains bounded, if $k$ is restricted to be of bounded degree over $\R$ or $\Q_p$.

By the inclusion, $A \subset M_n(k)$, $k^n$ is a right $A$-module. For the rest of this section, we fix an identification of  right $A$-modules
$\iota: A \rightarrow k^n$; in other words,
\begin{equation}\label{modulemap}(ab)^{\iota} = a^{\iota} . b, a \in A, b \in A\end{equation}
where, on the right hand side, we regard $b\in M_n(k)$. 
The identification $\iota$ is unique up to multiplication by elements of $A^{\times}$.  We write $\iota$ on the right
to be as consistent as possible with our other notations.

In this setting, we have described (cf. \S \ref{CNA}) two norms: the canonical norm on the algebra $A$, to be denoted $N_A$,
 and the canonical norm on the algebra $k^n$, to be denoted $N_{0}$. We have denoted their unit balls by $\order_{A}$ and $\order_{k^n}$ respectively. These unit balls coincide with the maximal
 compact subgroups of these algebras in the non-archimedean case.
  
 We are going to introduce an element $h \in \GL_n(k)$ with quantifies the relation between these norms. 
Choose $h \in \GL_n(k)$ so that:
 \begin{equation} \label{hdef} h N_0(x) (:= N_0(xh)) = N_{\algebra}(x^{\iota^{-1}}), \ \ x \in k^n.\end{equation}  This is possible because, by choice, the norm $N_{\algebra}$
corresponds to a vertex of the building of $\GL(\algebra)$, i.e.\ takes values in $q^{\Z}$ in the nonarchimedean case. 
Observe that choice of $h$ depends on the choice of $\iota$;
given $\iota$, the quantity $|\det h|$ is uniquely determined.

This definition implies that  ${\order_A}^{\iota} =  \order_{k^n} h^{-1}$; thus
$$\vol_A(( \order_{k^n})^{\iota^{-1}}) = \vol_A(\order_A) |\det h|,$$ or:
\begin{equation} \label{hdef2} \frac{\iota_* \vol_A}{\vol_{k^n}} = \frac{\vol_A(\order_A)}{\vol_{k^n}(\order_{k^n})}
|\det h|.\end{equation}

Similarly, we define $h_A \in \GL(A)$ by the following rule:
$$ (y h_A)^{\iota} = y^{\iota} h, y \in A.$$
This means that 
\begin{equation}\label{itmeans}
h_{A}\iota N_0= N_A=\iota h N_{0},
\end{equation}
 and so
$|\det h_A|_k = |\det h|_k$.

With these conventions, and those of \S \ref{normconventions}, 
\begin{equation} \label{NormCompat}t \iota^{-1} N_A = \iota^{-1} t N_A, (t \in A^{\times})\end{equation}
Let us be completely explicit, because $t$ is acting
in two different ways on the two sides of this equation. 
 According to the conventions
set out in \S \ref{normconventions}, the left-hand
side is the norm on $k^n$ defined by $\mathbf{x} \mapsto
N_A(\mathbf{x} t \iota^{-1})$; in particular, $t$ is acting as an endomorphism of $k^n$. The right-hand side is the norm on $k^n$ defined
by $\mathbf{x} \mapsto N_A(\mathbf{x} \iota^{-1} t)$; here $t$
is acting by right multiplication on $A$. 
The coincidence of the two sides follows from \eqref{modulemap}. 
\subsection{Discriminant vs. discriminant} \label{subsec:thediscriminant}
We have two notions of ``discriminant" attached to the local data $A\subset M_{n}(k)$:

On the one hand, we have the {\em absolute} discriminant, $\disc(\algebra)$, of the $k$-algebra $A$:

\begin{itemize}
\item[-]If $k$ is non-archimedean, it is  given by
$$\disc(\algebra) = [\order_{A}^*:\order_{A}]=\frac{\vol_{A}(\order_{A}^*)}{\vol_{A}(\order_{A})},$$
 Here $\order_A^{*}$ denote the dual lattice of $\order_{A}$ in $A$
$$\order_A^{*} := \{a \in A: |\tr_{A/k}(a \order_A)| \leq 1\}\supset \order_{A}.$$ 
\item[-]If $k$ is archimedean, we set $\disc(\algebra)=1$.
\end{itemize}
 In particular in either case, we have
\begin{equation} \label{volnorm} \disc(\algebra)  \asymp_{n}
\frac{\vol_{A}(\order_A^*)}{\vol_{A}(\order_A)}.\end{equation}

On the other hand, in our previous discussion \S \ref{HsetsGLn},
we have defined a notion of discriminant $\disc(\data)$ which is {\em relative} to the embedding $A\subset M_{n}(k)$.
 We shall presently compare the two notions and for this we shall interpret $\disc(\data)$ in more geometric terms.

Let $\Lambda$ be the set of $x \in A$ with operator norm $\leq 1$ with respect to the norm
$N_0$; here we regard $k^n$ as an $A$-module via the {\em right} multiplication
of $A \subset M_n(k)$.  If $k$ is nonarchimedean, $\Lambda$ is an {\em order} (and thus is contained in $\order_{A}$); indeed, $\Lambda = A \cap M_n(\order_k) $. 
Let $\Lambda^*$ be the dual to $\Lambda$,
\begin{equation}\label{dualdef}\Lambda^* = \{ y \in A:  | \tr_{A/k}(y \Lambda) | \leq 1\}.\end{equation}

\begin{Lemma} \label{dvexpl} We have (compare with \refs{volnorm})
$$\disc(\data) \asymp_n \frac{\vol_{A}(\Lambda^*)}{\vol_{A}(\Lambda)}$$
Moreover, $\asymp_{n}$ may be replaced by {\em equality} if $k$ is nonarchimedean,
of residue characteristic exceeding $n$. 
\end{Lemma}  
\proof

The definition of $\disc(\data)$ is explicated in \eqref{embeddinglocal}. Choose, first of all,
a $k$-basis $f_0 = 1, f_1, f_2, \dots, f_{n-1}$ for $A$ so that
the unit cube $C$ on basis $f_i$, i.e.
\begin{equation}\label{unitcube}\sum_{i} \lambda_i f_i,\
|\lambda_i| \leq 1,\end{equation} is equal to $\Lambda$ (nonarchimedean case)
and comparable to $\Lambda$ as a convex body (archimedean case). 

If $\bar{f}_i$ denotes the projection of $f_i$ to $A/k$, then
$\det B(\bar{f_i}, \bar{f_j})$ and $\det(\tr(f_i f_j))$
differ by $(2n)^{n-1}$ (see \cite[4.1.3]{ELMV1} or \S \ref{Proof-equiv}). 
Let $f_0^*, f_1^*, \dots, f_{n-1}^* \in A$ be the dual basis to the $f_i$, that is to say: $\tr(f_i f_j^*) = \delta_{ij}$. 

Then the unit cube $C^*$ on basis $f_i^*$ 
equals $\Lambda^*$ (nonarchimedean case) and
is comparable with $\Lambda^*$ as a convex body (archimedean case). 

On the other hand,  $\vol(C)/\vol(C^*) = \det \tr(f_i f_j)$. 
Our claimed conclusion follows. 
\qed

\begin{Lemma} \label{mindisc}
Suppose $k$ is non-archimedean, of the residue characteristic greater than $n$ and that  $\disc(\data) = 1$.
  Then $\algebra$ is unramified (i.e. a sum of unramified field extensions of $k$)  and $\order_{k^n}$ is stable under $\order_{A}$.
  In particular, $(\order_{k^n})^{\iota^{-1}} = \lambda \order_A$, for some $\lambda \in A^{\times}$. 
\end{Lemma}
\proof Since $k$ is non-archimedean, we have the chain of inclusions
\begin{equation}\label{inclusion}\Lambda \subset \order_{A} \subset \order_{A}^* \subset \Lambda^*;
\end{equation}
our assumption and the prior Lemma shows that equality holds. 

That is to say, $\order_{A}$ is self-dual w.r.t. the trace form; this implies that $\algebra$ is unramified over $k$.   

As for the latter statement, $\order_{k^n}$ is stable by $\Lambda$ by definition,
and therefore by $\order_{A}$ since $\Lambda = \order_{A}$. 
It is equivalent to say that $(\order_{k^n})^{\iota^{-1}}$ is stable under $\order_A$; whence the final statement. 
\qed

Therefore, the quantity $\disc(\data)$ measures the {\em distance of the data $\data$
from the most pleasant situation.}

\begin{Lemma} \label{localvolumes}
If $k$ nonarchimedean, let $\Lambda^{\times}$
be the units of the order $\Lambda$. 
Then $$\frac{\vol_{A}(\Lambda^{\times})}{\vol_{A}(\order_\algebra^{\times})} \geq c(n) \max(1-n/q, q^{-n})  \left( \frac{\disc(\data)}{\disc(A)} \right)^{-1/2}.$$  
Here $c(n)=1$ if the residue characteristic of $k$ exceeds $n$. 
\end{Lemma}

\proof We have by \refs{inclusion}
$$[\Lambda^*:\order_{A}^*]=\frac{\vol_{A}(\Lambda^*)}{\vol_{A}(\order_{A}^*)}=[\order_{A}:\Lambda]=\frac{\vol_{A}(\order_{A})}{\vol_{A}(\Lambda)}$$
hence
 $$\frac{\vol_{A}(\Lambda)}{ \vol_{A}(\order_A)} \asymp 
\left(\frac{\disc(\data)}{\disc(A)} \right)^{-1/2},$$ where $\asymp$ may be replaced by equality when the residue characteristic
is greater than $n$.  
On the other hand, our normalizations of measure are so that $\vol_{A}(\order_A^{\times}) = \vol_{A}(\order_A)$.

It remains (cf. \eqref{volnorm}) to show $$\vol_{A}(\Lambda^{\times}) \geq \max(q^{-n},1-n/q) \vol_{A}(\Lambda).$$ 
That is effected by the following comments:

\begin{enumerate}
\item \label{one}  Let $\pi$ be any uniformizer of $k$.
 Then $1+ \pi \Lambda \subset \Lambda^{\times}$;
\item \label{two} The fraction of elements in $\Lambda$
which are invertible is $\geq 1-n q^{-1}$.
\end{enumerate}
Both of these statements would remain valid when $\Lambda$ is replaced by any
sub-order of $\order_{A}$, containing $\order_k$. 

\eqref{one} follows, because one may use the Taylor series expansion
of $(1 + \pi \lambda)^{-1}$, for $\lambda \in \Lambda$, to invert it. 

Note that $\order_{A}$ has at most $n$ maximal ideals. 
If $\lambda \in \Lambda$ does not belong to any of these ideal,
it is annihilated by a monic polynomial of degree $\leq n$ with  coefficients in $\order_{k}$ and constant term in $\order_{k}^{\times}$. 
Therefore $\lambda^{-1} \in \Lambda$ also. 
The volume of each maximal ideal is, as a fraction of the volume of $\order_{A}$, at most $q^{-1}$. This shows \eqref{two}. 
\qed

\subsection{Local bounds.}
For any Schwartz function $\Psi$ on $k^n$, we shall write $\Psi_{\algebra}$ for the function $\iota \Psi$ on $A$,
 that is to say,$$\Psi_{\algebra}\ :\ {x} \mapsto \Psi({x}^{\iota}).$$ 

Our final goal is to discuss
bounds for integrals of the form $\int_{x\in A^\times} \Psi_A(x) \psi(x) d^{\times}_A x$, where $\psi$ is a character of $\algebra^{\times}$
(we will eventually normalize it to make it independent of $\iota$). 
We wish to find useful bounds for this, when $\Psi$ is fixed and the data $\data$ is allowed to vary. We will do so in terms of
the following norms:
let $\check{\Psi}$ to be the Fourier transform of $\Psi$, w.r.t. the character $e$,
i.e.\ $\check{\Psi}(y) = \int_{y \in k^n} \Psi(y) e(xy) dy$. 
Put
$$\|\Psi\|=
 \max ( \int_{k^n} |\Psi(x)|d\vol_k(x),   \int_{k^n}  |\check{\Psi}(x)|
d\widehat{\vol}_k(x)  ).$$

The following proposition presents bounds. The reader should ignore
the many constants and focus on the
fact that these bounds decay as $\disc(A)$ or $\disc(\data)$ 
increase. In our present language, this is the analogue of the discussion of \eqref{fdef}.

\begin{Proposition} \label{localbounds} (Local bounds)
Let $\psi$ be a unitary character of $A^{\times}$, $\Psi$
a Schwartz function on $k^n$. Set:
$$I(\Psi) =  |\det h|^{-1/2} \int_{A^\times} \Psi(x^{\iota}) |{x}|_A^{1/2} \psi(x) d_A^{\times}x,$$ 
(note that $|I(\Psi)|$ is independent of the choice of $\iota$). 
Then:
\begin{enumerate}
\item 
\begin{equation} \label{first} |I(\Psi)| \ll_{n}  C_{\Psi} \cdot C_V \cdot  
 \left( \frac{\disc(\data)}{\disc(A)}\right)^{-\frac{1}{16n^2}} \end{equation}
where $C_V = \vol(\order_A^{\times})$, 
and moreover, we may take $C_{\Psi} = 1$ when $\Psi$ is the characteristic function
of $\order_{k^n}$.

\item \begin{equation} \label{second} |I(\Psi)|  \ll_{e,n}  C_V \|\Psi\| 
(\disc(\psi)
\disc(A))^{-1/4} , 
  \end{equation}
where
$C_V =
\left( \frac{\vol(\order_A)}{\vol(\order_k)^n} \right)^{1/2} 
$ (the notion of $\disc(\psi)$ is defined in \S \ref{discdef}.) 
\end{enumerate}

In inequality \eqref{first} the implicit constant depends at most on $n$ and the
degree\footnote{We have already remarked that the implicit constants in this section may depend on this degree without explicit mention; thus this is not denoted explicitly in \eqref{first} or \eqref{second}.} of $k$ over $\R$ or $\Q_p$; in \eqref{second} it depends at most on these and on, in addition, 
 the additive character $e$ of $k$. 
\end{Proposition}

 \subsection{Local harmonic analysis and the local discriminant}
We now work on the building of $\PGL_n(k)$
and on $\PGL(A)$. In particular, to simplify notations, if $N_{1}$ and $N_{2}$
are norms either on $A$ or $k^{n}$, 
$\dist(N_{1},N_{2})$ refers to the distance between their respective homothety classes, $\dist([N_{1}],[N_{2}])$, as defined in \S \ref{local}.  

Let us recall that the norm $N_0$, defined by $N_0(x_1, \dots, x_n) = \max|x_i|$, defines a point in the building of $\PGL_n(k)$, and the norm $N_A$
defines a point in the building of $\PGL(A)$. 

\begin{Lemma} \label{extreme}
Write $\algebra = \oplus_i K_i$ as a direct sum of fields. 
 Let $\|\cdot\|_i=\|\cdot\|_{K_{i}}$ be the absolute value on $K_i$ extending
$|\cdot|$ on $k$. \footnote{Thus $\|\cdot\|_{i}=|\cdot|_{K_{i}}^{1/[K_{i}:k]}$, 
in the case when $k \neq \C$; when $k=\C$ we have simply $K_i = \C$
and the absolute value on both $k$ and $K_i$ is the usual absolute value on $\C$, according to \S \ref{NormCon}.}
For $t \in \algebra^{\times}$, set  $$\|t\| := \max_{i} \|t_{i}\|_i/\min_{i} \|t_{i}\|_i,$$
where $t=(t_{1},\dots,t_{i},\dots)$ with $t_{i}\in K^\times_{i}$.
Then, for any $t \in \algebra^{\times}$, we have 
$$\dist( t N_{A}, N_{A})
\geq \frac{1}{2} \log \|t\|.$$
 \end{Lemma}
\proof 
It is a consequence of the definition of $\dist$
and \eqref{note}.  

Let $K$ be any of the fields $K_{i}$ let $N_K$ be the canonical norm attached to $K$. 
It is easy to see that for any $t\in K^\times$,
\begin{equation}\label{inegNK}
N_K(t^{-1})^{-1} \leq \|t\|_K \leq N_K(t).
\end{equation}

Given $t=(t_{1},\dots,t_{i},\dots)\in A^\times$, $t_{i}\in K_{i}^\times$. Let $\imax$ and $\imin$ be, respectively,
 those values of $i$ for which $\|t\|_i$ is maximized and minimized. Let $x_{max}\in A$ be the element whose $\imax$-th component
 is $1$ and whose other components are $0$ and let $x_{min}$ be the element whose $i_{min}$-th component
  is $t_{\imin}^{-1}$ and whose other components are $0$. Then by \eqref{note} applied to $v_{1}=x_{max},\ v_{2}=x_{min}$ and by
  \refs{maxNA} and \refs{inegNK}, one has
\begin{align*}\dist( t N_{A}, N_{A})
&\geq \frac{1}{2} \log \bigl(\frac{N_{A}(tx_{max})}{N_{A}(x_{max})}\frac{N_{A}(x_{min})}{N_{A}(tx_{min})}\bigr)\\
&=\frac{1}{2} \log (N_{K_{\imax}}(t_{\imax})N_{K_{\imin}}(t^{-1}_{\imin}))\\
&\geq \frac{1}{2} \log(\frac{\|t_{\imax}\|_{\imax}}
{\|t_{\imin}\|_{\imin}})=\frac{1}{2} \log \|t\|
\end{align*}

 \qed 

The following Lemma shows that the discriminant of local torus data is 
related to distance-measurements on the building. 
\begin{Lemma} \label{Deltabound} One has the lower bound
$$\inf_{t \in \algebra^{\times} } \dist(N_0, t \iota^{-1} N_A) \geq \frac{1}{4n} \log(\frac{\disc(\data)}{\disc(A)})+O_n(1);$$ here one may ignore the $O_n(1)$ term when $k$ is nonarchimedean and of residue characteristic larger than $n$. 
\end{Lemma} 
\proof
The action of $\PGL_n(k)$ on the building is proper
and so the infimum is attained.  
Let $t_0$ attain the infimum and put $\Delta =  \dist(t_0 \iota^{-1} N_A,  N_0)$.  

We are going to use the characterization of $\disc(\data)$
from Lemma \ref{dvexpl}. 

Adjusting $t_0$ as necessary by an element of $k^{\times}$, we may assume:
\begin{equation} \label{cork}e^{-\Delta} N_0 \leq  t_0 \iota^{-1} N_A \leq e^{\Delta}  N_0 \end{equation}

Suppose that $x \in k$ satisfies $|x| = \exp(-2  \Delta)$. Such an $x$ exists 
(cf. \eqref{note} and the subsequent comment; $\Delta$ is a half-integral multiple
of $\log q$ in the non-archimedean case). If $y \in \algebra$
has operator norm $\leq 1$ with respect to $N_{0}$, then \eqref{cork} shows that $xy$ has operator
norm $\leq 1$ with respect to $t_0 \iota^{-1} N_A$, and vice versa. 

The set of $a \in A$ which have operator norm $\leq 1$ w.r.t. $N_0$
is exactly $\Lambda$.

The set of $a \in A$ which have operator norm $\leq 1$ w.r.t. $t_0 \iota^{-1} N_A$ is by definition the set of $a$
such that for all $x\in k^{n}$
$$t_0 \iota^{-1} N_A(xa)\leq t_0 \iota^{-1} N_A(x),$$
that is by \refs{NormCompat} 
$$N_{A}(x^{\iota^{-1}}at_{0})\leq N_{A}(x^{\iota^{-1}}t_{0}).$$
Using that $x^{\iota^{-1}}at_{0}=x^{\iota^{-1}}t_{0}a$ and changing $x^{\iota^{-1}}t_{0}$ to $t$, we see that this set
 is the set of $a\in A$ satisfying for any $t\in A$
$$N_{A}(ta)\leq N_{A}(t)$$
which is precisely $\order_{A}$.

We conclude:
$$x \order_{\algebra} \subset \Lambda\subset \order_{A}, \ \ \ \order_{A}^*\subset\Lambda^* \subset x^{-1} \order_{\algebra}^*.$$ 
The second equation is obtained by duality from the first; here $\Lambda^*$ 
and $\order_{\algebra}^*$ are dual in the sense of \eqref{dualdef}. 
Thereby,
$$\frac{\vol_{A}(\Lambda^*)}{\vol_{A}(\Lambda)} \leq \exp(4 n \Delta) \frac{\vol_{A} (\order_{\algebra}^*)}{\vol_{A}( \order_{\algebra})},$$
whence the result (cf. \eqref{volnorm} and Lemma \ref{dvexpl}).  \qed

\begin{Lemma} Let $R \geq 0$. Then:
\begin{equation}\label{volbound}\frac{\vol\{ t \in \algebra^{\times}/k^{\times}:
 \log \|t\| \in [R, R+1] \}}{\vol_{A}(\order_{\algebra}^{\times})/\vol_{k}(\order_k^{\times})}\ll_n (1+R)^{n-1} \end{equation}
\end{Lemma}
Here, $\|t\|$ is as in the statement of Lemma \ref{extreme}. 

\proof
The archimedean case may be verified by direct computation. 

Consider $k$ nonarchimedean. We may write $\algebra = \oplus_{i=1}^{r} K_i$ as a sum of $r$ fields. The map
$$a = \oplus a_i \mapsto \frac{\log \|a_i\|_i}{\log(q)}$$
gives an isomorphism of $\algebra^{\times}/\order_{A}^{\times}$
with a finite index sublattice $Q$ of $\frac{1}{n!} \Z^r$. 
Moreover, $k^{\times}/\order_k^{\times}$ is identified with 
the sublattice $Q' \subset Q$ generated by $(1,1,1, \dots, 1) \in \Z^r$. 
Finally, the norm $\log \|t\|$ on $\algebra^{\times}/k^{\times}$
descends to a function on $Q/Q'$ and is described explicitly as:
$$(\mu_1, \dots, \mu_r) \in Q \mapsto \log(q) \left(\max_{i} \mu_i - \min_{i} \mu_i\right)
.$$
The LHS of \eqref{volbound} is thereby bounded by:
$$\# \{\mu \in \frac{1}{n!} \Z^r/\Z: \left(\max_{i} \mu_i - \min_{i} \mu_j \right) \leq
\frac{R+1}{ \log 2}\},$$
which is bounded as indicated, since $r \leq n$. 
\qed

\begin{Lemma} \label{mgrnds}
For any $\alpha \in (0,1)$, 
\begin{equation} \label{basicest}
 \frac{ \int_{t \in A^{\times}/k^{\times}}\exp(-\alpha \dist( N_0 , t \iota^{-1} N_{A} ))}{\vol_{A}(\order_{\algebra}^{\times})/\vol_{k}(\order_k^{\times})}
   \ll_{\alpha,n} \left( \frac{\disc(\data)}{\disc(\algebra)}\right)^{-\frac{\alpha}{8n}}. \end{equation}
\end{Lemma}

\proof 
Let $\Delta, t_0$ be as in Lemma \ref{Deltabound}.
Using it,  \eqref{NormCompat}, Lemma \ref{extreme}, and  the triangle inequality:
\begin{multline} \dist( N_0,  t \iota^{-1} N_{\algebra} ) \geq 
\dist(t \iota^{-1} N_{\algebra} , t_0 \iota^{-1} N_{\algebra})
- \dist(t_0 \iota^{-1} N_{\algebra} ,  N_0) \\
\geq \frac{1}{2} \log(\|t/t_0\|) - \Delta. \end{multline}

To estimate \eqref{basicest}, we split the $t$-integral
into regions when $\log \|t/t_0\| \leq [4 \Delta]$ and $\log \|t/t_0\| \in
 [R, R+1]$, where $R$ ranges through integers $\geq [4 \Delta]$. 
Here $[4 \Delta]$ is the greatest integer $\leq 4 \Delta$. 
Thereby, the left-hand side of \eqref{basicest} is bounded by
$$ C(n)
\left( (1+\Delta)^n \exp(- \alpha \Delta) + \sum_{R \geq [4 \Delta]} \exp(
\alpha(\Delta- R/2)) (1+R)^{n-1} \right) $$
To conclude, we bound $\Delta$ from below using Lemma \ref{Deltabound}. 
\qed

\subsection{The action of $\GL(A)$ on $L^2(A)$.}Let $V$ comprise $-n/2$-homogeneous functions on $\algebra$, i.e.
$$V = \{f : A \rightarrow \C: f(\lambda x) = |\lambda|_k^{-n/2} f(\mbfx), \lambda \in k,  \mbfx\in \algebra.\}$$

The group $\PGL(\algebra)$ acts on $V$, via
$$gf(\mbfx) = f(\mbfx g)  |\det g|_k^{1/2}$$ 
The space $V$ possesses a (unique up to scaling) natural $\GL(\algebra)$-invariant inner product.   We shall normalize it as follows:
for any Schwartz function $\Phi$ on $\algebra$, let $\widetilde{\Phi}$ be
its projection to $V$, defined as $$\widetilde{\Phi}(\mbfx) =
 \int_{\lambda \in k^\times} \Phi(\lambda\mbfx ) |\lambda|_k^{n/2} d^{\times}\lambda.$$
   We normalize the inner product
 $\langle \cdot, \cdot \rangle$ on $V$ so that for any $v \in V$: 
$$\langle v, \widetilde{\Phi} \rangle = \int_{x \in \algebra} v(\mbfx) 
\overline{\Phi(\mbfx)} d_{A}\mbfx.$$
In particular, let $\Phi_1, \Phi_2$ be Schwartz functions on $\algebra$. We have:
\begin{align} \label{Hintinpractice} \int_{t \in \algebra^{\times}/k^{\times}} \langle t \cdot \widetilde{\Phi_{1}}, \widetilde{\Phi_{2}} \rangle &=
\int_{t \in \algebra^{\times}} \int_{ \mbfx\in \algebra} \Phi_{1}(\mbfx t)|t|_{A}^{1/2} \overline{ \Phi_{2}(\mbfx)} d_{A}\mbfx d_{A}^{\times}t  \\& = 
(\int_{y \in A^\times} \Phi_{1}(\mbfy) |\mbfy|_A^{1/2} d^\times_A\mbfy) \overline{ (\int_{\mbfx\in A^\times} \Phi_{2}(\mbfx) |\mbfx|_A^{-1/2} d_A \mbfx) }\nonumber
 \\& =   \zeta_A(1)
(\int_{y \in A} \Phi_{1}(\mbfy) |\mbfy|_A^{-1/2} d_A\mbfy) \overline{ (\int_{\mbfx\in A} \Phi_{2}(\mbfx) |\mbfx|_A^{-1/2} d_A \mbfx) }\nonumber \end{align}
Let us note that we use, in the above reasoning and at various other points in the text, the evident fact that the measure of $A - A^{\times}$ is zero.

Let $K \subset \GL(\algebra)$ be a maximal compact subgroup, corresponding
to the stabilizer of $\iota N_0$, i.e.\ 
in the nonarchimedean case, the stabilizer of the lattice
$(\order_{k^n})^{\iota^{-1}}$. 
Let $\Xi_0: \GL(\algebra) \rightarrow \C$ be the Harish-Chandra spherical function with respect to $K$. 
For two vectors $v_1, v_2 \in V$, and $\sigma \in \GL(\algebra)$, we have the bound (see \cite{CHH}
and also equation \eqref{HCS}):
\begin{multline}\label{MC}\langle \sigma v_1, v_2 \rangle \leq (\dim K v_1)^{1/2} (\dim K v_2)^{1/2} \Xi_0(\sigma)^{1/n}  \|v_1\|_{2} \|v_2\|_{2}
\\ \ll_{\alpha} (\dim K v_1)^{1/2} (\dim K v_2)^{1/2} 
\exp(-\alpha\dist(\sigma\iota N_0, \iota N_0)) \|v_1\|_{2} \|v_2\|_{2}.
\end{multline}
for any $\alpha<1/n$. Here distances are measured between homothety classes of norms.

\subsection{Proof of the first estimate in Proposition \ref{localbounds}.}
Since $\psi$ is a unitary character, it is sufficient (by taking absolute values) to assume that $\psi$ is trivial and $\Psi$, non-negative.

Let $\Psi$ be a non-negative Schwartz function on $k^n$;
put $\Psi_A := \Psi \circ \iota$, 
$\Phi$  the characteristic function of the unit ball of $\iota N_0$,
and $\Phi_2 = 1_{\order_{A}}$. 

These are all Schwartz functions on $A$. In the nonarchimedean case
 $\Phi$ is the characteristic function of 
$(\order_{k^n})^{\iota^{-1}}$.
Also, the definition of $h_A$ (see \eqref{hdef2}) shows that $\Phi_2 = h_A \Phi$. Indeed, $h_A \Phi$ is the characteristic function of $$\{x \in A: N_0( (x h_A)^{\iota}) \leq 1\}=\{x \in A: N_A(x) \leq 1\}=\order_{A}.$$
 Consequently,
$h_A \widetilde{\Phi} = |\det h|_k^{1/2} \widetilde{h_A \Phi} = |\det h|_k^{1/2} \widetilde{\Phi}_2$.

We shall proceed in the case when $k$ nonarchimedean, the archimedean case being similar (the only difference: one needs to decompose $\Psi$ as a sum of $K$-finite functions in the archimedean case, and the implicit constant will be bounded by a Sobolev norm of $\Psi$).

Let us observe that there is a constant $C_{\Psi}\geq 0$,
equal to $1$ when $\Psi = \Phi$, so that 
$$\langle \widetilde{\Psi}_A, \widetilde{\Psi}_A \rangle \leq C^2_{\Psi}
\langle \widetilde{\Phi}, \widetilde{\Phi}\rangle.$$

Indeed for some $\lambda_{\Psi}\in k$, one has
$\Psi(x)\leq \|\Psi\|_{\infty}1_{\order_{k^n}}(\lambda_{\Psi}x),\ x\in k^n$;
this bounds $\Psi_A$ in terms of $\Phi$ and leads
to the above-claimed bound.

For $t \in A^{\times}/k^{\times}$, 
\begin{multline} \label{est1}|\det h|_k^{1/2}
\langle t \widetilde{\Psi_A}, \widetilde{\Phi_2} \rangle = \langle  t \widetilde{\Psi_A}, h_A \widetilde{\Phi} \rangle= 
 \langle h_A^{-1} t \widetilde{\Psi_A}, \widetilde{\Phi} \rangle \\  \ll_{\alpha} \|\widetilde{\Psi_A} \|_{2} \|\widetilde{\Phi}\|_{2}
\dim(K \widetilde{\Psi_A}) \exp(-\alpha \dist( h_A^{-1} t \iota N_0,  \iota N_0)), \ \ \alpha < \frac{1}{n}.\end{multline}
We have applied  \eqref{MC} with $v_1 =\widetilde{\Psi_A},  v_2 = \widetilde{\Phi}, 
\sigma = h_A^{-1} t$; observe that our choice of $\Phi$ means $\dim(K \widetilde{\Phi})=1$. 

We observe that, using the definitions \eqref{hdef}, \eqref{hdef2} and the compatibility
\eqref{NormCompat}, we have $\dist(h_A^{-1} t \iota N_0, \iota N_0) = \dist(N_0, t^{-1} \iota^{-1}  N_A)$.
To write out every step, this follows from the chain of equalities 
\begin{gather*}\dist(h_A^{-1} t \iota N_0,  \iota N_0) = \dist(t \iota N_0, h_A  \iota N_0) 
= \dist( t \iota N_0, N_A)   = \dist(N_0, \iota^{-1} t^{-1} N_A)\\
 = \dist(N_0, t^{-1} \iota^{-1} N_A)
= \dist(t N_0, \iota^{-1} N_A)  = \dist(N_0, t^{-1} \iota^{-1} N_A).
\end{gather*}

 We now integrate over $t \in A^{\times}/k^{\times}$.
By \eqref{Hintinpractice} and  \eqref{est1}  we have, for any $\alpha < \frac{1}{n}$, 
\begin{multline}\label{claim133} \left| \int \Psi(\mbfx^{\iota})  |\mbfx|_A^{-1/2} d_A\mbfx  \right|
\left| \int \Phi_2(\mbfx)  |\mbfx|_A^{-1/2} d_A\mbfx  \right|   = 
\zeta_A(1)^{-1} \int_{A^{\times}/k^{\times}} \langle t \widetilde{\Psi}, \widetilde{\Phi} \rangle  \\
 \ll_{n, \alpha}C_{\Psi}
 |\det h|_k^{-1/2} (\dim \GL_n(\order_k) . \Psi)  \|\widetilde{\Phi}\|_2^2
  \int_{t \in \algebra^{\times}/k^{\times}} e^{-\alpha \dist(N_0, t^{-1} \iota^{-1} N_A )}  \end{multline}
We show below that
\begin{equation} \label{claim134} \frac{\|\widetilde{\Phi}\|_2^2}{\int_{\mbfx\in A} \Phi_2(\mbfx) |\mbfx|_A^{-1/2} d_A\mbfx  } 
\ll_n 
|\det h|_k
 \vol_{k}(\order_k^{\times}).\end{equation}
Combining \eqref{claim133} and \eqref{claim134} 
with Lemma \ref{mgrnds}, establish the first claim of Proposition
\ref{localbounds}.  

To prove \eqref{claim134},  proceed as follows:
First of all, 
$$\int_{\mbfy} \Phi_2(\mbfy) |\mbfy|_A^{-1/2} d_A \mbfy=\int_{\order_{A}}|\mbfy|_A^{-1/2} d_A \mbfy=\vol_{A}(\order_{A}^\times)\zeta_{A}(1/2).$$

 Noting that $\widetilde{\Phi}(\mbfx) 
= \iota N_0(\mbfx)^{-n/2}
\int_{|\lambda|_k \leq 1} |\lambda|_k^{n/2} d^{\times} \lambda$,
we see: 
\begin{align} \nonumber \label{rep1} \| \widetilde{\Phi}\|_2^2 &= 
\int_{|\lambda|_k \leq 1} |\lambda|_k^{n/2} d^{\times} \lambda \cdot \int_{\iota N_0(\mbfx) \leq 1}
 \iota N_0(\mbfx)^{-n/2} d_A \mbfx\\ \nonumber
 &
 =\vol_{k}(\order_{k}^\times)\zeta_{k}(n/2)\int_{\iota N_0(\mbfx) \leq 1}
 \iota N_0(\mbfx)^{-n/2} d_A \mbfx
\end{align}
 Noting that $\iota N_0=h_A^{-1}N_A$, we have
 $$ \int_{\iota N_0(\mbfx) \leq 1}
 \iota N_0(\mbfx)^{-n/2} d_A \mbfx= \int_{N_{A}(\mbfx h_{A}^{-1}) \leq 1}
 N_{A}(\mbfx h_{A}^{-1})^{-n/2} d_A \mbfx,$$ 
 so that, making the change of variable $x'=xh_{A}^{-1}$, the previous integral equals
 $$|\det h_{A}|_{k}\int_{N_{A}(x)\leq 1}N_{A}(x)^{-n/2}d_{A}x=|\det h|_{k}\int_{\order_{A}}N_{A}(x)^{-n/2}d_{A}x.$$
 For $k$ nonarchimedean, the last integral equals ($\pi$ denote an uniformizer of $k$)
\begin{align*}
\sum_{j\geq 0}|\pi^j|_{k}^{-n/2}\int_{N_{A}(x)=|\pi^j|_k}d_{A}x&=\sum_{j\geq 0}|\pi^j|_{k}^{-n/2}\int_{N_{A}(x\pi^{-j})=1}d_{A}x\\
&=\vol(\{x,\ N_{A}(x)=1\})\sum_{j\geq 0}|\pi_{k}^j|^{n/2}\\
&=\vol(\{x,\ N_{A}(x)=1\})\zeta_{k}(n/2).
\end{align*}
 Combining these, the left-hand side of \eqref{claim134}
is bounded by:
$$\ll_n
  |\det h|_k \vol_{k}(\order_{k}^\times)\frac{\vol_A(\{x \in A: N_A(x) = 1\})}{\vol_A(\order_A^{\times})}
$$
The last ratio is easily seen to be bounded above by $(1 + n/(q-1))$; in particular,
it is bounded above in terms of $n$ and the claim \eqref{claim134} follows. 
\qed

\subsection{Proof of  the second estimate in Proposition \ref{localbounds}.}

For any character $\psi$ of $A^\times$
we set $\psi_s(x) = \psi(x) |x|_A^s$. Let us comment
that the notation $\psi^{-1}_s$ always denotes the character
$x \mapsto \psi^{-1}(x) |x|_A^s$. In other words, we apply the operation
of twisting by $|x|_A^{s}$ {\em after} the operation of inverting $\psi$.

The following result is proved in Tate's thesis.  See \cite[(3.2.1), (3.2.6.3), (3.4.7)]{Tate-Corvallis}.
\begin{Lemma}\label{LFElemma} (Local functional equation)
 Let $\Phi$ be a Schwartz function on $A$, 
and set $$\widehat{\Phi}(x) = \int_{y \in A} e_{A}(xy) \Phi(y) d_A y.$$
Then, for a unitary character $\psi$ of $A^{\times}$, 
 \begin{equation} \label{LFE}\epsilon(A, \psi, s, e_{A}) \frac{\int_{A^\times} \Phi(x) \psi_s(x)d_A^{\times}x}{L(A, \psi,s)}
  = \frac{\int_{A^\times} \widehat{\Phi}(x) \psi^{-1}_{1-s}(x) d_A^{\times}x}{L(A, {\psi^{-1}},{1-s})} \end{equation} 
where $s \mapsto \epsilon(A, \psi,s, e_{A})$ is a holomorphic function of exponential type. 
More precisely, both sides of \eqref{LFE} are holomorphic, and:
\begin{enumerate}
\item If $\widehat{\vol}_{A}$ denotes the Haar measure dual to $\vol_A$ under the Fourier transform\footnote{With our choice of normalizations, $\widehat{\vol}_A = \vol_A$, but we prefer to phrase the Lemma in a fashion that is independent of choice of measures.},  $|\epsilon(\algebra, \psi,s, e_{A})|^2  = 
\frac{\vol_A}{\widehat{\vol}_A}$ when $\Re(s) = 1/2$. 
\item $|\epsilon(\algebra, \psi,s, e_{A})| = 
|\epsilon(\algebra, \psi, 0,e_{A})| (\delta(e)^n\disc(\psi)
\disc(\algebra))^{-\Re(s)}$, where $\delta(e)$ is a positive constant depending on the additive character $e$ of $k$.
\end{enumerate}
\end{Lemma}

We remark that one may take $\delta(e)=1$ in the unramified case, i.e. when $k$ is nonarchimedean 
and $e$ is an unramified character of $k$. 

We now proceed to the proof of the second estimate in Proposition \ref{localbounds}. Recall that $\psi$ is unitary.  To ease our notation,
we suppose (as we may do, without loss of generality) that $\|\Psi\|=1$. 

For $\Re(s) = 1$, we have:
\begin{equation} \label{leftbound} \left| \int_{\algebra^\times} \Psi_A(x) \psi_s(x) d^{\times}_A x \right| \leq   \zeta_{\algebra}(1) 
\frac{\iota_* \vol_A}{\vol_{k^n}}\ll_{n}\frac{\iota_* \vol_A}{\vol_{k^n}}
 \end{equation}

For $\Re(s) = 0$ we apply \eqref{LFE} to reduce to \eqref{leftbound}, obtaining:
\begin{equation} \label{rightbound2}\left|\frac{ 
 \int_{\algebra^\times} \Psi_{A}(x) \psi_s(x) d^{\times}_A x } {L(\algebra, \psi,s) / L(\algebra,{\psi^{-1}},{1-s})} \right|  \ll_{e,n}
(\disc(\psi) \disc(\algebra))^{-1/2},  
\end{equation}where the constant implied depend at most on the additive character $e$ of $k$ and on $n$

We may now apply the maximum modulus principle to interpolate between 
\eqref{leftbound} and \eqref{rightbound2}. The simplest
thing to do would be to apply the maximum modulus principle
to the holomorphic quotients that occur in \eqref{LFE}. 
This would be fine for nonarchimedean places; however, for archimedean places,
this would run into some annoyances owing to the decay of $\Gamma$-factors.
We proceed in a slightly different way. 

Let $F(s)$ be an analytic function in a neighbourhood of the strip $0 \leq \Re(s) \leq 1$ so that $F(s) L(\algebra, \psi,s)$ is holomorphic.  We shall choose
$F(s)$ momentarily. 
Let $U$ be the right-hand side of \eqref{leftbound}, and $V$
the right-hand side of \eqref{rightbound2}. 
 
Take $$h(s) := U^{-s}V^{-(1-s)}F(s) \int_{\algebra^\times}
 \Psi_{A}(x) \psi_{s}(x) d_{A}^{\times}x,$$ 
Note that $h(s)$ is holomorphic in $0 \leq \Re(s) \leq 1$. We have
$$|h(s)| \ll_{e,n} \begin{cases} |F(s)|, \Re(s) = 1. \\
\frac{|F(s)  L(A,\psi,s)|}{|L(A,{\psi}^{-1},{1-s})|}, \Re(s) = 0. \end{cases}
$$
\begin{enumerate}
\item $k$ nonarchimedean.
We choose $F(s) = L(\algebra, \psi,s)^{-1}$. 
Then $$\sup_{\Re(s) = 1} |F(s)|  \ll_{n}1$$
whereas $\sup_{\Re(s) = 0}  \frac{|F(s)  L(A,\psi,s)|}{|L(A,{\psi}^{-1},{1-s})|}$ is also bounded by $\zeta_A(1)$. 

Therefore, by the maximum modulus principle, one has for $\Re(s) = 1/2,$ $|h(s)| \ll_{e,n} 1$. 
On the other hand, for $\Re(s) = 1/2$, $|F(s)| \geq \zeta_A(1/2)^{-1}$. 
Therefore, for $\Re(s) = 1/2$,
$$\left| \int_{A^\times} \Psi_{A}(x) \psi_{s}(x) d_{A}^{\times}x \right| 
\ll_{n,e} \sqrt{UV}. $$
\item $k$ archimedean. In explicit terms,
$L(\algebra, \psi_s)$ is a product of a finite number of $\Gamma$-factors $$\prod_{i}\Gamma_{K_{i}}({s+\nu_{i}}),$$ where $\Re(\nu_{i})$ are non-negative integers,
and 
\begin{equation} \label{gammaKdef} \Gamma_{K}(s) =\begin{cases} \pi^{-s/2} \Gamma(s/2) = \Gamma_{\R}(s) &\hbox{ if $K=\R$}\\
2(2\pi)^{-s} \Gamma(s)=\Gamma_{\R}(s)\Gamma_{\R}(s+1)&\hbox{ if $K=\C$}.
\end{cases}\end{equation}

We take $F(s) = \prod_{\Re(\nu_{i}) = 0} (s+\nu_{i}) (s+\nu_{i} - 100)^{-1}.$
Then $\sup_{\Re(s) = 1} |F(s)|$ and $\sup_{\Re(s) = 0}  \frac{|F(s)  L(A,\psi,s)|}{|L(A,\check{\psi},{1-s})|}$
 are both bounded above by functions of $[A:\R]\leq 2n$.  The first is clear;
for the second:
$$ \frac{|F(s)  L(A,\psi,s)|}{|L(A,\check{\psi},{1-s})|}
= \prod_{\Re(\nu_i) = 0} \frac{s + \nu_{i}}{(s+\nu_{i}-100)} \frac{\Gamma_{K_{i}}({s+\nu_i})} {\Gamma_{K_{i}}({1-s + \ov\nu_i})}.
\prod_{\Re(\nu_i) \neq 0}
 \frac{\Gamma_{K_{i}}({s+\nu_i})} {\Gamma_{K_{i}}({1-s + \ov\nu_i})}.
$$
It is not hard to see that the right-hand side is, indeed,
bounded above when $\Re(s) = 0$, by a function of $[A:\R]$. 

We conclude that, for $\Re(s)= 1/2$, $|h(s)| \ll_{e,n} 1$.
But, for $\Re(s) = 1/2$, we also see $|F(s)| \gg_{n} 1$. 
We conclude that for $\Re(s)=1/2$:
$$\left| \int \Psi_{A}(x) \psi_{s}(x) d_{A}^{\times}x \right| 
\ll_{n,e} \sqrt{UV}.$$
\end{enumerate}

We have therefore shown that, for $\|\Psi\| = 1$, 
\begin{multline*}|\det h|^{-1/2} \left| \int \Psi_A(x) \psi(x) |x|_{A}^{1/2}d_A^{\times}(x) \right|
\\
\ll_{e,n}  |\det h|^{-1/2}
\left(\frac{\iota_* \vol_A}{\vol_{k^n}}\right)^{1/2}
(\disc(\psi) \disc(\algebra))^{-1/4} 
\end{multline*}
Taking into account \eqref{hdef2}, we see that
the proof of the second assertion of Proposition \ref{localbounds} is complete.

\section{Eisenstein series: definitions and torus integrals} \label{EisSec}
In this section,
we define the {\em Eisenstein series} on $\GL_n$ and
give a formula  (Lemma \ref{lem:int}) for their
integrals over tori. This formula will later be used to derive \eqref{HeckeFormula}. This section is included merely to make the paper self-contained.
 Indeed these computations go back to Hecke (see also \cite{Wielonsky}).

\subsection{Eisenstein series -- definition and meromorphic continuation.}

Throughout this section, we follow the notation of \S \ref{Notation}.

For each place $v$ of $F$, let $\Psi_v$ be a Schwartz function\footnote{Recall that this has the usual meaning if $v$ is archimedean, and means: locally constant of compact support, otherwise.}
   on $F_v^n$. We suppose that, for almost all $v$, the function $\Psi_v$ coincides with the characteristic function of $\order_{F,v}^n$.  Let  $\Psi :=
\prod_{v} \Psi_v$ be the corresponding Schwartz function on $\adele^n$. 

Put, for $g \in \GL_n(\adele)$ and $\chi \in \Omega(C_F)$, 
\begin{equation} E_{\Psi}(\chi,g) = \int_{t \in \adele^{\times}/F^{\times}}
\sum_{\mathbf{v} \in F ^ n - \{ 0 \}} \Psi(\mathbf{v} t  g) \chi(t) d ^{\times}t. \end{equation} 
The integral is convergent when $\Re \chi$ is sufficiently large. Note that $E_{\Psi}(\chi,g)$ has central character $\chi^{-1}$.  

To avoid conflicting notations, we shall set occasionally, for $s\in \Cc$
$$E_{\Psi}(\chi,s,g):=E_{\Psi}(\chi_{s},g)=E_{\Psi}(\chi|\cdot|_{\adele}^s,g).$$

Let $[\cdot, \cdot]$ be a nondegenerate bilinear pairing on $F^n$.
The pairing $[\cdot, \cdot]$  gives also a nondegenerate bilinear pairing $\adele^n \times \adele^n \rightarrow \adele$. 

Let $g^*$ be defined so that
$[\v_1 g, \v_2g^*] = [\v_1, \v_2]$. Therefore $|\det g|_{\adele} |\det g^*|_{\adele} = 1$ for $g \in \GL_n(\adele)$. 
The pairing defines a Fourier transform:
$$\widehat{\Psi}(\v^*) = \int_{\adele^n} \Psi(\mathbf{v}) e[\mathbf{v}, \mathbf{v}^*]\ d\mathbf{v};$$
in particular the Fourier transform of $\v \mapsto \Psi(\v t g)$ is
$$\v \mapsto |t|_{\adele}^{-n} |\det g|_{\adele}^{-1}  \widehat\Psi(\v t^{-1} g^*).$$
Recall that Poisson summation formula  shows that
$$\sum_{\mathbf{v} \in F^{n}} \Psi(\v) = \sum_{\mathbf{v} \in F^{n}} \widehat{\Psi}(\v).$$

\begin{Proposition} \label{meromorphic}
$E_{\Psi}(\chi,g)$ continues to a meromorphic function in the variable $\chi$, with simple poles when $\chi = 1$ and $\chi = |\cdot|_{\adele}^n$. 
One has
\begin{align*}\res_{\chi=1}E_{\Psi}(\chi,g)&=-\vol(\adele^{(1)}/F^{\times}) \Psi(0)\\
\res_{\chi=|\cdot|_{\adele}^n}E_{\Psi}(\chi,g)&=|\det g|_{\adele}^{-1}
  \vol(\adele^{(1)}/F^{\times}) \int_{\adele^n} \Psi(x)dx.
  \end{align*}  
Moreover,
$$|\det g|_{\adele} E_{\Psi}(\chi, g) = E_{\widehat{\Psi}}(\chi^{-1},n, g^{*}).$$ 
\end{Proposition}

\proof
Split the defining integral $E_{\Psi}(\chi,g)$ into $|t|_{\adele} \leq 1$ and $|t|_{\adele} \geq 1$. Apply Poisson summation to the former, and then substitute $t \leftarrow t^{-1}$. The result, valid for $\Re \chi \gg 1$, is:
\begin{multline}  \int_{|t|_{\adele} \geq 1} d^{\times}t \left( \sum_{v \in F^n} \Psi(\v t g) - \Psi(0)\right) \chi(t) \\  + 
|\det g^*|_{\adele} \int_{|t|_{\adele} \geq 1} d^{\times}t  \left(\sum_{\v \in F^n}   \widehat{\Psi}(\v t g^*) - \widehat{\Psi}(0)\right)  \chi_n^{-1} (t)
\\ - \Psi(0) \int_{|t|_{\adele} \geq 1} \chi(t) d^{\times}t + \widehat{\Psi}(0)
\int_{|t|_{\adele} \leq 1} \chi_{-n}(t) d^{\times}t 
\end{multline}
The last two terms can be explicitly evaluated. If $\chi$ is of the form $|x|_{\adele}^s$, they are equal to $-\frac{\Psi(0)}{s}\vol(\adele^{(1)}/F^{\times})$ and $|\det g^{*}|_{\adele}\frac{\int \Psi}{s-n} \vol(\adele^{(1)}/F^{\times})$, respectively; otherwise, they are identically zero.   The former two terms
define holomorphic functions of $\chi$. 
\qed
\subsection{Torus integrals of Eisenstein series}
Put $\Psi_K = \Psi(x^{\iota} g_{\data})$, a function on $\adele_K$. 
Let us recall we have fixed global torus data $\data = (K \subset M_n(F), g_{\data} \in \GL_n(\adele))$.

\begin{Lemma} (Integration of Eisenstein series over a torus) 
  \label{lem:int}. 

Let $\psi \in \CharK$
be so that $\chi  =  \psi|_F$.  Then:

\begin{equation} \label{Eis} \int_{\T_K(F) \backslash \T_K(\adele)} 
E_{\Psi}(\chi,t g_{\data}) \psi(t) dt =  \int_{\adele_K^{\times}} \Psi_K(y) \psi(y) d ^{\times} y \end{equation}

\end{Lemma}
Let us observe that, owing to the restriction $\chi = \psi|_F$, the map $t \mapsto E_{\Psi}(\chi,tg) \psi(t)$ indeed
defines a function on $\T_K(F) \backslash \T(\adele)$. 

The integral on the right-hand side can be expressed as a product over places of $K$; for almost all places, the resulting (local) integral equals an $L$-function. This is explicitly carried out in \eqref{localL}. 
Thus, the Lemma indeed gives the reduction of torus integrals of Eisenstein
series to $L$-functions, as discussed in \S \ref{outline}. 

\proof 
We unfold. 
\begin{multline}
\int_{\T_K(F) \backslash \T(\adele)} \psi(t) E_{\Psi}(\chi,t g_{\data}) dt  \\  = 
\int_{u \in \adele^{\times}/F^{\times}} du \int_{t \in \T_K(F) \backslash \T_K(\adele)}  \, dt \psi(t) 
\sum_{\x \in F^n - \{0\}} \Psi^{g_{\data}}( \x. u . t ) \chi(u)
\\ =  \int_{t \in (\adele_K^{\times}/K^{\times})/( \adele^{\times}/F^{\times})}  \int_{u \in \adele^{\times}/F^{\times}}
 \sum_{\mbfx \in K^{\times}} \Psi((\mbfx . u . t)^{\iota} g_{\data}) \chi(u) \psi(t) du dt
\\ = \int_{t \in \adele_K^{\times}/K^{\times}} \sum_{\mbfx \in K^{\times}} \Psi_K(\mbfx. t) 
\psi(t) d^{\times}t   = \int_{t \in \adele_K^{\times}} \Psi_K(t) \psi(t) d^{\times}t 
\end{multline}  \qed

\begin{Lemma} \label{lem:CNF} (Class number formula)
The measure of $\T_K(F) \backslash \T_K(\adele)$
equals $$n \frac{\mathrm{Res}_{s=1}\zeta_K(s)}{\vol(\adele^{(1)}/F^{\times})} .$$  
\end{Lemma}
\proof
We set $g=1$, $\psi = |\cdot|_{\adele_K}^s, \chi= |\cdot|_{\adele}^{ns}$ and take residues in \eqref{Eis}
as $s \rightarrow 1$. 

We first remark that, for almost all $v$, the local integral $\int_{K_v^{\times}} \Psi_{K,v}(t) |t|_v^{s} d^{\times}t $ equals the local zeta function $\zeta_{K,v}(s) := \prod_{w|v} \zeta_{K_w}(s)$. Taking residues yields:
$$\frac{\vol (\T_K(F) \backslash \T_K(\adele)) \cdot \vol(\adele^{(1)}/F^{\times})}{n} \cdot  \int_{\adele^n} \Psi =  \mathrm{Res}_{s=1} \zeta_K(s) . \prod_{v} \frac{\int_{K_v^{\times}} \Psi_{K,v} (t) |t|_v d^{\times}t}{\zeta_{K,v}(1)}$$
where almost all factors in the infinite product are identically $1$. 
The result follows from the choice of the measure. 
 \qed

\section{Eisenstein series: estimates} \label{core}

Let us explain by reference to \S \ref{ConnectionToLFunctions} the contents of this section:

In \S \ref{EisSec} we have established the general form of \eqref{HeckeFormula}.
We are going to assume known a subconvexity bound \eqref{Hyp}, which
one can see as a generalization of the bound \eqref{globalsubconvex}) from the introduction.
The results of \S \ref{toruslocal} in effect
establish the analog of \eqref{localsubconvex}.

Putting these together, we shall obtain in the present section -- Proposition \ref{centralbound} -- a slightly disguised form of \eqref{siegel}.  This disguised form is translated to a more familiar $S$-arithmetic context in the next section.

\subsection{Assumed subconvexity}

Our result makes the assumption of a certain {\em sub-convexity estimate}. 

In order to state what this means, we need to recall the notion
of archimedean conductor. 
For a character $\omega$ of a archimedean local field $k$, we define the archimedean conductor
\begin{equation}\label{Cinftydef}
C(\omega)=\prod_{i}(1+|\nu_{i}|)
\end{equation}
where the $\nu_{i} \in \C$ are so that
$L(\omega, s)=\prod_{i}\Gamma_{\Rr}(s+\nu_{i}), \hbox{say}$. 
(See  \eqref{gammaKdef} for definitions of $\Gamma_{\Rr}$.)
If $\omega$ is unitary, then $\Re(\nu_{i})\in \frac{1}{2}\N.$
For $\chi \in \Char$,  and $v$ archimedean, we put $C_v(\chi)$ to be the  archimedean conductor
of $\chi|_{F_v}$, and let $ C_{\infty}(\chi) =
\prod_{v |\infty} C_v(\chi)$. 
Similarly, one defines $C_{\infty}(\psi)$ for $\psi \in \CharK$.

Given a unitary character $\psi \in \Omega(C_K)$, we shall assume known the following bound
\begin{equation} \label{Hyp} |L(K, s, \psi)| \ll C_{\infty}(\psi_s)^N \disc(K)^{1/4-\theta}\disc(\psi)^{1/4- \theta}, \Re(s) = 1/2\end{equation}
for some constants $N$, $\theta>0$ which depend only on $F$ and $n=[K:F]$.

 The validity of \eqref{Hyp} is a consequence of the generalized Riemann hypothesis.
The generality in which \eqref{Hyp} is known unconditionally is fairly slim, but it is enough for some applications. For a recollection of what is known unconditionally, see Appendix \ref{subconvex-a}.

\subsection{The main estimate}
In this section, we use notations as in \S \ref{Notation}.
 We regard $F,n$ as fixed throughout; thus we allow implicit constants $\ll$ to depend both on $n$ and $F$.
   In particular, any discriminants depending only on $F$, e.g.\ $\disc_v(F)$, will often be incorporated into $\ll$ notation. 

For typographical simplicity, we write $D_{\psi}, D_{\data}, D_K, D_F$
in place of $\disc(\psi)$, $\disc(\data)$, etc. in the following Proposition. 
\begin{Proposition} \label{centralbound} 
Let $\data$ be global torus data, given by $K \subset M_n(F)$
and $g_{\data} \in \adele_K^{\times} \backslash \GL_n(\adele)$. 
Let $\Psi$ be a Schwartz function on $\adele^n$
and $\psi \in \CharK$ a normalized unitary character, $\chi = \psi|_F$. 

Suppose known \eqref{Hyp}. There exists $\beta > 0$,
depending on $n$ and the exponent $\theta$ of \eqref{Hyp}, so that:
for $\Re(s) = n/2$ and any $M \geq 1$, one has 
\begin{multline} \label{mysticgarden} \left|  \frac{|\det g_{\data}|_{\adele}^{s/n}}{\vol(\T_K(F) \backslash \T_K(\adele))} 
 \int_{\T_K(F) \backslash \T_K(\adele)} \psi_{s/n}(t) E_{\Psi}(\chi,s,  tg ) dt  \right| 
\\
\ll_{\Psi,M} \frac{C_{\infty}(\psi_{s/n})^N}{
 C_{\infty}(\chi_s)^{M}}D_{\psi}^{-\beta} D_{\data}^{-\beta}\end{multline}
\end{Proposition}

The main point of this is the decay in $D_{\psi}, D_{\data}$; the reader
should ignore the various factors of $C_{\infty}$, which are a technical matter.
In words,
\eqref{mysticgarden} asserts that varying sequence
of homogeneous toral sets, on $\GL_n$, become equidistributed ``as far as Eisenstein series are concerned,'' if we suppose the pertinent subconvexity hypothesis.

In terms of the discussion of the introduction, it is \eqref{mysticgarden} that proves \eqref{siegel}. \footnote{
\eqref{mysticgarden} also delivers uniformity in the $\psi$-variable; e.g.\ \S \ref{sparseequid} would use this aspect.}

\proof
Let us begin by clarifying 
volume normalizations. 
As in \S \ref{Notation}, we fix an identification $\iota: K \rightarrow F^n$ which 
is an isomorphism for the right $K$-module structures (see \S \ref{ssnf}).  

Define, for each $v$,
local torus data $A_v$ and identifications $\iota_v$ according to the discussion around \eqref{iotavdef}.  
It is important to keep in mind that $A_v$ is {\em not} $(K \otimes F_v)$ but rather its conjugate by $g_v$;
similarly $\iota_v$ is {\em not} simply ``$\iota \otimes_F F_v$'' but is rather twisted by $g_v$.

The discussion of \eqref{hdef} yields elements $h_v \in \GL_n(F_v)$; e.g. for $v$ nonarchimedean, we have:
$\order_{F_v}^n h_v^{-1} = (\order_{A_v})^{\iota_v}$.  

Observe that 
\begin{equation}\label{hv}\prod_v |\det h_v|_v|\det g_{\data,v}|_v \asymp (\disc(K/F))^{1/2}.\end{equation}
To verify \eqref{hv}, note that $\iota: K \rightarrow F^n$ induces
$$\iota_{\adele}: \adele_K = K \otimes_F \adele \rightarrow \adele^n$$
This identification is measure-preserving, because, with our choice of measures, both $\adele/F$ and $\adele_K/K$
have measure $1$. 
In view of the definition \eqref{iotavdef} of $\iota_v$,
this remark implies
$\prod_{v} |\det g_{\data,v}|_v \frac{\iota_{v*} \vol_{A_v}}{\vol_{F^n_v}}  = 1.$
Therefore, taking product of \eqref{hdef2} over all places $v$, 
$$\prod_{v} \frac{\vol(\order_{K_v})}{\vol(\order^n_{F_v})} = \prod_{v} 
|\det g_{\data,v}|_v^{-1} |\det h_v|_v^{-1}$$
which, in combination with \eqref{measure}, yields our claim. 
(Recall that $\order_{F,v}$ and $\order_{K,v}$ are defined, at archimedean places, by the unit balls for suitable norms, cf. \S \ref{subsec:Disc2}).

Without loss of generality $\Psi = \Psi_{\infty} \times \prod_{v} \Psi_v$ where, for each finite place $v$,
$\Psi_v$ is the characteristic function of an $\order_v$-lattice in $F_v^n$.  
(In the general case, one may express $\Psi$ as a sum of such, the implicit
cost being absorbed into the $\ll_{\Psi}$.)

Let $B$ be the union of the following sets of places:
\begin{enumerate}
\item $B_{\infty}$: $v$ is archimedean.
\item $B_{ram}$:  $\disc_v(\psi) \disc_v(\data) > 1$, or $F_v$ is ramified over $\Q$, or the residue field at $v$ has size $\leq n$. 
\item $B_{\Psi}$: $\Psi_v$ does not coincide with the characteristic function of $\order_{F,v}^n$.

\end{enumerate}

Let us note that  
\begin{equation} \label{omit}\exp(|B|) \ll_{F,\Psi} (D_{\psi} D_{\data})^{\varepsilon},\end{equation}
this being a consequence of the fact that the number of prime factors of an integer $N$ is $o(\log N)$.

 We denote by $L^{(B)}$ an $L$-function with the omission of those places inside $B$. Suppose $\Re(s)= n/2$. In view of \eqref{Eis}, 
$|\det g|_{\adele}^{1/2} \left|\int_{\T_K(F) \backslash \T_K(\adele)} E_{\Psi}(\chi,s,tg) \psi_{s/n}(t)\right|$
factors as: 
\begin{multline} \label{above}
\prod_{v} 
|\det g_{\data,v}|_v^{1/2} \left|\int_{K_v} \Psi_{K,v} (t) \psi_{s/n}(t) d^{\times} t \right| \\ = 
| L^{(B)}(K, \psi_{s/n})|  \cdot \prod_{v \notin B}
|\det g_{\data,v} \det h_v|_v^{1/2}
\left| \prod_{v \in B} |\det g_{\data,v}|_v^{1/2} \int_{K_v} \Psi_{K,v} (t) \psi_{s/n}(t) d^{\times} t\right| 
\\ = | L^{(B)}(K, \psi_{s/n})|  \cdot \prod_{v}
|\det g_{\data,v} \det h_v|_v^{1/2}
\left| \prod_{v \in B} |\det h_{v}|_v^{-1/2} \int_{K_v} \Psi_{K,v} (t) \psi_{s/n}(t) d^{\times} t\right| 
\end{multline}

Here $\Psi_{K,v}(x) := \Psi_{v}(x^{\iota} g_{\data,v})$,
so that $\prod_{v} \Psi_{K,v} = \Psi_K$. 
Moreover, we have
used the following evaluation for $v \notin B$:
For such $v$, $\Psi_v$ is the characteristic function of $\order_{F,v}^n$,
and Lemma \ref{mindisc} implies that $(\order_{F,v}^n)^{\iota_v^{-1}}
= \lambda_1 \order_{A_v,v}$
 (some $\lambda_1 \in A_v^{\times}$); because of the definition
\eqref{iotavdef} of $\iota_v$, this means that 
there is $\lambda \in K_v^{\times}$ with $(\lambda \order_{K,v})^{\iota} g_{\data,v}
= \order_{F,v}^n$. 
 Thus
 $x \mapsto \Psi_{v}(x^{\iota} g_{\data,v})$ coincides with the characteristic function
of $\lambda \order_{K,v}$. 
Because, by assumption, both $\psi$  and $K/F$ are unramified at such $v$, 
\begin{equation} \label{localL}|\det g_{\data,v}|_v^{s/n} \left|\int \Psi_{K,v}(t) \psi_{s/n}(t) d^{\times} t \right| =| L_v(K, \psi, s/n)| |\lambda \det(g_{\data,v})|_v^{s/n}\end{equation}
By \eqref{iotavdef} and the discussion preceding \eqref{hdef2}, 
$(\order_{K,v})^{\iota} g_{\data,v}= (\order_{A_v})^{\iota_v} = \order_{F_v}^n h_v^{-1}.$ Comparing this with $(\lambda \order_{K,v})^{\iota} g_{\data,v}
= \order_{F,v}^n$, we deduce that
$|\lambda_v|_v = |\det(h_v)|_v$. The left-hand side of \eqref{localL} thereby
has the same absolute value as 
$| \det g_{\data,v} \det h_v|_v^{s/n} L_v(K, \psi_{s/n})$, concluding
our justification of \eqref{above}.

From the assumed subconvexity bound \eqref{Hyp}, together with \eqref{omit},
\begin{equation} \label{boundA} |L^{(B)}(K, \psi_{s/n})| \ll_{\Psi} 
C_{\infty}(\psi_{s/n})^N (D_K D_{\psi}) ^{1/4- \theta}, \ \ \ \Re(s) = n/2 \end{equation}

 Proposition \ref{localbounds}, with our measure normalizations, 
and after identifying a $K_v$-integral to an $A_v$-integral in the obvious way,
shows for arbitrary $M \geq 1$:
\begin{multline} \label{boundB} \left| \prod_{v \in B} |\det h_v|_v^{-1/2} \int_{K_v} \Psi_{K,v} (t) \psi_{s/n}(t) d^{\times} t \right|
\\
\ll_{M, \Psi,\varepsilon} (D_{\psi} D_{\data})^{\varepsilon}
\begin{cases}C_{\infty}(\chi_s)^{-M} D_K^{-1/4} (D_{\psi} D_K)^{-1/4}, \\
 D_K^{-1/2} (D_{\data}/D_K)^{-\frac{1}{16n^2}}.\end{cases}\end{multline}

The factor $C_{\infty}(\chi_s)^{-M}$ we have interposed on the right-hand side
amounts to ``integrating by parts'' at archimedean places,
before applying Proposition \ref{localbounds}; it will be useful
for convergence purposes later. 
Indeed, suppose $v$ is archimedean. The integral of \eqref{boundB}
is unchanged if we replace $\Psi_v(x)$ by $\psi_{s/n}(\lambda) \Psi_v(\lambda x)$, for $\lambda \in F_v^{\times}$.  On the other hand,
$\psi|_F = \chi$, so $\psi_{s/n}(\lambda) =  \chi_v(\lambda) |\lambda|_v^s \ (\lambda \in F_v^{\times})$. 
Thus if $\nu$ is any probability measure on $F_v^{\times}$,
the integral of \eqref{boundB} is unchanged by the substitution:
$$\Psi_v \mapsto \Psi_v', \Psi_v'(x)= \int_{\lambda} \Psi_v(\lambda x)
\chi_s(\lambda) d\nu(\lambda).$$
Now compare $\|\Psi_v'\|$ and $\|\Psi_v\|$, the norm $\|\cdot\|$
being defined as in the 
 the second statement of Proposition \ref{localbounds}.
An elementary computation shows that, for a smooth 
measure $\nu$, we must have:
\begin{equation}\label{psivprimebound}\|\Psi_v'\| \ll_{\Psi_v} \Cond_v(\chi_s)^{-M}. \end{equation}
Let us note that the implicit constants here depend on higher derivatives of $\Psi_v$; this is permissible.\footnote{ In explicit terms,
\eqref{psivprimebound} for $v$ real amounts to a bound of the type:
$$\int_{x \in \R} \left| \int_{1/2 \lesssim \lambda \lesssim 2} 
  f(\lambda x) |\lambda|^{it} d \lambda \right| \ll_f (1+|t|)^{-M},$$
for a Schwartz function $f$, which is easily verified by integration by parts.}

Combining \eqref{boundA}, \eqref{boundB} and
\eqref{hv}, we see that for $\Re(s) = n/2$:

\begin{multline}\left| |\det g_{\data}|_{\adele}^{1/2} \int_{\T_K(F) \backslash \T_K(\adele)} E_{\Psi}(\chi,s, tg) \psi_{s/n}(t) d^{\times}t  \right|
 \\ 
\ll_{M, \Psi,\varepsilon}
 C_{\infty}(\psi_{s/n})^N (D_{\psi} D_{\data})^{\varepsilon}
 \begin{cases}C_{\infty}(\chi_s)^{-M}
 D_K^{-\theta} D_{\psi}^{-\theta}, \\
 D_{\psi}^{1/4} D_K^{-\theta} (D_{\data}/D_K)^{-\frac{1}{16n^2}}  \end{cases} \end{multline}

Pick $0 < p < 1$.  Using the obvious fact $\min(U,V) \leq U^{p} V^{1-p}$, 
we may replace the right-hand side (ignoring
 $\varepsilon$-exponents) by:
\begin{multline}C_{\infty}(\psi_{s/n})^{N} C_{\infty}(\chi_s)^{-pM} 
D_K^{-\theta} D_{\psi}^{-\theta p + \frac{1-p}{4}}
(D_{\data}/D_K)^{-\frac{1-p}{16n^2}}
\\ \leq 
C_{\infty}(\psi_{s/n})^{N} C_{\infty}(\chi_s)^{-pM} 
 D_{\data}^{-a} D_{\psi}^{-b}  \end{multline}
where $a = \min(\frac{1-p}{16n^2}, \theta), b = \theta p - \frac{1-p}{4}$. For $p$ sufficiently close to $1$;
 these are all positive. Making $M$ arbitrarily large, and using Lemma \ref{lem:CNF}
together with bounds for Dedekind $\zeta$-functions, yields the desired conclusion.

\section{The reaping: a priori bounds} 
\label{explication}

In this section, we translate Proposition \ref{centralbound} into a form
that very explicitly generalizes \eqref{siegel}.

 The result is Proposition \ref{siegelgen}. The work has already been done; this section simply translates between adelic and $S$-arithmetic.

We begin by explicating the connection of the Eisenstein series $E_{\Psi}$ with the classical ``Siegel-Eisenstein'' series, in the case
when the base field is $\Q$. We then carry out the analogue in an $S$-arithmetic setting over an arbitrary base field $F$. 

\subsection{Explications over $\Q$}
The (Siegel)-Eisenstein series on $\PGL_n(\Z) \backslash \PGL_n(\R)$ often appears in the following guise.
Let $f$ be a Schwartz function on $\R^n$. To each $L \in \PGL_n(\Z) \backslash \PGL_n(\R)$
thought of as a lattice $L \subset \R^n$ of covolume $1$, we associate the number 
\begin{equation} \label{Siegel}E_f(L) := \sum_{v \in L-\{0\}} f(v).\end{equation} 

We shall explicate the connection of this construction with the Eisenstein series that we defined previously.

Specialize to the case $F = \Q, \chi =|\cdot|_{\adele}^s$. We take $\Psi_v$ to coincide with the characteristic function of
$\Z_v^n$ for all finite $v$, and $\Psi_{\infty} =f$. 
Define $f_{s}$ on $\R^n$ via the rule
$$f_{s} (\mathbf{v}) = \int_{t \in \R^{\times}} |t|^{s} f(\mathbf{v} t) d^{\times}t.$$
Then $f_s$ satisfies the transformation property $f_{s}(\lambda \mathbf{v}) = |\lambda|^{-s} f_s(\mathbf{v})$.
Taking into account the fact that the natural map $\R_{>0} \times \prod_{p} \Z_p^{\times} \rightarrow \adele_{\Q}^{\times}/\Q^{\times}$ is a homeomorphism, 
we see that for $g_{\infty} \in \GL_n(\R)$, 

$$|\det g_{\infty}|^{s/n} E_{\Psi}(|\cdot|^s, g_{\infty}) = |\det g_{\infty}|^{s/n} \sum_{\mathbf{v} \in \Z^n. g_{\infty}} f_s(\mathbf{v})$$

Note that, by Mellin inversion, $f =  \int_{s} f_s ds$, where the $s$-integration is taken over a
line of the form $\Re(s) = \sigma \gg 1$.  Consequently, 
the Siegel-Eisenstein series \eqref{Siegel} corresponds to the function 
on $\PGL_n(\Q) \backslash \PGL_n(\adele)$ defined by 
$$g \mapsto \int_{\Re(s) = \sigma \gg 1} |\det g|_{\adele}^{s/n} E_{\Psi}(s, g) ds.$$

\subsection{$S$-arithmetic setup.} \label{Sarithsetup}
We revert to the general setting of a number field $F$. 
We shall henceforth pass from an adelic setup, to an $S$-arithmetic setup.

Fix, therefore, a finite set of places $S$, containing all archimedean ones.
Set $F_S = \prod_{v \in S} F_v$. We assume that $S$ is chosen
so large that $$\adele^{\times} = F^{\times} F_S^{\times}  \prod_{v \notin S} \order_{F,v}^{\times}.$$ 

Under these assumptions, we may identify $\adele^{\times}/F^{\times} \prod_{v \notin S} \order_{F,v}^{\times}$ 
to $F_S^{\times}/\order^{\times}$, where $\order := F \cap \prod_{v \notin S} \order_{F,v}$. By Dirichlet's theorem, the quotient of $\order$ by roots of unity
comprises a free abelian group of rank $|S| - 1$. 

Similarly, making use of the strong approximation theorem for the group $\SL_n$, we can identify $$\PGL_n(F) \backslash \PGL_n(\adele)/ \prod_{v \notin S} \PGL_n(\order_{F,v})$$
 to the quotient $$\PGL_n(\order) \backslash \PGL_n(F_S).$$
If $\mu$ is a measure on $\PGL_n(F) \backslash \PGL_n(\adele)$,
we shall often abuse notation and identify $\mu$
with the projected measure on $\PGL_n(\order) \backslash \PGL_n(F_S)$. 

\subsection{The $S$-arithmetic Eisenstein series.}
It will take us a little work to unravel the Eisenstein series into a form
which we can easily use in the $S$-arithmetic case.  There will be some complications arising from the
 failure of strong approximation for $\PGL_n$.

Let $F_S^{(1)}$ consist of those elements of $F_S^{\times}$
with $|x| = 1$. Then $F_S^{\times}$ can be identified with $F_S^{(1)} \times \R_{>0}$. In this way, we can identify a character of the compact group
$F_S^{(1)}/\order^{\times}$, to a character of $F_S^{\times}/\order^{\times}$:
by extending trivially on $\R_{>0}$.  Thus
the group of all normalized characters of $C_F$, unramified away from $S$,
is identified to the character group of $F_S^{(1)}/\order^{\times}$. 

We have already normalized measures on $F_S^{\times}$.
We normalize the measure on $F_S^{(1)}$ by differentiating
along fibers
of the map $F_S^{\times} \rightarrow \R_{>0}$ given by $x \mapsto |x|$; here we equip
$\R_{>0}$ with the measure $\frac{d\lambda}{\lambda}$.
The measure of $F_S^{(1)}/\order^{\times}$ then equals the measure
of $\adele^{(1)}/F^{\times}$. 

We now construct the $S$-arithmetic version of the Eisenstein series. 
Let $\Psi_S$ is a Schwartz function on $F_S$; put
$\Psi = \Psi_S \times \prod_{v \notin S} 1_{\order_v^n}$. 
Let $\chi$ be a character of $F_S^{(1)}/\order^{\times}$, identified, via the remarks above, with a character of $\adele^{\times}/F^{\times}$. 

Set
$$\bar{E}_{\Psi}(\chi,g)  = 
\int_{\Re(s) = 2n}
E_{\Psi}(\chi,s, g) |\det g|_{\adele}^{s/n} \frac{ds}{2 \pi i}.$$ Then, for $g \in GL_n(F_S)$,  $E_{\Psi}(\chi,g)$ equals:
 \begin{multline}
 \int_{t \in \adele^{\times}/F^{\times}:
|t|_{\adele}^n = |\det g|_{\adele}} \chi^{-1}(t)
 \sum_{\mathbf{v} \in F^{n} - \{0\}} \Psi(\mathbf{v} t^{-1} g ) 
  \\ = \int_{t \in F_S^{\times}/\order^{\times}:
|t|_{F_S}^n = |\det g|_{F_S}} \chi^{-1}(t)
 \sum_{\mathbf{v} \in \order^n - \{0\}} \Psi_{S}(\mathbf{v} t^{-1} g ) \end{multline}
In the first expression, the $t$-integral is taken w.r.t. the measure that is transported
from the measure on $\adele^{(1)}/F^{\times}$;
in the second expression, the $t$-integral here is taken over a compact abelian Lie group
of dimension $|S|-1$, with respect to the measure previously discussed.

Fix $a \in F_S^{(1)}$ and put, for $g \in\GL_n(F_S)$, 
\begin{multline} \label{Sarithmetic} \bar{E}_{\Psi, a}(g) = 
\frac{1}{\vol(F_S^{(1)}/\order^{\times})}\sum_{\chi: F_S^{(1)}/\order^{\times} \rightarrow S^1}
\chi(a \det g) \bar{E}_{\Psi}(\chi^n,g) \\ = \sum_{\stackrel{t^n = a \det g}{t \in F_S^{\times}/\order^{\times}}}
\sum_{\mathbf{v} \in \order^n - \{0\}} \Psi_S(\mathbf{v} t^{-1} g) \end{multline}
Indeed, the $\chi$-sum restricts to those $t \in F_S^{\times}/\order^{\times}$ so that
$t^n$ and $a \det g$ differ by an element of $\R_{>0}$;
however, since $|t^n|_{F_S} = |a \det g|_{F_S}$, this forces $t^n = a \det g$. 
The $t$-sum is finite, for the quotient $F_S^{(1)}/\order^{\times}$
is a compact abelian Lie group. 

The function $\bar{E}_{\Psi, a}(g)$ defines a function on $\PGL_n(\order) \backslash \PGL_n(F_S)$; it is the $S$-arithmetic version of our degenerate Eisenstein series. 
\subsection{Bounding the mass of $\bar{E}_{\Psi,a}$.}
In order to bound the $\mu_{\data}$-measures of functions of the
type $\bar{E}_{\Psi,a}$, via Proposition \ref{centralbound},
we require a subconvex bound for the $L$-functions $L(K, \psi)$,
when $\psi$ is pulled back from a fixed character of $F$ via the norm map, i.e.
we require the following estimate for $\Re(s) = 1/2$:
\begin{equation}\label{Hypmod}|L(K,\chi\circ \Norm_{K/F}, s)| \ll
\left( C_{\infty}(\chi_{s})D_{\chi}\right)^N  D_K^{1/4-\eta}, \ \chi \in \Omega(C_F) \mbox{ unitary.} \end{equation}
for some constants $N,\eta>0$ which depend at most on $F$ and $n=[K:F]$ . The quantity $C_{\infty}(\chi)$
is defined as in \refs{Cinftydef}. Notice that unlike \refs{Hyp}, we do not require in \refs{Hypmod} a subconvex bound in the 
``$D_{\chi}$-aspect''.

The bound \eqref{Hypmod} is known in more cases than 
\eqref{Hyp}. For instance, it is known\footnote{even in the $D_{\chi}$-aspect} when $F=\Q$ and $[K:\Q]\leq 3$ (cf. Appendix \ref{subconvex-a}). 
\begin{Proposition}\label{siegelgen}
Suppose \eqref{Hypmod} is known. Then, for homogeneous toral data $\data$, 
\begin{equation} \label{mu} \left| \mu_{\data}( \bar{E}_{\Psi,a}) - 
\delta \int_{\adele_F^n} \Psi  \right| \ll_{\Psi}  \disc(\data)^{-\question}\end{equation}
for some $\beta>0$  -- depending 
on $n$ and the exponent $\eta$ of \eqref{Hypmod};
and where $\delta=\delta_{a} \in \Q$ belongs to a finite set of rational numbers,
depending on $S, F$ and $n$.
\end{Proposition}
The $\delta$ arise
from  the ``connected components.'' We could be a little more precise about their value, but there is no point. 
This result implies the generalization, to an arbitrary base field $F$
and an $S$-arithmetic setting, of \eqref{siegel}, discussed in the Introduction. 
\proof 
The data $\data$ is defined by a field $K \subset M_n(F)$
and an element $g_{\data} \in \GL_n(\adele)$.

Recall
\begin{multline} \bar{E}_{\Psi}(\chi,g)  = 
\int_{\Re(s) \gg 1}
E_{\Psi}(\chi,s, g) |\det g|_{\adele}^{s/n} \frac{ds}{2 \pi i}, \\
\bar{E}_{\Psi, a}(g) = 
\frac{1}{\vol(F_S^{(1)}/\order^{\times})}\sum_{\chi: F_S^{(1)}/\order^{\times} \rightarrow S^1}
\chi(a \det g) \bar{E}_{\Psi}(\chi^n,g) 
\end{multline}

As we have already commented, we are going to identify $\mu_{\data}$ with its projection 
to $\PGL_n(\order) \backslash \PGL_n(F_S)$. This being so, let us consider
$\mu_{\data}(\bar{E}_{\Psi,a})$. 
Shift contours to $\Re(s) = n/2$ in the defining integrals
and apply the bounds of Proposition \ref{centralbound}. 
(There are no concerns with convergence; the support of $\mu_{\data}$ is compact.)
 The function $E_{\Psi}(\chi^n,s, g)$
has a pole (by Proposition \ref{meromorphic}) precisely when $\chi^n$ is the trivial character and $s \in \{0,n\}$; moreover,
 Proposition \ref{meromorphic} computes the residue in those cases. 

The result is:
\begin{equation} \label{blah2}\mu_{\data}(\bar{E}_{\Psi,a}) =   \sum_{
\stackrel{\chi: F_S^{(1)}/\order^{\times} \rightarrow S^1}{
\chi^n= 1}} \mu_{\data}
(\chi(a \det g))  \int_{\adele^n} \Psi(\mbfx) d\mbfx  +  \mathrm{Error},\end{equation}
(note that $g \mapsto \chi(a \det g)$ indeed defines a function on $\PGL_n(\order) \backslash \PGL_n(F_S)$), and:
\begin{multline}\label{Errorbound}|\mathrm{Error}|
\ll
\frac{1}{\vol(F_S^{(1)}/\order^{\times})}\times\\
\sum_{\chi: F_S^{(1)}/\order^{\times} \rightarrow S^1} \max_{\Re(s) = 1/2}(1+|s|)^2  \left|
\mu_{\data}\left(
E_{\Psi}(\chi^n,{s}, \cdot) \chi_{s/n} \circ \det
\right)\right|.\end{multline}

For $\chi$ as above, put  $\psi := \chi \circ N_{K/F}$, 
a character of $C_K$. We have
$$|\mu_{\data}\left(
E_{\Psi}(\chi^n,s, \cdot) \chi_{s/n} \circ \det
\right)| = |\det g_{\data}|_{\adele}^{s/n}
\left| \frac
{\int_{\T_K(F) \backslash \T_K(\adele)} \psi_{s/n}(t) E_{\Psi}(\chi^n,s, 
 tg_{\data} ) dt }
{\vol(\T_K(F) \backslash \T_K(\adele))}  \right|.$$
By (the proof of) Proposition \ref{centralbound}, we have under \refs{Hypmod}  for any $M>1$ and some $\beta,N>0$
$$|\mu_{\data}\left(E_{\Psi}(\chi^n,s, \cdot) \chi_{s/n} \circ \det
\right)| \ll_{M,\Psi} (C_{\infty}(\chi^n_s))^{-M}D_{\chi}^ND_{\data}^{-\beta}.$$

We have utilized the notation $\chi^n_s := (\chi^n)_s$, the character
$x \mapsto \chi(x)^n |x|^s$. 
We have used the fact that, in the present context, $C_{\infty}(\chi^n_s)$ and
$ C_{\infty}(\psi_{s/n})$ are bounded within powers of each other.

 Taking $M$ large enough, we obtain the following bound:\footnote{The number of possibilities for $\chi|_{F_{\infty}}$ is, in general, infinite. However, 
if  $\chi$ contributes non-trivially to $\mathrm{Error}$, 
then $D_{\chi}$ is bounded above depending on $\Psi$. 
Moreover, the number of such $\chi$s with $\inf_{\Re(s) = 1/2} C_{\infty}(\chi^n_s) \leq T$ is bounded polynomially in $T$. }
 $$\mbox{Error}\ll_{\Psi}D_{\data}^{-\beta}.$$ 

Let us now analyze the right-hand side of \eqref{blah2}. 
The set of elements of $C_F$ of the form $\Norm_{K/F}(x) . \det g_{\data}$
(for some $x \in C_K$)
is a coset of a subgroup of $C_F$ of finite index. 
It projects to a coset of a subgroup of $F_S^{\times}/\order^{\times}$
which contains the $n$th powers, which we may identify
with a coset of a subgroup of the finite group $F_S^{\times}/\order^{\times} (F_S^{\times})^n$.   Call this subgroup $Q$ and the
 pertinent coset $bQ$.  
If $\chi: F_S^{\times}/\order^{\times}\rightarrow S^1$ is so that $\chi^n=1$, then $$\mu_{\mathscr{D}}(\chi(a \det g))=\frac{1}{|Q|} \sum_{x \in Q} \chi(a b x).$$ 
Therefore, the coefficient of $\int \Psi$, on the right-hand side 
of \eqref{blah2}, is given by:
$$ \sum_{\chi^n = 1} \frac{1}{|Q|} \sum_{x \in Q} \chi(a b x)$$
In particular, this lies
in a finite set of rational numbers. \qed

\section{Proof of Theorem \ref{RealTheorem} and Theorem \ref{Volumes}}
\label{Proofs}

Let $F$
be a number field, and let $\mathscr{D}_i$ be a sequence of homogeneous toral data on $\PGL_n$ over $F$.  Let $Y_i$ be the associated homogeneous toral sets, and $\mu_{i}=  \mu_{\mathscr{D}_i}$ the corresponding probability measures. 
Recall that $\mathscr{D}_i$ is defined
by a torus $\mathbf{T}_i \subset \PGL_n$ together with $g_i \in \PGL_n(\adele)$. Let $K_i$ be the corresponding (degree $n$) field extensions.
See \S \ref{HsetsGLn}. 

When convenient we may drop the subscript $i$, referring simply to $\mathscr{D}, \mu_{\mathscr{D}}, K, \mathbf{T}$, etc. 

\subsection{Bounds on the mass of small balls and the cusp.}

Fix a set of representatives $1=a_1, \dots, a_r \in F_S^{(1)}$ for  
$F_S^{\times}/(F_S^{\times})^n \order^{\times}$.  Such representatives
may indeed be chosen in $F_S^{(1)}$.

Take $x \in \PGL_n(\order) \backslash \PGL_n(F_S)$; 
we say that a lattice in $F_S^n$ (i.e. an $\order$-submodule of rank $n$)
{\em corresponds to $x$} if it is of the form
$\order^n . g_0 t_0^{-1}$
where $g_0$ is a representative for $x$ in $\GL_n(F_S)$;
and  $t_0 \in F_S^{\times}/\order^{\times}$ is so that $t_0^n =
\det(g_0) a_i$, some $1 \leq i \leq r$. 
There are only finitely many lattices corresponding to a given $x$.

We say a set $\mathcal{Q}$ in $F_S^n$ is nice if 
there exists a Schwartz function $\Psi_S$
on $F_S^n$ so that $\Psi_S \geq 1$ on $\mathcal{Q}$ with $\int \Psi_S \leq 2 \vol(\mathcal{Q})$.

\begin{Lemma} \label{refinement}
Let $\mathcal{Q} \subset F_S^n$ be nice; set
$\mathcal{L}_{\mathcal{Q}} \subset \PGL_n(\order) \backslash \PGL_n(F_S)$ to comprise those
$x\in \PGL_n(\order) \backslash \PGL_n(F_S)$ so that
a lattice corresponding to $x$ contains 
an element of $\mathcal{Q}$. 
 
If \eqref{Hypmod} is known, then 
  \begin{equation} \label{eq:refinement}
 \mu_{\mathscr{D}}(\mathcal{L}_{\mathcal{Q}}) 
 \ll_{F,S,n} \vol(\mathcal{Q}) + 
O_{\mathcal{Q}}(\disc(\data)^{-\question}) 
\end{equation}
\end{Lemma}
\proof 
Indeed, choose $\Psi_S$ as remarked. As is clear from \eqref{Sarithmetic},  the function $\sum_{i=1}^{r} \bar{E}_{\Psi, a_i}$
dominates the characteristic function of $\mathcal{L}_{\mathcal{Q}}$.  The result is a consequence of Proposition \ref{siegelgen}. 
\qed 

Let $\epsilonvector = (\epsilon_v)_{v \in S}$ be a choice of $\epsilon_v \in (0,1)$ for each $v \in S$. Set $\|\epsilonvector\| = \prod_{v \in S} |\epsilon_v|_v$. 
For each $v \in S$,
let $B_v(\epsilon_v)$ be an $\epsilon_v$-neighbourhood of the identity
in $\PGL_n(F_v)$.   Here, we equip $\PGL_n(F_v)$ with the metric
that arises from the adjoint embedding $\PGL_n \hookrightarrow M_{n^2}$;
and we equip $M_{n^2}$ with the metric that arises from norm: supremum of all matrix entries.

Let $$B_S(\epsilonvector) = \prod_{v \in S}  B_v(\epsilon_v) \subset \PGL_n(F_S).$$
The following is a consequence of Lemma \ref{refinement}, for a suitable choice of $\mathcal{Q}$; we leave the details to the reader. 
\begin{Lemma}\label{balllemma}
(Bounds for the mass of small balls).  
Suppose \eqref{Hypmod} is known. 
 Let $x_0 \in \PGL_n(\order) \backslash \PGL_n(F_S)$. Then
 $$\mu_{\mathscr{D}} (x_0 B_S(\epsilonvector)) \ll_{F,S,n} \|\epsilonvector\|^n + 
O_{\epsilonvector}(\disc(\data)^{-\question}) $$
Moreover, the implicit constant in $O_{\epsilonvector}(\dots)$ is bounded
uniformly when $x_0$ belongs to any fixed compact. 

 \end{Lemma}

Now Let $N_{0,v}$ be the standard norm on $F_v^n$ (cf. \S \ref{local}). 
  For $g \in  \GL_n(\adele)$, we set
 $$\height(g)^{-1} =   |\det(g)|_{\adele}^{-1/n} \inf_{\lambda \in F^n - \{0\}}
  \prod_{v} N_{0,v}(\lambda g_v).$$  

 If $K_{\max} = \prod_{v} \mathrm{Stabilizer}(N_{0,v})$ is the standard maximal compact subgroup of $\GL_n(\adele)$, and $\bar{K}_{\max}$ its image
  in $\PGL_n$, then $\height$ descent to a proper map from\ $\PGL_n(F) \backslash \PGL_n(\adele)
  /\bar{K}_{\max}$ to $\R_{>0}$.  In particular, sets of the type $\height^{-1}\big([R,\infty)\big)$,
for large $R>0$, define ``the cusp.''
 
 The next result is again a consequence of Lemma \ref{refinement}.
  \begin{Lemma} \label{escapelemma}
(Bounds for the cusp.) Suppose \eqref{Hypmod} is known. 
  $$\mu_{\data}(\height^{-1}[R, \infty)) \ll R^{-n}+O_{R}(\disc(\data)^{-\question})$$
  \end{Lemma}

\subsection{Proof of Theorem \ref{Volumes}.}
 The volume of the homogeneous toral set associated to $\mathscr{D}$
is defined in \eqref{volume}. We take the set $\Omega_0$
to be the product $\prod_{v|\infty} \Omega_v \times
\prod_{v \, \operatorname{finite}} \PGL_n(\order_{F,v})$. 
Here, we set $\Omega_v \subset \PGL_n(F_v)$,
for $v$ archimedean, to equal the image, in $\PGL_n$, of $\exp(E_v)$; here
$$E_v:=\{Y \in M_n(F_v): \mbox{$Y$ has operator norm $\leq \frac{1}{10}$}\}$$ 
and the operator norm is taken w.r.t. the canonical norm on $F_v^n$. 
Let us note that $\exp: E_v \rightarrow \exp(E_v)$ is a diffeomorphism
onto its image, being inverted by $\log$.

We claim that
\begin{equation}  \label{measurechoice}
\log \vol(Y) = \sum_{v} \log \vol\{t \in \T(F_v):
g_{v}^{-1} t g_{v} \in \Omega_0\}^{-1}  + o(\log \disc(\data)),\end{equation}
where, on the right hand side, the measure on $\T(F_v)$
is normalized as indicated in \S \ref{Notation}; and we understand
$\T(F_v)$ as being embedded in $\T(\adele)$ in the natural way. 

To verify \eqref{measurechoice} we need to consider our measure normalizations.
In the definition of ``$\vol(Y)$'', in \S \ref{Voldef},
we endowed $\T(F) \backslash \T(\adele)$ with a probability measure.
If we normalize the measures on $\T(F_v)$ according to \S \ref{Notation},
the product measure is not a probability measure
on $\T(F) \backslash \T(\adele)$; its mass is given by
Lemma \ref{lem:CNF} to be a certain $\zeta$-value. 
By a result of Siegel, $\log |\zeta_{K}(1)| = o(\log \disc K)$,
and, by Lemma \ref{mindisc}, $\disc(K) \ll \disc(\data)$. 
This establishes \eqref{measurechoice}.

Let $v$ be a finite place. Let us recall 
that we defined an order $\Lambda_{v} \subset g_v^{-1} K_{v} g_v$
via $\Lambda_{v} = g_v^{-1} K_{v} g_v \cap M_n(\order_{F,v})$ (see 
\S \ref{subsec:thediscriminant}). 
Therefore, $\Lambda_{v}^{\times} = g_v^{-1} K_{v} g_v \cap \GL_n(\order_{F,v})$. 
Thus
$\{t \in \T(F_v): g_v^{-1} t g_v \in \Omega_0 \}$ is identified,
via $t \mapsto g_v^{-1} t g_v$, to $\Lambda_{v}^{\times} F_v^{\times}/F_v^{\times}$.
Since $\Lambda_{v} \cap F_v = \order_{F,v}$, 
we see, given the measure normalizations of \S \ref{Notation}, 
that $\vol\{t \in \T(F_v): g_v^{-1} t g_v \in \Omega_0 \}
= \vol(\Lambda_{v}^{\times})/\vol(\order_{F,v}^{\times})$. 
Taking into account
Lemma \ref{dvexpl}, Lemma \ref{localvolumes} and \eqref{measure}, this becomes:
\begin{equation}\label{finitevol}\log \vol\{t \in \T(F_v): g_v^{-1} t g_v \in \Omega_0 \}
=\frac{1}{2} \log \disc(\data_v) + O_F(1), \ \ \
\mbox{$v$ finite.}\end{equation}
Here the error term $o_F(1)$ is identically zero if the residue characteristic of $v$
is larger than $n$, $F$ is unramified at $v$,  and $\disc_v(\data_v)=1$. 

Now we establish an (approximate) analog of \eqref{finitevol} at archimedean places. 
We claim that, for archimedean $v$, 
\begin{multline} \label{fredy} \log \vol\{t \in \T(F_v): g_v^{-1} t g_v \in \Omega_0 \}
=  \log \vol(\Lambda_{v}) + O_F(1) \\ = \frac{1}{2} \log \disc_v(\data_v) + O_F(1)\end{multline}
The second equality follows readily from Lemma \ref{dvexpl}, so we need
to verify the first equality. 
Put $A_v = g_v^{-1} (K \otimes F_v) g_v \subset M_n(F_v)$.
Because $A_v$ is a {\em subalgebra}, we have
$\exp(A_v) \subset A_v$, and $\log(A_v) \subset A_v$
when $\log$ is defined. Therefore, 
$$A_v \cap \exp(E_v) = \exp(A_v \cap E_v)$$

The set $\{t \in \T(F_v): g_v^{-1} t g_v \in \Omega_0\}$
is identified, via $t \mapsto g_v^{-1}  t g_v$, to
$F_v^{\times} \exp(A_v \cap E_v)/F_v^{\times}$.
Its measure is therefore easily seen to be bounded above and below
by constant multiples by the $A_v$-measure of
$\exp(A_v \cap E_v) . \{F_v^{\times} \cap \exp(E_v)\}$. 
This set contains $\exp(A_v \cap E_v)$ and is contained in
$\exp(A_v \cap 2 E_v)$.  The measure of $\exp(A_v \cap E_v)$
and $\exp(A_v \cap 2 E_v)$ are bounded above and below by constant multiples of the volume of $\Lambda_{v}$
by constants, for the map $\exp: A_v \rightarrow A_v^{\times}$ preserves (up to a constant) measure. 
This establishes \eqref{fredy}.

Combining \eqref{measurechoice}, \eqref{finitevol} and \eqref{fredy},
we see that
$$\log \vol(Y) = \frac{1}{2} \sum_{v} \log \disc(\data_{v}) + o_F(\log 
\disc \data).$$
The conclusion of Theorem \ref{Volumes} follows. 
\qed

\subsection{Proof of Theorem \ref{RealTheorem}.}
Let $v$ be a place as indicated in the proof of Theorem \ref{RealTheorem}.
Let $S$ be a finite set of places of $F$ as in \S \ref{Sarithsetup};
enlarging $S$, we may suppose $v \in S$ without loss of generality.

Let $H_i = g_{i,v}^{-1} \mathbf{T}_i(\Q_v) g_{i,v}$. 
The measure $\mu_i:= \mu_{\mathscr{D}_i}$, upon projection
to \linebreak$\PGL_3(\order) \backslash \PGL_3(F_S)$, 
is invariant under $H_i$.  Denote by $\barmu_i$ this projection.

Let $\barmu_{\infty}$ be any $\mathrm{weak}^*$ limit of the measures $\barmu_{\mathscr{D}_i}$,
which we may assume is the projection of a limit $\mu_{\infty}$ of the original sequence.
Because the bounds of Lemma~\ref{escapelemma} are uniform in $\mathscr{D}_i$,
the measure $\barmu_{\infty}$ is a probability measure.

It will suffice to show that 
$\barmu_{\infty}$  is $\SL_3(F_v)$-invariant. In fact, it then follows that $\barmu_\infty$
is $\SL_3(F_S)$-invariant by irreducibility of the lattice $\PGL_3(\order) \subset \PGL_3(F_S)$ and,
 $S$ being arbitrary, that $\mu_\infty$ is $\SL_3(\adele)$-invariant. 
 Once $\mu_{\infty}$ is $\SL_3(\adele)$-invariant, it is determined by its projection
 to the compact group $\adele_F^{\times}/F^{\times} (\adele_F^{\times})^3$. 
 We are reduced to showing that limits of homogeneous measures on a compact
 abelian group are of the same type, which is easy. 
 
We shall use the following observation,
which is a consequence of Lemma~\ref{refinement}:
Let $\mathbf{P} \subset \PGL_3$ be the stabilizer of a line in $F^3$;
let $\mathbf{P}'$ be the stabilizer of a plane in $F^3$. 
Thus $\mathbf{P}, \mathbf{P}'$ are $F$-parabolic subgroups.
Let $Z \subset \mathbf{P}(F_w) \backslash \PGL_3(F_w), Z' \subset \mathbf{P}'(F_w) \backslash \PGL_3(F_w)$ be the $F_w$-points of varieties of dimension $\leq 1$.  Then:
\begin{equation} \label{UB}\barmu_{\infty} (
\left(\mathbf{P}(\order) \backslash \mathbf{P}(F_S)\right) . Z) = 0,\ 
\barmu_{\infty}\left(\mathbf{P'}(\order) \backslash \mathbf{P}'(F_S). Z' \right) =0. 
\end{equation}
Indeed, the first assertion of \eqref{UB} follows directly from
Lemma \ref{refinement}, taking for $\mathcal{Q}$ a suitable sequence of sets. The second assertion
may be deduced from the first by applying the outer automorphism (transpose-inverse) of $\PGL_3$
to the entire situation. 

{\em Case 1.} Suppose $\disc_v(\mathscr{D}_i) \rightarrow \infty$.
Let $\mathfrak{h}_i = \mathrm{Lie}(H_i)$;
let $\mathfrak{h}_{\infty}$ be any limit of $\mathfrak{h}_i$
inside the Grassmannian of $\mathfrak{pgl}_3$.  It
is a 2-dimensional commutative Lie subalgebra.  The
measure $\barmu_{\infty}$ is invariant by $\exp(\mathfrak{h}_{\infty})$. 
Necessarily $\mathfrak{h}_{\infty}$ contains a nilpotent element.

Identify $\mathfrak{pgl}_3$ with trace-free $3 \times 3$ matrices. 
There are two conjugacy classes of nontrivial nilpotents in $\mathfrak{pgl}_3$
according to the two possible Jordan blocks. Suppose that $\mathfrak{h}_\infty$ contains 
a conjugate of $\left(\begin{array}{ccc} 0 & 1 & 0 \\ 0 & 0 & 1\\ 0 & 0 & 0 \end{array}\right)$
(i.e.\ a generic nilpotent element). Then since the centralizer of this generic element is two dimensional,
it follows  in this case by commutativity of $\mathfrak{h}_{\infty}$
that $\mathfrak{h}_\infty$ is this centralizer. I.e.\  $\mathfrak{h}_{\infty}$ contains
in any case a conjugate  $\mathfrak{n} $ of the Lie algebra spanned by 
$\left( 
\begin{array}{ccc} 0 & 0 & 1 \\ 0 & 0 & 0 \\ 0 & 0 & 0 \end{array}\right) $.

We make the following observation:
Suppose $\mathfrak{j}$ is a proper Lie subalgebra 
of $\mathfrak{sl}_3$ containing $\mathfrak{n}$. 
Then $\mathfrak{j}$ is reducible over $\bar{F}_v^3$,
i.e.\ fixes a line or a plane over the algebraic closure.
Indeed, the only proper Lie subalgebra of $\mathfrak{sl}_3$
that acts irreducibly over the algebraic closure
is $\mathfrak{so}_q$, for $q$ a nondegenerate quadratic form. 
But $\mathfrak{so}_q$ does not contain any conjugate of $\mathfrak{n}$.

By \cites{MT, Ratner}, $\barmu_{\infty}$ may be expressed as
a convex linear combination of Haar probability measures $\mu_\iota$ on closed orbits $x_\iota H_\iota$
with $\iota$ belonging to some probability space $I$;
here $H_\iota$ is a closed subgroup of $\PGL_3(F_S)$.
  Moreover, all the measures $\mu_\iota$ are ergodic
under the action of $N=\exp(\mathfrak{n})$.

Suppose $\barmu_{\infty}$ is not $\SL_3(F_v)$-invariant. Then for a positive proportion of 
the $\iota$, say for $\iota\in I'$,
the subgroup $H_\iota$ does not contain $\SL_3(F_v)$. Therefore, 
$\barmu_{\infty}$ dominates the convex combination:
$$\barmu_{\infty} \geq \int_{I'} \mu_\iota.$$

Fix some $\iota\in I'$ and let $x_\iota=\PGL_3(\order) g$,
where we may assume that $x_\iota$ has dense orbit in $x_\iota H_\iota$
under the action of $N$. We claim:
\begin{equation}\label{claim}\mbox{There exists a proper $F$-subgroup $\mathbf{J} \subset \PGL_3$ so that $x_{\iota} H_{\iota} \subset
(\mathbf{J}(\order) \backslash \mathbf{J}(F_S)). g$.}\end{equation}   

For the proof of the claim, consider first $H_\iota$ as a subgroup of a product of real
and $p$-adic Lie groups: if $\bar{S}$
is the set of places of $\Q$ below $S$, then we consider
 $\PGL_3(F_S) = \prod_{p \in \bar{S}} \prod_{w|p, w \in S}
\PGL_3(F_w)$. 
It may be seen that, in a neighbourhood of the identity $H_\iota$ is itself a product of real and $p$-adic subgroups. 
We {\em define} the Lie algebra $\mathfrak h$ of $H_\iota$
to be the product of the real Lie algebra and the various $p$-adic ones;
this is, by definition, a $\Q_{\bar{S}}$-submodule
of $\oplus_{w \in S} \mathfrak{pgl}_3(F_w)$.  (Here, and in what follows,
we use $\mathfrak{pgl}_3(k)$ to denote the $k$-points of the vector space $\pgl_3$.)

In general, if the map $S \rightarrow \bar{S}$ is not bijective,
$\mathfrak{h}$ may not be a direct sum of its projections to the
$\mathfrak{pgl}_3(F_w)$. 
We claim, however, that the projection of $\mathfrak{h}$
to $\mathfrak{pgl}_3(F_v)$ is a proper subalgebra. 
Indeed,  $\mathfrak n$ is a Lie subalgebra 
of $\mathfrak h \cap \pgl_3(F_v)$; 
it follows that all conjugates of $\mathfrak n$
by elements of $H_\iota$ again belong to $\mathfrak h \cap \pgl_3(F_v)$. 
Were the projection of $\mathfrak{h}$ to $\mathfrak{pgl}_3(F_v)$
surjective, it would follow -- by the simplicity of $\pgl_3(F_v)$ as a module over itself -- that $\pgl_3(F_v)$ is contained in $\mathfrak h$; contradiction.

Next let $\mathbf{J}'$ be the Zariski closure of $gH_{\iota}g^{-1}\cap \PGL_3(\order)$; by definition, this is an $F$-algebraic subgroup of $\PGL_3$. 
$\mathbf{J}'$ preserves the projection of $\Ad(g) \mathfrak{h}$ to $\mathfrak{pgl}_3(F_v)$,
and is therefore a proper subgroup of $\PGL_3$. 

Just as in the Borel density theorem it follows that $g N g^{-1}$ is contained in $\mathbf{J}'(F_v)$.
In fact, by Chevalley's theorem there is an algebraic representation $\phi$ of $\PGL_3$ on $V$ and a vector $v_\phi\in V_\phi$
such that $\mathbf{J}'$ is the stabilizer of the line generated by $v_\phi$. Fix some
 parametrization $n_t$ of $N$ as a one-parameter unipotent group.
Note that the line spanned by $\phi(gn_tg^{-1})(v_\phi)$ approaches the line spanned by an
 eigenvector $v_N$ of $\phi(gNg^{-1})$
if $|t|\to\infty$.
 By our assumption on $x_\iota$ we have a sequence $t_k\in F_v$ with $|t_k|_v\to\infty$
for which $x_\iota n_{t_k}\to x_\iota$ as $k\to\infty$. Therefore, there exists some sequence
$\gamma_k\in\PGL_3(\order)\cap gH_\iota g^{-1}$ and $g_k\in\PGL_3(F_v)$ with $ g n_{t_k} =\gamma_k g g_k  $ and 
$g_k$ approaching the identity.
This implies that $\phi(g n_{t_k}^{-1}g^{-1})(v_\phi)=\phi( gg_k^{-1} g^{-1}\gamma_k^{-1}) (v_\phi)$
both approaches $v_N$ and $v_\phi$, i.e.\ that $v_N=V_\phi$ and so $gNg^{-1}\subset\mathbf{J'}(F_v)$. 

To prove \eqref{claim}, we proceed as follows. Let $\mathbf{J}''\subset\mathbf{J}'$
 be the preimage, in $\mathbf{J}'$, of the commutator subgroup of $\mathbf{J}'/R_u(\mathbf{J}')$. Since $\mathbf{J}''$
is $F$-algebraic without $F$-characters, it follows that $\PGL_3(\order)\mathbf{J}''(F_S)$
is closed. Therefore, the same holds for $\PGL_3(\order)\mathbf{J}''(F_S)g$ which is invariant
under $N$. We see that $\PGL_3(\order)\mathbf{J}''(F_S)g$ contains $x_\iota H_\iota$ by our choice of $x_\iota$.
We can take $\mathbf{J} := \mathbf{J}''$.

Next we claim that $\mathbf{J}$ is contained in an $F$-parabolic subgroup $\mathbf{P}$. 
For this we need to show that $\mathbf{J}$ preserves a line in $F^3$ or a line
in the dual $(F^3)^*$.  Indeed, as noted before
$\Lie \mathbf{J} \otimes_F \bar{F}_v$ preserves a line or a dual line in $\bar{F}_v^3$,
and so $\Lie \mathbf{J}$ preserves a line or a dual line over the algebraic closure. If the Galois conjugates of this line are not contained in a plane, then
 then $\mathbf{J}$
would be a torus, contradicting the fact it contains unipotents. 
Otherwise, the Galois conjugates of the line span either a line or a plane;
 this yields a preserved line\footnote{Implicitly, we use Hilbert's theorem 90.}
in $F^3$ or in $(F^3)^*$. 

Therefore,
$$
 \PGL_3(\order)\mathbf{J}(F_S)g \subset \PGL_3(\order) \mathbf{P}(F_S) g,
$$ 
for some $F$-parabolic subgroup $\mathbf{P}$, the stabilizer of a line or a dual line. 

Moreover, we know that $\mathfrak{n} \subset \Ad(g^{-1}) \mathfrak{p}$ where $\mathfrak{p}$
is the Lie algebra of $\mathbf{P}$. 
Thus the fixed line (or dual line) for the parabolic $\Ad(g_v^{-1}) \mathbf{P}(F_v) $ is also preserved by $\mathfrak{n}$
acting on $F^3$ or $(F^3)^*$, i.e.\ belongs to the kernel of $\mathfrak{n}$. 
The kernel has dimension $2$, and so we see that the coset
$\mathbf{P}. g_v$ belongs to a one-dimensional subvariety of $\mathbf{P}(F_v) \backslash \PGL_3(F_v)$.

Applying \eqref{UB} and noting that there are only countably many $F$-parabolic
subgroups, we derive a contradiction.

{\em Case 2.} 
There exists a place $v$ so that the associated tori
$\mathbf{T}_i$ are all $F_v$-split. We may assume 
that $\disc_v(\mathscr{D}_i)$ remain bounded.  
The subgroups $H_i$ then remain in a compact set within the space of tori
in $\PGL_3(\Q_v)$. Let $H$ be any limit of the subgroups $H_i$. 
Then $\barmu_{\infty}$ is $H$-invariant and $H$ is an $F_v$-split
torus inside $\PGL_3(F_v)$. 

By Lemma \ref{balllemma},
every $H$-ergodic component of $\barmu_{\infty}$ has positive entropy with 
respect to the action of a regular element in $H$.  It follows by \cite[Theorem 2.6]{EL-ICM}
(which generalizes \cite{EKL} to the $S$-algebraic setting)
that $\barmu_{\infty}$ is $\SL_3(F_v)$-invariant.   \qed

\appendix
\section{Recollections on subconvexity}\label{subconvex-a}

In this section, we are going to briefly recall the {\em subconvexity problem}
for $L$-functions and some of the progress towards it. We refer to \cite{IS} for a more complete description of the subconvexity problem.

Let $L(\pi,s)$ denote an $L$-function attached to some arithmetic object $\pi$ (of degree $n\geq 1$),
$$L(\pi,s)=\prod_{p}L(\pi_{p},s)=\prod_{p}\prod_{i=1}^n(1-\alpha_{\pi,i}(p)p^{-s})^{-1}.$$
$L(\pi,s)$ is expected to have meromorphic continuation to $\Cc$ with (under an appropriate normalization) finitely many poles located on the lines 
$\Re s=0,1$. It satisfies a functional equation
 of the form
$$q_{\pi}^{s/2}L(\pi_{\infty},s)L(\pi,s)=\omega(\pi)q_{\pi}^{(1-s)/2}\ov{L(\pi_{\infty},1-\ov s)L(\pi,1-\ov s)}.$$
Here $$L(\pi_{\infty},s)=\prod_{i=1}^n\Gamma_{\Rr}(s+\mu_{\pi,i}),\ \Gamma_{\Rr}(s)=\pi^{-s/2}\Gamma(s/2),$$
$q_{\pi}\geq 1$ is an integer (the conductor of $\pi$) and $|\omega(\pi)|=1$. 

A subconvex bound (in the conductor aspect), is a bound of the form
\begin{equation}\label{defsubconvexbound}L(\pi,s)\ll (C_{\infty}(\pi_s))^{N}q_{\pi}^{1/4-\theta}
\end{equation}
for some absolute constants $N>0$ and $\theta>0$.
Here we denote by $C_{\infty}(\pi_s)$ the quantity
\begin{equation}\label{conductorinfinite}
C_{\infty}(\pi_s)=\prod_{i=1}^d(1+|\mu_{\pi,i}+s|),
\end{equation}

The bound \refs{defsubconvexbound} is named subconvex by comparison
 with the (easier) {\em convexity} bound  -- which may be deduced from the Phragm\'en--Lindel\"of convexity principle --
in which the exponent $1/4-\delta$ is replaced by any exponent $>1/4$.
 
In this paper the main class of $L$-functions for which we consider the subconvexity problem are the {\em Dedekind} $\zeta$-function
 of a number field $K$: the {\em Dedekind $\zeta$-function} of $K$ is a
function of a single complex variable. For $\Re(s) > 1$ it is defined by the rule
$$\zeta_K(s) = \sum_{\mathfrak{a} \subset \order_K} \Norm_{K/\Qq} (\mathfrak{a})^{-s},$$ the sum being taken
over the non-zero ideals of $\order_K$. It extends to a meromorphic function of $s$ with a simple pole at $s=1$. In that case the conductor
of $\zeta_{K}$ is the (absolute value of the) discriminant of $K$:

\begin{Hypothesis} \label{subconvex}
Let $K$ be a number field of fixed degree $n$.
There exists $\theta, N > 0$ (depending at most on $n$) such that for $\Re s=1/2$,
$$\zeta_K(s) \ll_{n}  |s|^ N \disc(K)^{1/4 - \theta}.$$
\end{Hypothesis}

By now Hypothesis \ref{subconvex} is known for a restriced class of number fields $K$: 

\begin{itemize}
\item[-] when $K$ is an abelian extensions of $\Q$ of fixed degree ($n$ say); this follows from the Kronecker-Weber theorem and from Burgess's subconvex bound for Dirichlet $L$-functions \cite{Bu}.
More generally, Hypothesis 5.1 holds if $K$ varies through the abelian extensions of a {\em fixed} number field $F$ by the work of fourth named author \cite{Ve}. 

\item[-] when $K$ is a cubic extension of $\Q$: when $K$ is not abelian, $\zeta_{K}(s)$
factors as $\zeta(s)L(\rho,s)$ where $\rho$ is a dihedral ($2$-dimensional) irreducible complex Galois representation
 of $\Gal(\overline{\Qq}/\Qq)$; more precisely $L(\rho,s)$ is the $L$-function of a cubic {\em ring class} character of the unique quadratic extension
 contained in the closure of $K$. By quadratic base change, $L(\rho,s)$ is the $L$-function of a $GL_{2,\Qq}$-automorphic form and the subconvex bound for 
 the latter class of $L$-functions follows from the works of Duke, Friedlander, Iwaniec \cite{DFI8} and Blomer, Harcos and the third author \cite{BHM}. By  the work of the third and fourth named authors \cite{MV}, this now holds when $K$ is a cubic extension of a {\em fixed} number field $F$. 

\item[-] More generally, by the above quoted works, Hypothesis \ref{subconvex} is known if $K/\Q$ is contained in a ring class field of an arbitrary quadratic extension of an arbitrary ground field $F$. 
\end{itemize}

We also need  to consider the $L$-function, $L(K,\psi,s)$, associated to a Hecke Gr\"ossencharacter of $K$ $\psi$ (in other words a character of the id\`eles of $K$, $\adele_{K}^\times/K^\times$). The conductor of $L(K,\psi,s)$ is the product of $\disc(K)$ and the ``discriminant of $\psi$.'' (Usually, the conductor of $\psi$ is defined as a certain integral ideal of $K$; the norm of this ideal is the discriminant of $\psi$.)

\begin{Hypothesis} \label{subconvexchi}
Let $K$ be a number field of degree $n$ and $\psi$ a unitary character of the id\`eles of $K$. 
Then there exists $\theta, N > 0$  (depending at most on $n$) so that for $\Re s=1/2$,
$$L(K,\psi,s) \ll_{n} C_{\infty}(\psi,s)^N (\disc(\psi)\disc(K))^{1/4-\theta},$$
where $\disc(\psi)$ denote the conductor of $\psi$.
\end{Hypothesis}
Hypothesis \ref{subconvexchi} is known in even fewer cases:
\begin{itemize}
\item[-] when $K$ is an {\em fixed} number field and $\psi$ is varying : this is a consequence of Burgess work if $K=\Qq$ and of \cite{Ve} in general.
\item[-] when $K$ is a (possibly varying) quadratic extension of the base field $F$: again this follows (by quadratic base change) from the works \cites{DFI8,BHM,MV}.
\item[-] when $K/F$ is an extension of given degree which is either, abelian, cubic or contained
 in a ring class field of an arbitrary quadratic extension of $F$ and $\psi$ {\em factors through the norm map}: that is
$\psi=\chi\circ N_{K/F}$ for some Hecke character (over $F$). In that case $L(K,\psi,s)$
 (viewed as an $L$-function ``over" $F$) equals the twist of $\zeta_{K}(s)$ by the character $\chi$
  and the subconvex bound follows from a combination of the above quoted works. 
  \end{itemize}

\def\cprime{$'$}

\end{document}